%% file: Paper.tex
\documentclass[1p,preprint]{elsarticle}

\usepackage{graphicx}
\usepackage{subfigure}
\usepackage{cancel}
\usepackage{color}
\usepackage{amsfonts}
\usepackage{amsmath,amssymb}
\usepackage{booktabs}
\usepackage{rotating}
\usepackage{hyperref}
\usepackage{lineno}

\newcommand{\bq}{\begin{equation}}
\newcommand{\eq}{\end{equation}}
\newcommand{\bqs}{\begin{equation*}}
\newcommand{\eqs}{\end{equation*}}

\newcommand{\derv}[2]{\frac{\partial #1}{\partial #2}}
\newcommand{\ra}[1]{\renewcommand{\arraystretch}{#1}}
\newcommand{\reals}{\mathbb{R}}

\newcommand{\leavethisout}[1]{}

\journal{Journal of Computational Physics}

\begin{document}

\begin{frontmatter}

\title{Riemann Solver for a Kinematic Wave Traffic Model with
  Discontinuous Flux\tnoteref{t4}} 

\author[sfu]{Jeffrey K. Wiens\corref{cor1}}
\ead{jwiens@sfu.ca}

\author[sfu]{John M. Stockie}
\ead{jstockie@sfu.ca}
\ead[url]{http://www.math.sfu.ca/~stockie}

\author[sfu]{JF Williams}
\ead{jfwillia@sfu.ca}

\address[sfu]{Department of Mathematics, Simon Fraser University, 8888
  University Drive, Burnaby, BC, Canada, V5A 1S6}
\cortext[cor1]{Corresponding author.}

\begin{abstract}
  We investigate a model for traffic flow based on the
  Lighthill--Whitham--Richards model that consists of a hyperbolic
  conservation law with a discontinuous, piecewise-linear flux.  A
  mollifier is used to smooth out the discontinuity in the flux function
  over a small distance $\epsilon\ll 1$ and then the analytical solution
  to the corresponding Riemann problem is derived in the limit as
  $\epsilon\rightarrow 0$.  For certain initial data, the Riemann
  problem can give rise to \emph{zero waves} that propagate with
  infinite speed but have zero strength.  We propose a Godunov-type
  numerical scheme that avoids the otherwise severely restrictive CFL
  constraint that would arise from waves with infinite speed by exchanging
  information between local Riemann problems and thereby incorporating
  the effects of zero waves directly into the Riemann solver.
  Numerical simulations are provided to illustrate the behaviour of zero
  waves and their impact on the solution.  The effectiveness of our
  approach is demonstrated through a careful convergence study and
  comparisons to computations using a third-order WENO scheme.
\end{abstract}

\tnotetext[t4]{This work was supported by two Discovery Grants from the
  Natural Sciences and Engineering Research Council of Canada (NSERC).}


\leavethisout{
}

\end{frontmatter}

\linenumbers



\input{Sec1}
\input{Sec2}

\input{Sec3}

\input{Sec4}

\input{Sec5}

\input{Sec6}

\bibliography{Paper}
\bibliographystyle{abbrv}

\end{document}

%% file: Sec1.tex
\section{Introduction}

In the 1950's, Lighthill and Whitham \cite{LighthillWhitham1955} and
Richards \cite{Richards1956} independently proposed the first
macroscopic traffic flow model, now commonly known as the LWR model.
Although this model has proven successful in capturing some aspects of
traffic behaviour, its limitations are well-documented and many more
sophisticated models have been proposed to capture the complex dynamics
and patterns observed in actual vehicular traffic \cite{Helbing2001}.
Despite this progress, the LWR model remains an important and
widely-used model because of its combination of simplicity and
explanatory power.

The LWR model consists of a single scalar nonlinear conservation law in
one dimension
\bq
   \rho_t + f(\rho)_x = 0, \label{eqn:ConvLaw}  
\eq
where $\rho(x,t)$ is the traffic density (cars/m),
\bqs
  f(\rho) = \rho \, v(\rho)
  \label{eqn:LWRFlux}
\eqs
is the traffic flow rate or flux (cars/sec), and $v(\rho)$ is the local
velocity (m/sec).  The most commonly used flux function is \bq f(\rho) =
u_{\text{\emph{max}}}\,\rho \left(1 -
  \frac{\rho}{\rho_{\text{\emph{max}}}} \right)
\mbox{,} \label{eqn:GreenshieldsFlux} \eq which was obtained by
Greenshields \cite{Greenshields1935} in the 1930's by fitting
experimental measurements of vehicle velocity and traffic density.
Here, $u_{max}$ is the maximum free-flow speed, while
$\rho_{\text{\emph{max}}}$ is the maximum density corresponding to
bumper-to-bumper traffic where speed drops to zero.  The LWR model
belongs to a more general class of \emph{kinematic wave traffic models}
that couple the conservation law Eq.~\eqref{eqn:ConvLaw} with a variety
of different flux functions.

Extensive studies of the empirical correlation between flow rate and
density have been performed in the traffic flow literature.  This correlation
is commonly referred to as the \emph{fundamental diagram} and is
represented graphically by a plot of flux~$f$ versus density~$\rho$ such
as that shown in Fig.~\ref{fig:BothFluxes}.  A striking feature of many
experimental results is the presence of an apparent discontinuity that
separates the free flow (low density) and congested (high density)
states, something that has been discussed by many authors, including
\cite{Ceder1976,Easa1982,Edie1961,Kerner2004}.  In particular, Koshi et
al.~\cite{KoshiEtAl1983} characterize flux data such as that shown in
Fig.~\ref{fig:BothFluxes} as having a \emph{reverse lambda} shape in
which the discontinuity appears at some peak value of the flux.

This behavior is also referred to as the \emph{two-capacity} or
\emph{dual-mode phenomenon} \cite{Bank1991a,Bank1991b} and has led to
the development of a diverse range of mathematical models.  Zhang and
Kim \cite{ZhangKim2005} incorporated the capacity drop into a
microscopic car-following model that generates fundamental diagrams with
the characteristic reverse-lambda shape.  Wong and Wong
\cite{WongWong2002} performed simulations using a multi-class LWR model
from which they also observed a discontinuous flux-density relationship.
Colombo \cite{Colombo2002} and Goatin \cite{Goatin2006} developed a
macroscopic model that couples an LWR equation for density in the free
flow state, along with a 2$\times$2 system of conservation laws for
density and momentum in the congested state; the phase transition
between these two states is a free boundary that is governed by the
Rankine-Hugoniot conditions.  Lu et al.~\cite{Lu2009} incorporated a
discontinuous (piecewise quadratic) flux directly into an LWR model, and
then solved the corresponding Riemann problem analytically by
constructing the convex hull for a regularized continuous flux function
that consists of two quadratic pieces joined over a narrow region by a
linear connecting piece.

There remains some disagreement in the literature regarding the
existence of discontinuities in the traffic flux, with some researchers
(e.g., Hall \cite{HallEtAl1986}) arguing that the apparent gaps are due
simply to missing data and can be accounted for by providing additional
information about traffic behaviour at specific locations.  Indeed,
Persaud and Hall \cite{Persaud-Hall-1989} and Wu \cite{Wu2002} contend
that the discontinuous fundamental diagram should be viewed instead as
the 2D projection of a higher dimensional smooth surface.  

We will nonetheless make the assumption in this paper that the
fundamental diagram is discontinuous.  Our aim here is not to argue the
validity of this assumption in the context of traffic flow, since that
point has already been discussed extensively by
\cite{Lu2009,WongWong2002,ZhangKim2005}, among others.  Instead our
objective is to study the effect that such a flux discontinuity has on
the analytical solution of a 1D hyperbolic conservation law, as well as
to develop an accurate and efficient numerical algorithm to simulate
such problems.

A related class of conservation laws, in which the flux $f(\rho,x)$ is a
discontinuous function of the spatial variable $x$, has been thoroughly
studied in recent years (see \cite{BurgerEtAl2008,BurgerKarlsen2008} and
the references therein).  Considerably less attention has been paid to
the situation where the flux function has a discontinuity in $\rho$.
Gimse~\cite{Gimse1993} solved the Riemann problem for a piecewise linear
flux function with a single jump discontinuity in $\rho$ by generalizing
the method of convex hull construction \cite[Ch.~16]{LeVequeRedBook}.
In particular, Gimse identified the existence of \emph{zero shocks},
which are discontinuities in the solution that carry no information and
have infinite speed of propagation.  We note that more recently,
Armbruster et al.~\cite{Armbruster2011} observed \emph{zero rarefaction
  waves} with infinite speed of propagation in their study of supply
chains with finite buffers (although they did not refer to them using
this terminology).

Gimse's results were improved on by Dias and
Figueira~\cite{DiasFigueira2005}, who used a mollifier function
$\eta_\epsilon$ to smooth out discontinuities in the flux function over
an interval of width $0<\epsilon\ll 1$ before constructing the convex
hull using standard techniques.  Solutions to the mollified problem were
proven to converge to solutions of the original problem in the limit as
$\epsilon\rightarrow 0$~\cite{DiasFigueira2005}.  Dias and Figueira's
framework has also been applied to problems involving fluid phase
transitions~\cite{DiasFigueira2004,DiasFigueiraRodrigues2005} and 
viscoelasticity~\cite{DiasFigueira2005b}.

In this paper, we apply Dias and Figueira's mollification approach to
solving a conservation law with a piecewise linear flux function
$f(\rho)$ in which there is a single discontinuity at $\rho=\rho_m$ (see
Fig.~\ref{fig:BothFluxes}).  The model equations and their relevance in
the context of traffic flow are discussed in
Section~\ref{sec:Equations}.  We introduce a mollified flux function
$f_\epsilon(\rho)$ in Section \ref{sec:RiemannProblem} and verify its
convexity, which then permits us to derive the analytical solution to
the Riemann problem using the method of convex hull construction.

In Section \ref{sec:ZeroWaves}, we consider the special case where
either of the two constant initial states in the Riemann problem equals
$\rho_m$, the density at the discontinuity point.  This is precisely the
case when a rarefaction wave of strength $O(\epsilon)$ and speed
$O(1/\epsilon)$ arises, which approaches a zero rarefaction in the limit
of vanishing $\epsilon$.  There are two issues that need to be addressed
regarding these zero waves.  First, we consider the convergence of the
mollified solution to that of the original problem, since Dias and
Figueira's convergence results \cite{DiasFigueira2005} do not consider
(nor easily extend to) the case when the left or right initial states in
the Riemann problem are identical to $\rho_m$.  Secondly, we discuss the
physical relevance of an infinite speed of propagation in the context of
traffic flow.

The remainder of the paper is focused on constructing a Riemann solver
that forms the basis for a high resolution finite volume scheme of
Godunov type.  Because zero waves travel at infinite speed, the usual
CFL restriction suggests that choosing a stable time step might not be
possible.  Some authors have avoided this difficulty by using an
implicit time discretization~\cite{MartinVovelle2008}, but this approach
introduces added expense and complication in the numerical algorithm.
Another approach employed in \cite{Armbruster2011} is to replace the
discontinuous flux by a regularized (continuous) function which joins
the discontinuous pieces by a linear connection over an interval of
width $\epsilon\ll 1$, after which standard numerical schemes may be
applied; however, this approach requires a small $\epsilon$ to achieve 
reasonable accuracy resulting in a severe time step
restriction.

We use an alternate approach that eliminates the severe CFL constraint
by incorporating the effect of zero waves directly into the local
Riemann solver.  In the process, we find it necessary to construct
solutions to a subsidiary problem that we refer to as the \emph{double
  Riemann problem}, which introduces an additional intermediate state
corresponding to the discontinuity value $\rho=\rho_m$.  A similar
approach was used by Gimse~\cite{Gimse1993} who constructed a
first-order variant of Godunov's method, although he omitted to
perform any computations using his proposed method.  We improve upon
Gimse's work in three ways: first, we solve the double Riemann problem
within Dias and Figueira's mollification framework; second, we implement
a high resolution variant of Godunov's scheme to increase the spatial
accuracy; and third, we provide extensive numerical computations and a
careful convergence study to demonstrate the effectiveness of our
approach.

%% file: Sec2.tex
\section{Mathematical Model}
\label{sec:Equations}

\begin{figure}[btp]
  \begin{center}
    \subfigure[]{\label{fig:discFlux} 
      \includegraphics[width=0.45\textwidth]{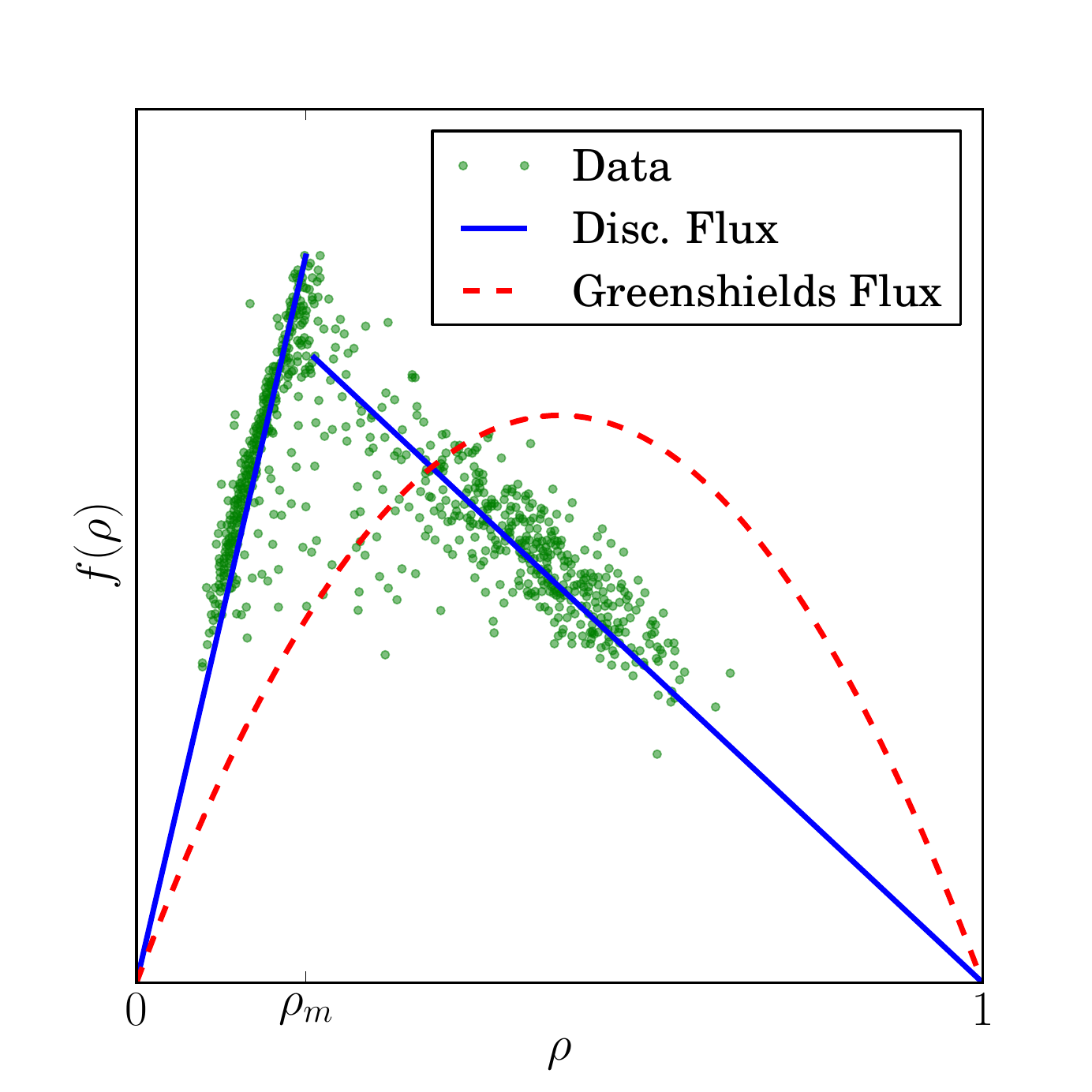}}
    \subfigure[]{\label{fig:smoothFlux}
      \includegraphics[width=0.45\textwidth]{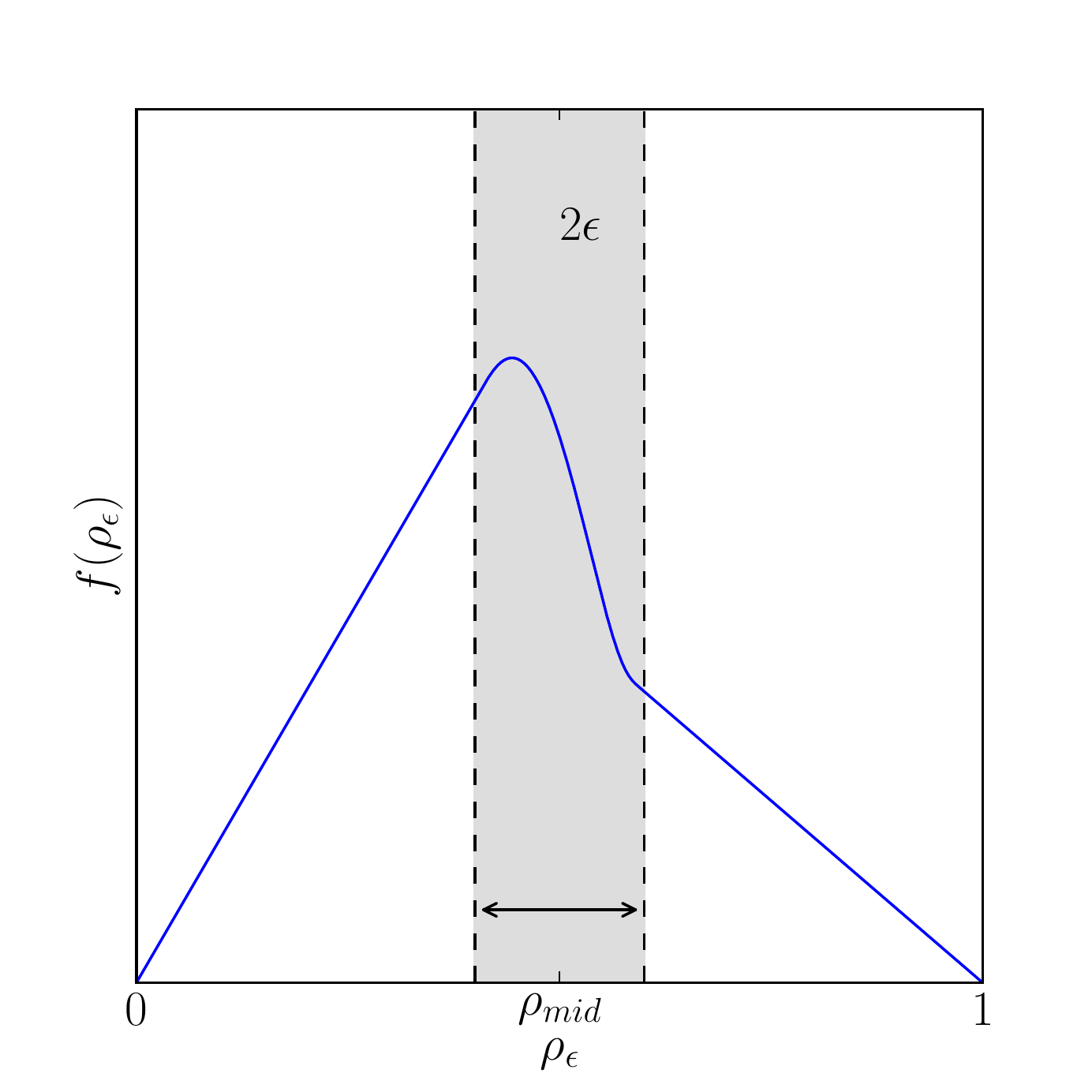}}
    \caption{In \subref{fig:discFlux}, the discontinuous ``reverse
      lambda'' flux function Eq.~\eqref{eqn:DiscFlux} is overlaid with
      empirical data extracted from \cite[Fig.~7]{HallEtAl1986}
      (reproduced with permission of Elsevier B.V.), along with the
      quadratic Greenshields flux~\eqref{eqn:GreenshieldsFlux}.  The
      mollified flux from Eq.~\eqref{eqn:MollifiedFlux} is depicted in
      \subref{fig:smoothFlux}.}
  \label{fig:BothFluxes}
  \end{center}
\end{figure}

We are concerned in this paper with the scalar conservation law
\bq 
  \rho_t + f(\rho)_x = 0,
  \label{eqn:DiscConsLaw} 
\eq
having a discontinuous flux function
\bq 
  f(\rho) =  \left\{
    \begin{array}{ll}
      \displaystyle g_{\mathit{f}}(\rho), & \text{if}~0 \leqslant \rho < \rho_{m}, \\
      \displaystyle g_{\mathit{c}}(\rho), & \text{if}~\rho_{m} \leqslant
      \rho \leqslant 1,
    \end{array}
  \right.
  \label{eqn:DiscFlux}
\eq
that is depicted in Fig.~\ref{fig:BothFluxes}\subref{fig:discFlux}.  The
vehicle density $\rho(x,t)$ is normalized so that $0 \leqslant \rho
\leqslant 1$, and $\rho_m$ is the point of discontinuity in the flux
$f(\rho)$.  We restrict the flux to be a piecewise linear
function in which the free flow branch has
\bqs
   g_{\mathit{f}}(\rho) = \rho,
   \label{eqn:DiscFlux-gf}
\eqs
and the congested flow branch has
\bqs
  g_{\mathit{c}}(\rho) = \gamma( 1 - \rho).
   \label{eqn:DiscFlux-gc}
\eqs
Experimental data suggests that $g_{\mathit{f}}(\rho_m) >
g_{\mathit{c}}(\rho_m)$, and so we impose the constraint
\bq 
   0 < \gamma < \frac{\rho_m}{1-\rho_m} \label{eqn:FluxConstraint} \mbox{.}
\eq

We utilize the mollifier approach of Dias and Figueira
\cite{DiasFigueira2005} in order to approximate the original equation
\eqref{eqn:DiscConsLaw} by
\bq
  \derv{\rho_\epsilon}{t} + \derv{f_\epsilon(\rho_\epsilon)}{x} =
  0, \label{eqn:MolConsLaw} 
\eq
where the mollified flux (pictured in
Fig.~\ref{fig:BothFluxes}\subref{fig:smoothFlux}) is 
\bq 
  f_\epsilon(\rho_\epsilon) = \rho_\epsilon + ( \gamma -
  (\gamma+1)\rho_\epsilon) \int_{\rho_{m}-\epsilon}^{\rho_\epsilon}
  \eta_\epsilon(s-\rho_{m})\, ds \mbox{,} 
  \label{eqn:MollifiedFlux}
\eq
with $0< \epsilon \ll 1$.  The mollifier function is given by
$\eta_\epsilon(s) = \frac{1}{\epsilon} \eta(s/\epsilon)$, where $\eta(s)$
is a \emph{canonical mollifier} that satisfies the following conditions:
\begin{enumerate}
  \renewcommand{\theenumi}{\roman{enumi}}
  \renewcommand{\labelenumi}{(\theenumi)}
  \renewcommand{\itemsep}{3pt}
\item $\eta \geqslant 0$; 
\item $\eta \in C^\infty(\reals)$ and is compactly supported on
  $[-1,1]$;
\item $\eta(-s)=\eta(s)$ for all $s\in \reals$; and
\item \label{item:iv} $\int_{-\infty}^\infty \eta(s)\, ds = 1$.
\end{enumerate}
The results in Dias and Figueira~\cite{DiasFigueira2005}
guarantee that any mollifier satisfying the above criteria converges to
the same unique solution in the limit $\epsilon \rightarrow 0$.  We use
the following mollifier
\bq
   \eta(s) = \left\{
	\begin{array}{ll}
	  C \exp \left(\frac{1}{s^2-1}\right), & \text{if}~|s| < 1, \\
	  0, &  \text{if}~|s| \geqslant 1,
	\end{array}
    \right. 
    \label{eqn:BumpMollifier}
\eq
where $C\approx 2.2522836\dots$\ is a constant determined numerically so
that condition (\ref{item:iv}) holds; this choice is made for reasons of
analytical convenience since the derivative $\eta^\prime$ can be written
in terms of $\eta$.

Because the mollified flux function is smooth, the conservation
law~\eqref{eqn:MolConsLaw} may now be solved using standard techniques.
We note that in the context of traffic flow, a potential problem arises
when applying the usual Ole\u{\i}nik entropy condition
\cite{Oleinik1964} as the selection principle to enforce uniqueness.
Although Ole\u{\i}nik's entropy condition does yield the
physically-correct weak solution in the context of fluid flow
applications, it does not always do so for kinematic wave models of
traffic flow (see LeVeque~\cite{LeVeque2001}, for example).  In
particular, applying Ole\u{\i}nik's entropy condition can lead to a
solution that is not anisotropic \cite{Daganzo1995}, corresponding to
the non-physical situation where drivers react to vehicles both in front
\emph{and behind}.
  
Zhang~\cite{Zhang2003} suggests two additional criteria on the flux
function to guarantee anisotropic flows in kinematic wave traffic
models.  First, the characteristic velocity should be smaller than the
vehicle speed $v(\rho) = f(\rho)/\rho$.  That is, we require
\bqs
  \frac{f(\rho)}{\rho} \geq f^\prime(\rho),~~ \forall ~ 0 \leq \rho
  \leq 1. 
\eqs
Secondly, all elementary waves must travel more slowly than the vehicles
carrying them, or
\bqs
  \min \left( \frac{f(\rho_l)}{\rho_l},~ \frac{f(\rho_r)}{\rho_r}
  \right) \geq  \frac{f(\rho_r) - f(\rho_l)}{\rho_r - \rho_l} ,
\eqs
for all $0 \leq \rho_l,\rho_r \leq 1$. As shown in \cite{Wiens2011}, our
flux function \eqref{eqn:DiscFlux} satisfies both of these conditions
and therefore it is reasonable to apply the Ole\u{\i}nik entropy
condition as the selection criterion for our traffic flow model.

%% file: Sec3.tex
\section{Exact Solution of the Riemann Problem with Mollified Flux}
\label{sec:RiemannProblem}

We next construct and analyze the solution $\rho_\epsilon(x,t)$ of the
mollified Riemann problem, which consists of the conservation law
\eqref{eqn:MolConsLaw} along with flux \eqref{eqn:MollifiedFlux} and
piecewise constant initial conditions
\bqs
    \rho_\epsilon(x, 0) = \left\{
	\begin{array}{ll}
	  \rho_l, & \text{if}~ x<0, \\
	  \rho_r, & \text{if}~ x\geqslant 0.
	\end{array}
    \right. 
\eqs
This problem can be solved using the method of convex hull construction
\cite[Ch.~16]{LeVequeRedBook} which requires knowledge of the inflection
points of the mollified flux $f_\epsilon(\rho_\epsilon)$.  Using
Eq.~\eqref{eqn:BumpMollifier}, the first and second
derivatives of the flux are
\begin{align}
   f_\epsilon^\prime(\rho_\epsilon) &= 1 +
   \left[\gamma-(\gamma+1)\rho_\epsilon\right] \,
   \eta_\epsilon(\rho_\epsilon-\rho_{m}) 
   - (\gamma +1)\int_{\rho_{m}-\epsilon}^{\rho_\epsilon}
   \eta_\epsilon(s-\rho_{m})\, ds , 
  \label{eqn:FirstDerivativeFlux}
  \\
  f_\epsilon^{\prime\prime}(\rho_\epsilon) &=  \kappa(\rho_\epsilon) \left[
    \left(\frac{\gamma}{\gamma+1}-\rho_\epsilon\right)(\rho_{m}-\rho_\epsilon)
    - \epsilon^2
    \left[\left(\frac{\rho_\epsilon-\rho_{m}}{\epsilon}\right)^2-1\right]^2\right] 
  \label{eqn:SecondDerivativeFlux} ,
\end{align}
where
\bqs
  \kappa(\rho_\epsilon) = \frac{2
    \eta_\epsilon(\rho_\epsilon-\rho_{m})}{\epsilon^2(\gamma+1)
    [(\frac{\rho_\epsilon-\rho_{m}}{\epsilon})^2-1]^2} 
  \label{eqn:Kappa_SecondDerivativeFlux} .
\eqs
Since the mollifier has compact support, we know that there can be no
inflection points outside the smoothing region of width $2\epsilon$;
that is, $f_\epsilon^{\prime\prime}(\rho_\epsilon) \equiv 0$ when
$|\rho_\epsilon - \rho_{m} | \geqslant \epsilon$.  Also, since
$\kappa(\rho_\epsilon)>0$, the convexity of $f_\epsilon$ is determined
solely by the sign of the quantity
\bq 
  {P}(\rho_\epsilon) =
  \left(\frac{\gamma}{\gamma+1}-\rho_\epsilon\right)(\rho_{m}-\rho_\epsilon)
  - \epsilon^2
  \left[\left(\frac{\rho_\epsilon-\rho_{m}}{\epsilon}\right)^2-1\right]^2 ,
  \label{eqn:ugly_polynomial}
\eq
which we analyze next.  

\begin{figure}[btp]
  \begin{center}
    \subfigure[]{\includegraphics[width=0.45\textwidth]{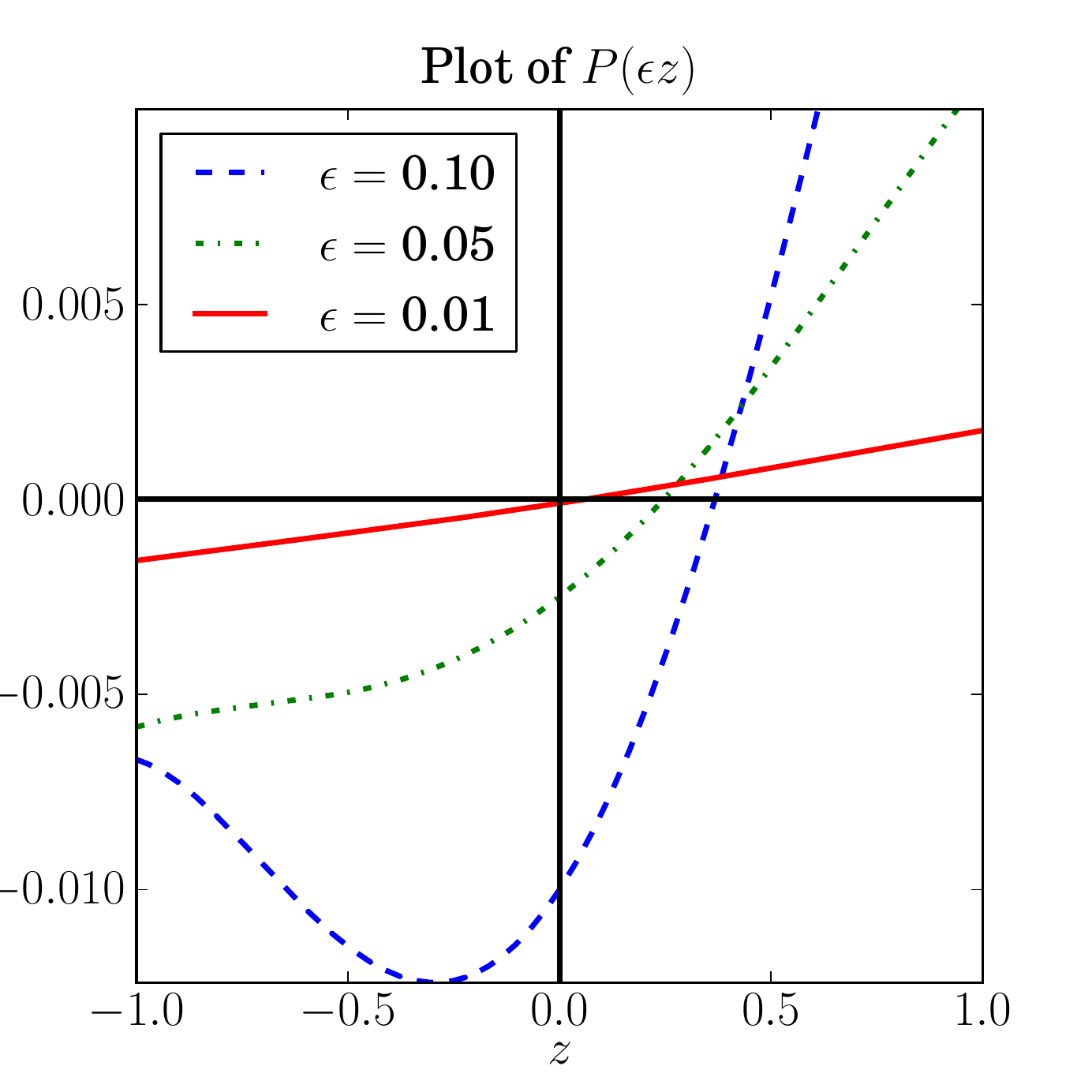}}
    \subfigure[]{\includegraphics[width=0.45\textwidth]{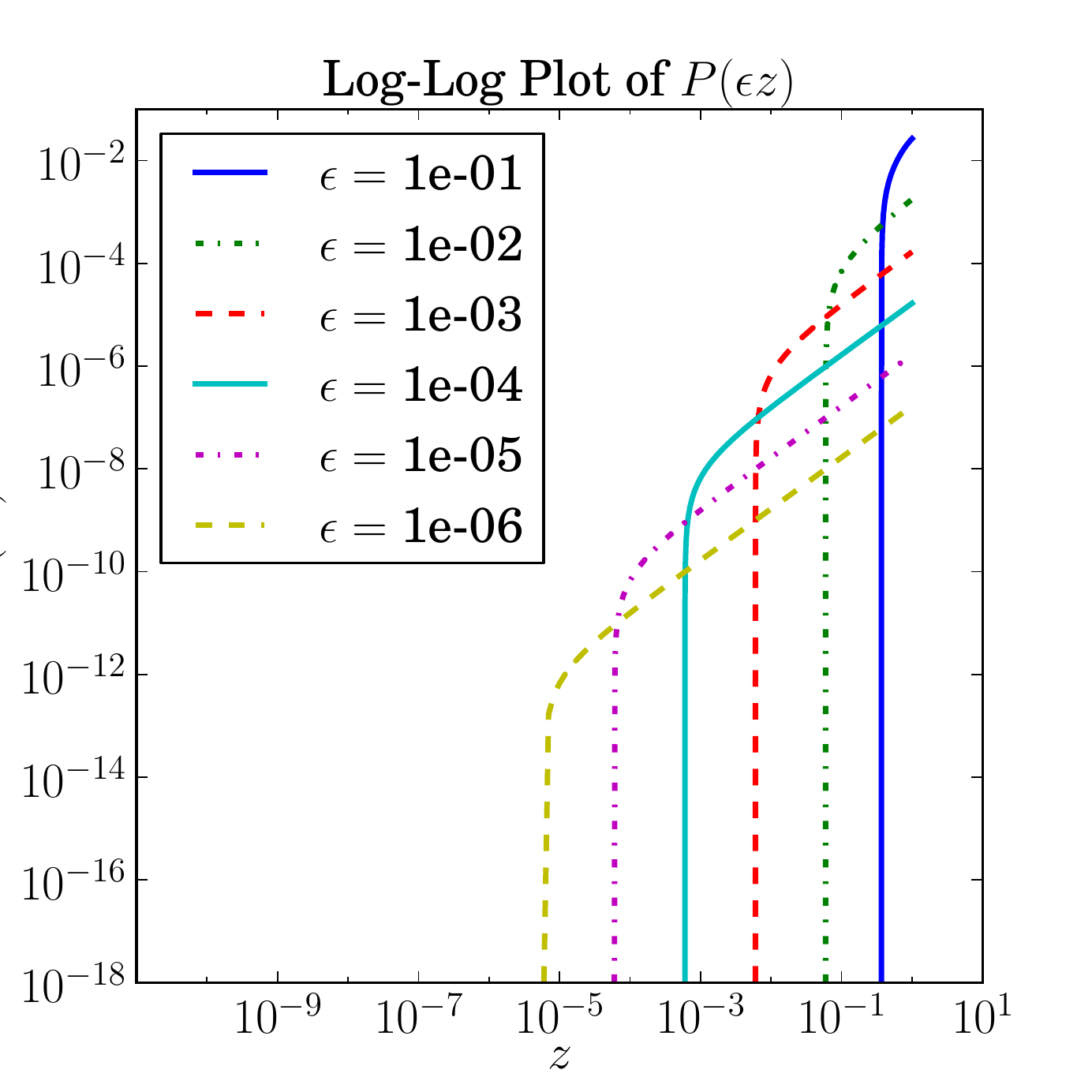}}
    \caption{Plots of the function $P(\epsilon z)$ and its roots for
      different $\epsilon$ when $\gamma = 0.5$ and $\rho_{m}=0.5$.}
  \label{fig:BothPolyPlots}
  \end{center}
\end{figure}

Since Eq.~\eqref{eqn:ugly_polynomial} is a quartic polynomial in
$\rho_\epsilon$, analytic expressions are available for the roots;
however, these are too complicated for our purposes.  Instead, we take
advantage of the scaling properties of the polynomial 
to simplify $P$ and rewrite Eq.~\eqref{eqn:ugly_polynomial} as
\bq
   P(y) = -\frac{1}{\epsilon^2} y^4 + 3 y^2 + M y - \epsilon^2 ,
  \label{eqn:nice_polynomial}
\eq
where $y = \rho_\epsilon - \rho_{m}\in[-\epsilon, \epsilon]$ and $M =
\rho_{m} -\gamma/(\gamma+1)$.  Clearly $|y|=|\rho_\epsilon - \rho_{m}| <
\epsilon$, and so it is natural to suppose that $y =
O(\epsilon^\lambda)$ with $\lambda \geq 1$, after which the various
terms in Eq.~\eqref{eqn:nice_polynomial} have the scalings indicated
below: 
\bq
   P(y) = 
	\underbrace{ -\frac{1}{\epsilon^2} y^4 \vphantom{-\frac{1}{\epsilon^2} y^4 }}_{O(\epsilon^{4\lambda-2})}
      + \underbrace{3 y^2 \vphantom{-\frac{1}{\epsilon^2} y^4}}_{O(\epsilon^{2\lambda})} 
      + \underbrace{M y \vphantom{-\frac{1}{\epsilon^2} y^4}}_{O(\epsilon^\lambda)}
      - \underbrace{\epsilon^2 \vphantom{-\frac{1}{\epsilon^2} y^4}}_{O(\epsilon^2)} \mbox{.}
  \label{eqn:order_polynomial}
\eq
Since $M$ is a constant that is independent of $\epsilon$, we can take
$M = O(1)$ as $\epsilon\rightarrow 0$.  Guided by the scalings in
Eq.~\eqref{eqn:order_polynomial}, the dominant terms in $P(y)$ are the
last two terms having orders $O(\epsilon^\lambda)$ and $O(\epsilon^2)$. 
When $\epsilon$ is sufficiently small,\ 
we may therefore neglect the remaining terms and determine the convexity 
of $f_\epsilon$ based on the sign of the simpler linear polynomial
\bqs
   P_1(y) = M y - \epsilon^2 \mbox{,}
  \label{eqn:firstorder_roots}
\eqs
which has a single root at $y=\epsilon^2/M$ corresponding to
$\rho_\epsilon = \rho_m + \epsilon^2/M$.  For values of $y$ close enough
to $\epsilon^2/M$, higher order terms in the polynomial $P(y)$ become
significant and could potentially introduce additional roots; however,
by continuing this method of dominant balance, we find that $P(y)$
maintains the single root at higher orders as well, which we demonstrate
numerically using the plots of $P$ summarized in
Fig.~\ref{fig:BothPolyPlots}.  Based on this argument and the
observation that inequality \eqref{eqn:FluxConstraint} requires $M>0$,
we can conclude that
\begin{align*}
   & f_\epsilon^{\prime\prime}(\rho_\epsilon) \geqslant 0 
   \quad \mbox{when} \quad
   \rho_\epsilon \in [\rho_{m},1]\mbox{,}  \\ 
   \text{and} \qquad 
   & f_\epsilon^{\prime\prime}(\rho_\epsilon) \leqslant 0 
   \quad \mbox{when} \quad
   \rho_\epsilon \in [0,\rho_{m}]\mbox{.}  
\end{align*}
Therefore, the mollified flux function has a single inflection point at
$\rho_\epsilon = \rho_m + O(\epsilon^2)$ as $\epsilon \rightarrow 0$,
where the slope $f_\epsilon^{\prime}(\rho_{m}) \rightarrow - \infty$.
Since the flux derivative $f^\prime(\rho)$ corresponds to the elementary
wave speed in the Riemann problem, this same point $\rho_m$ is also the
source of the infinite speed of propagation which will be analyzed in
more detail in Section~\ref{sec:ZeroWaves}.

Using this information, we can now construct the convex hull of the flux
function $f_\epsilon(\rho_\epsilon)$ which is then used to solve the
Riemann problem. There are three non-trivial cases to consider,
depending on the left and right initial states, $\rho_l$ and $\rho_r$.
Two of these cases (which we call A and B) lead to the emergence of a
new constant intermediate state with density $\rho_m$.  This ``plateau''
is a characteristic feature of solutions to our LWR model with
discontinuous flux.
\\
\\
\noindent{\bf Case~A: $\boldsymbol{\rho_r < \rho_{m} < \rho_l}$}.

\begin{figure}[btp]
  \begin{center}
    \subfigure[]{\label{fig:Case1_ConvexHull}\includegraphics[width=0.45\textwidth]{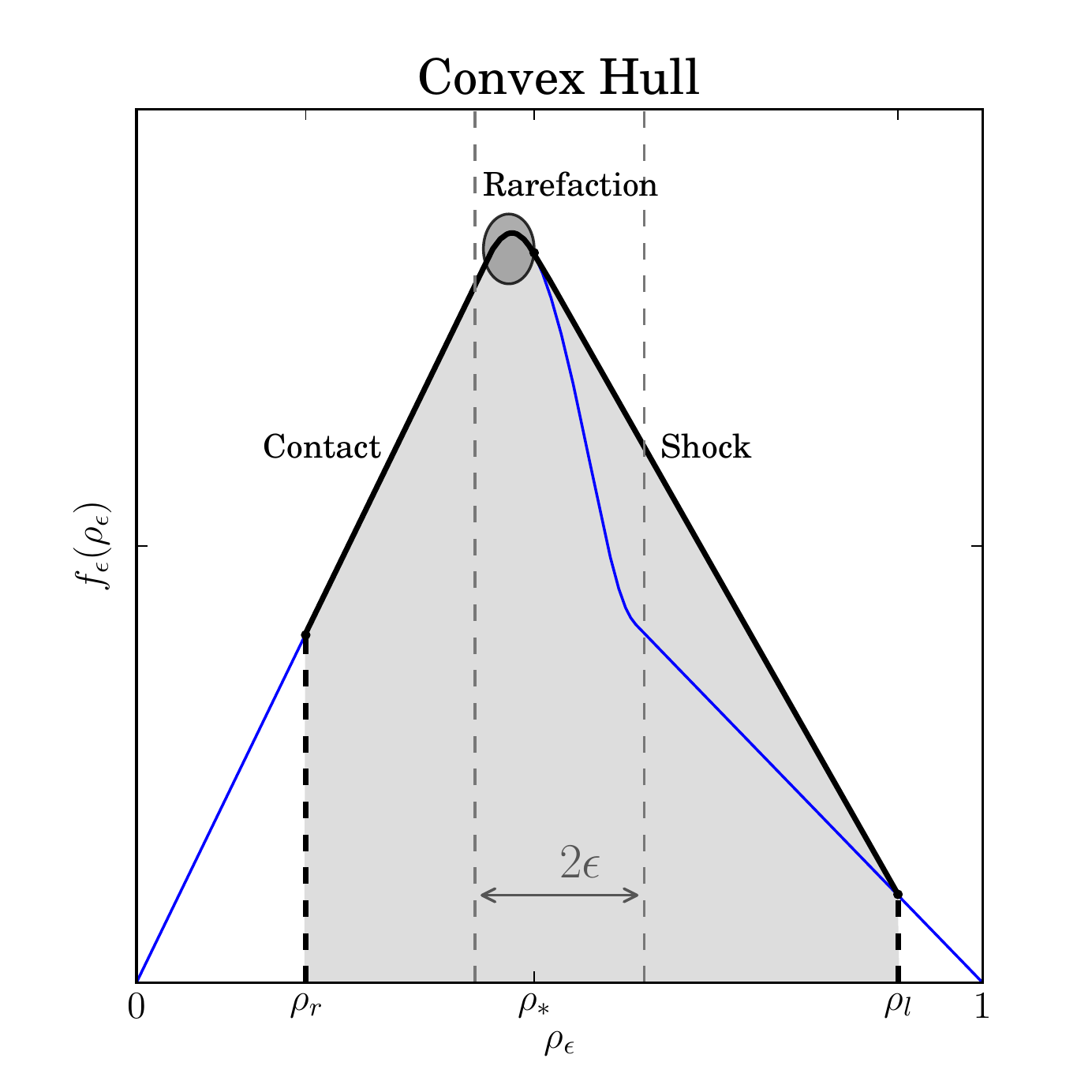}}
    \subfigure[]{\label{fig:Case1_ShockProfile}\includegraphics[width=0.45\textwidth]{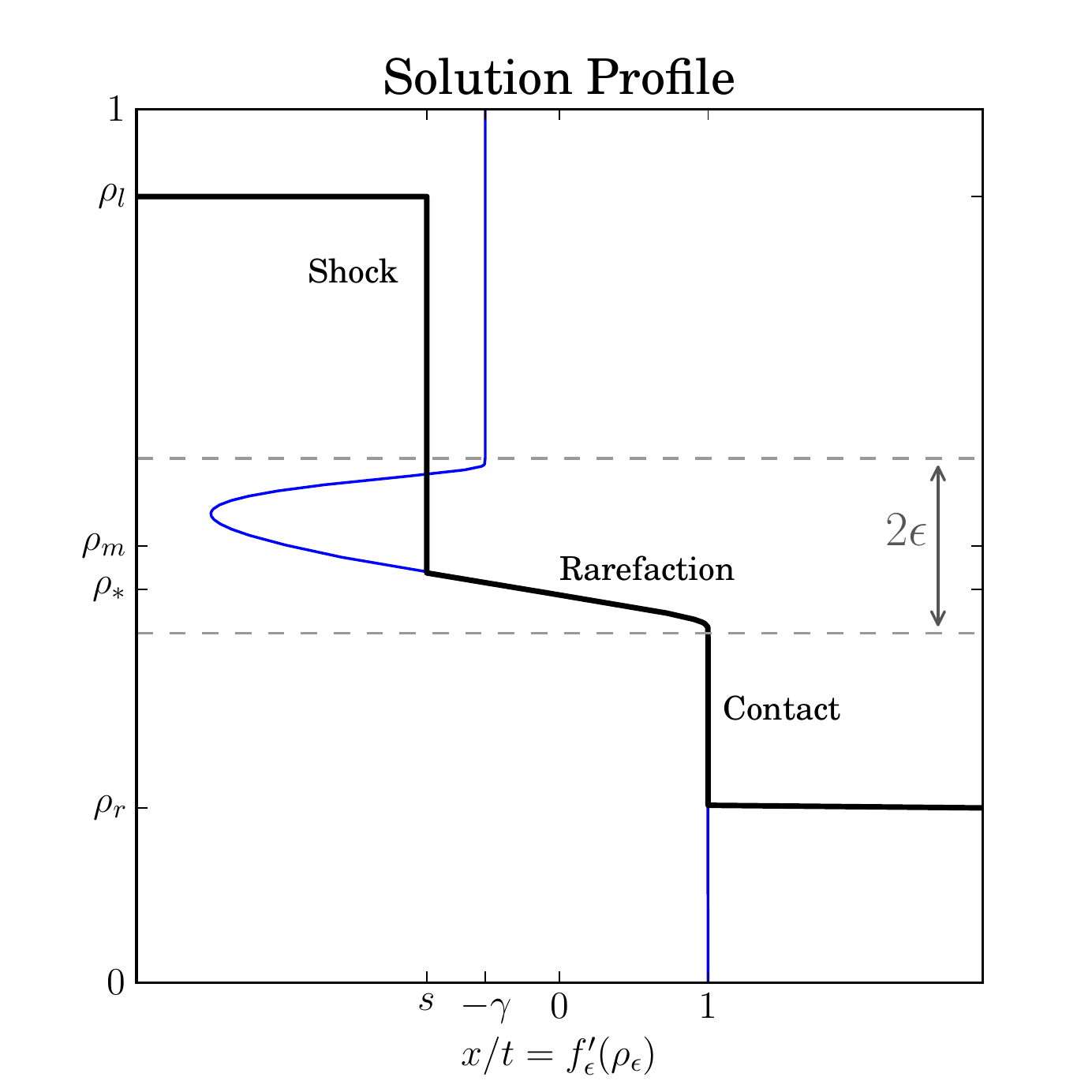}}
    \caption{The convex hull construction for Case~A.
      \subref{fig:Case1_ConvexHull} The thin (blue) line is the
      mollified flux function $f_\epsilon(\rho_\epsilon)$ and the thick
      (black) line is the convex hull.  The shaded oval region
      highlights the rarefaction wave. \subref{fig:Case1_ShockProfile}
      The solution $\rho_\epsilon(x/t)$ is shown as a thick (black)
      line, along with the corresponding flux derivative
      $f^\prime(\rho)$ as a thin (blue) line.  Parameter values are
      $\rho_l=0.9$, $\rho_r=0.2$, $\rho_{m}=0.5$, $\epsilon = 0.1$ and
      $\gamma = 0.5$.}
  \label{fig:Case1_RiemannSolution}
\end{center}
\end{figure}

Here we construct the smallest convex hull of the set
$\{(\rho_\epsilon,y): \rho_r < \rho_\epsilon < \rho_l \mbox{ and } y
\leqslant f_\epsilon(\rho_\epsilon)\}$, which as shown in
Fig.~\ref{fig:Case1_RiemannSolution}\subref{fig:Case1_ConvexHull} must
consist of three pieces. The first piece corresponds to a contact line
that follows $f_\epsilon(\rho_\epsilon)$ on the left up to the point
$\rho_* \in (\rho_{m}-\epsilon, \rho_{m})$ for which the shock and
characteristic speeds are equal; that is,
\bqs
   s = \frac{f_\epsilon(\rho_l) - f_\epsilon(\rho_*)}{\rho_l-\rho_*} =
   f^\prime_\epsilon(\rho_*) \mbox{.} 
\eqs
The third piece of the convex hull corresponds to a shock that connects
the states $\rho_*$ and $\rho_l$.  The middle piece, in between the
contact line and shock, gives rise to a rarefaction wave that follows
the curved portion of the flux in the neighbourhood of $\rho_m$.  Based
on this convex hull, we can then construct the solution profile shown in
Fig.~\ref{fig:Case1_RiemannSolution}\subref{fig:Case1_ShockProfile}. Since
the rarefaction wave consists of density values bounded between
$\rho_{m}-\epsilon$ and $\rho_{m}+\epsilon$, this wave flattens out and
degenerates to a constant intermediate state $\rho_{m}$ in the limit as
$\epsilon\rightarrow 0$.

To summarize, in the limit as $\epsilon \rightarrow 0$ the solution to
the Riemann problem when $\rho_r < \rho_{m} < \rho_l$ is
\bq
    \rho(x, t) = \left\{
	\begin{array}{ll}
	  \rho_l,  & \text{if}~x< s t,\\
	  \rho_{m},& \text{if}~s t \leqslant x \leqslant t, \\
	  \rho_r,  & \text{if}~x > t,
	\end{array}
    \right.
    \label{eqn:Case1_Solution}
\eq
consisting of a 1-shock moving to the left with speed
\bq
   s = \frac{f(\rho_l) - g_f(\rho_{m})}{\rho_l-\rho_{m}} < 0, 
  \label{eqn:Case1_ShockSpeed}
\eq
and a 2-contact moving to the right with speed 1.
\\
\\
\noindent{\bf Case~B: $\boldsymbol{\gamma/(\gamma+1) < \rho_l < \rho_{m}
    < \rho_r}$.} 

\begin{figure}[btp]
  \begin{center}
    \subfigure[]{\label{fig:Case2_ConvexHull}\includegraphics[width=0.45\textwidth]{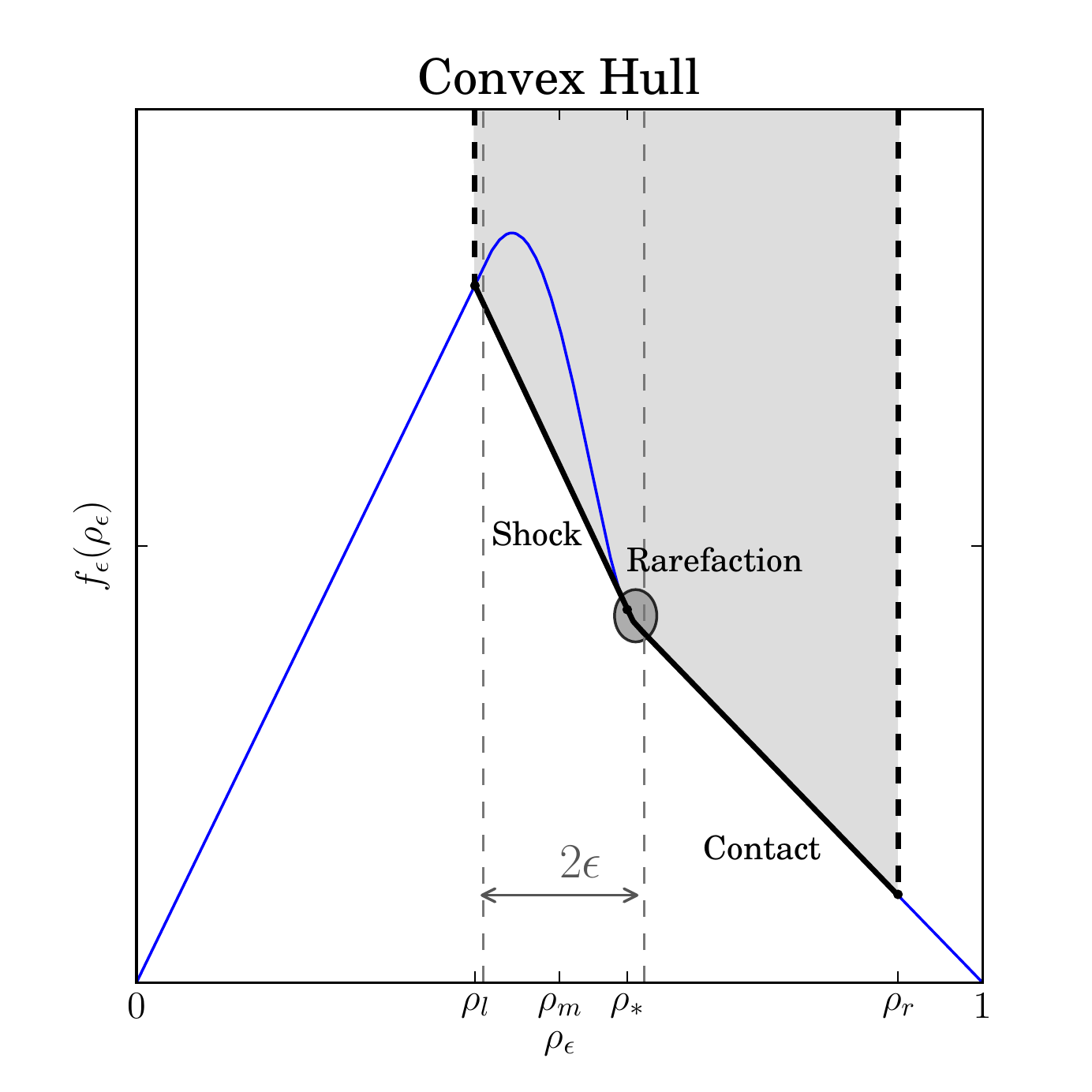}}
    \subfigure[]{\label{fig:Case2_ShockProfile}\includegraphics[width=0.45\textwidth]{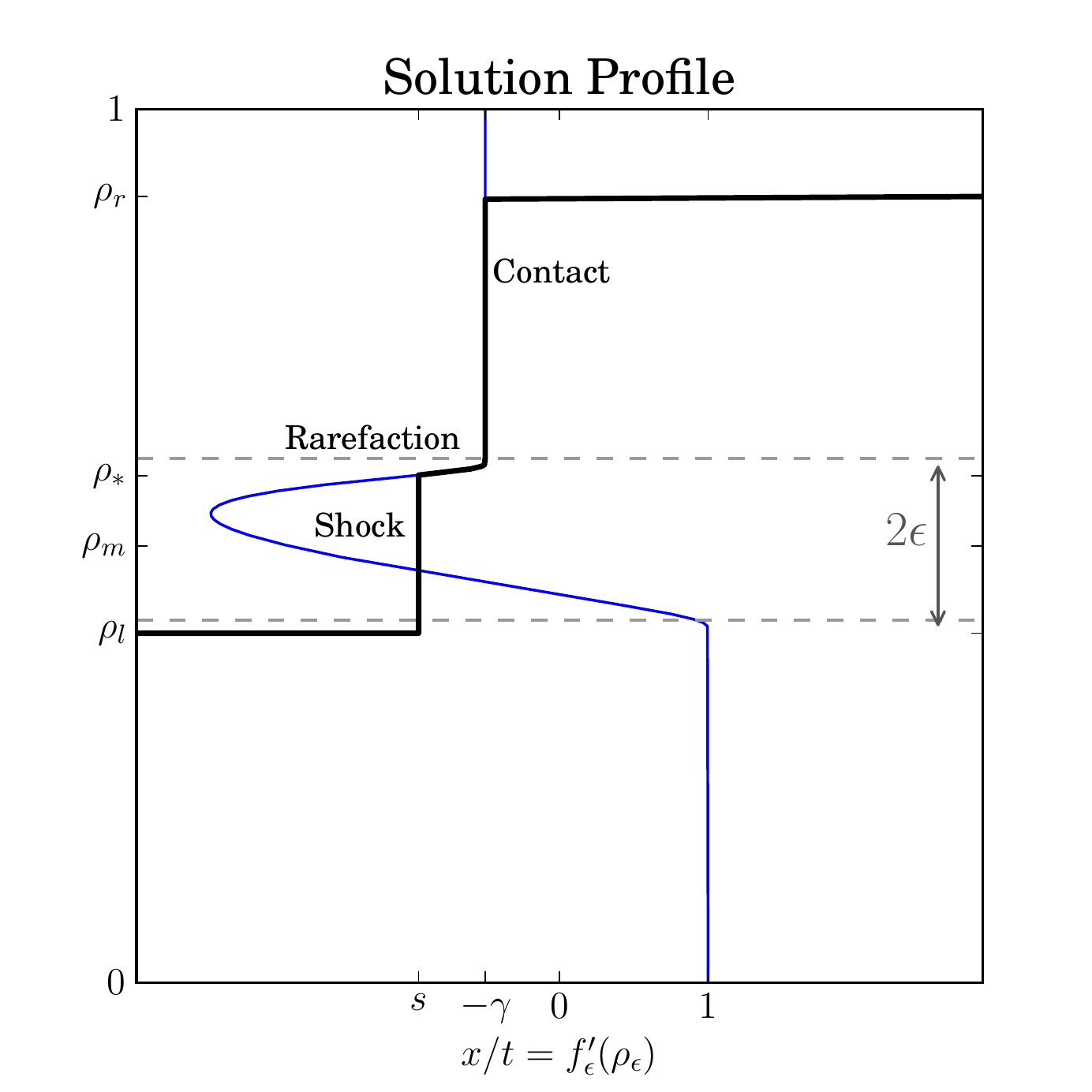}}
    \caption{The convex hull construction for Case~B.
      \subref{fig:Case2_ConvexHull} The thin (blue) line is the
      mollified flux function $f_\epsilon(\rho_\epsilon)$ and the thick
      (black) line is the convex hull.  The shaded oval region
      highlights the rarefaction wave. \subref{fig:Case2_ShockProfile}
      The solution $\rho_\epsilon(x/t)$ is shown as a thick (black)
      line, along with the corresponding flux derivative
      $f^\prime(\rho)$ as a thin (blue) line.  Parameter values are
      $\rho_l=0.4$, $\rho_r=0.9$, $\epsilon = 0.1$, $\gamma = 0.5$,
      $\rho_{m}=0.5$.}
    \label{fig:Case2_RiemannSolution}
 \end{center}
\end{figure}

In contrast with the previous case, we now have $\rho_l<\rho_r$ and so
the convex hull $\{(\rho_\epsilon,y) : \rho_l < \rho_\epsilon < \rho_r
\mbox{ and } y \geqslant f_\epsilon(\rho_\epsilon)\}$ lies above the
flux-density curve as shown in
Fig.~\ref{fig:Case2_RiemannSolution}\subref{fig:Case2_ConvexHull}.  The
first piece of the convex hull corresponds to a shock wave connecting
the states $\rho_l$ and $\rho_* \in (\rho_{m}-\epsilon, \rho_{m})$,
while the third piece follows the portion of the flux function in the
region $[\rho_*,\rho_r]$.  As before, the state $\rho_*$ is chosen so
that the shock and characteristic speeds are equal:
\bqs
   s = \frac{f_\epsilon(\rho_*) - f_\epsilon(\rho_l)}{\rho_*-\rho_l} =
   f^\prime_\epsilon(\rho_*)\mbox{.} 
\eqs
The resulting solution profile in the limit as $\epsilon \rightarrow 0$
is pictured in
Fig.~\ref{fig:Case2_RiemannSolution}\subref{fig:Case2_ShockProfile}:
\bq
    \rho(x, t) = \left\{
	\begin{array}{ll}
	  \rho_l,   & \text{if}~x < s t,\\
	  \rho_{m}, & \text{if}~s t \leqslant x \leqslant - \gamma t, \\
	  \rho_r,   & \text{if}~x > -\gamma t,
	\end{array}
    \right.
    \label{eqn:Case2_Solution}
\eq
which consists of a 1-shock moving to the left with speed
\bq
   s = \frac{g_c(\rho_{m}) - f(\rho_l)}{\rho_{m}-\rho_l} < 0,
  \label{eqn:Case2_ShockSpeed}
\eq
and a 2-contact moving to the left with speed $-\gamma$.  
\\
\\
\noindent{\bf Case~C: $\boldsymbol{\rho_l < \rho_{m} < \rho_r \;\;
    \mbox{and} \;\; \gamma/(\gamma+1) \geqslant \rho_l}$}.

\begin{figure}[btp]
  \begin{center}
    \subfigure[]{\label{fig:Case3_ConvexHull}\includegraphics[width=0.45\textwidth]{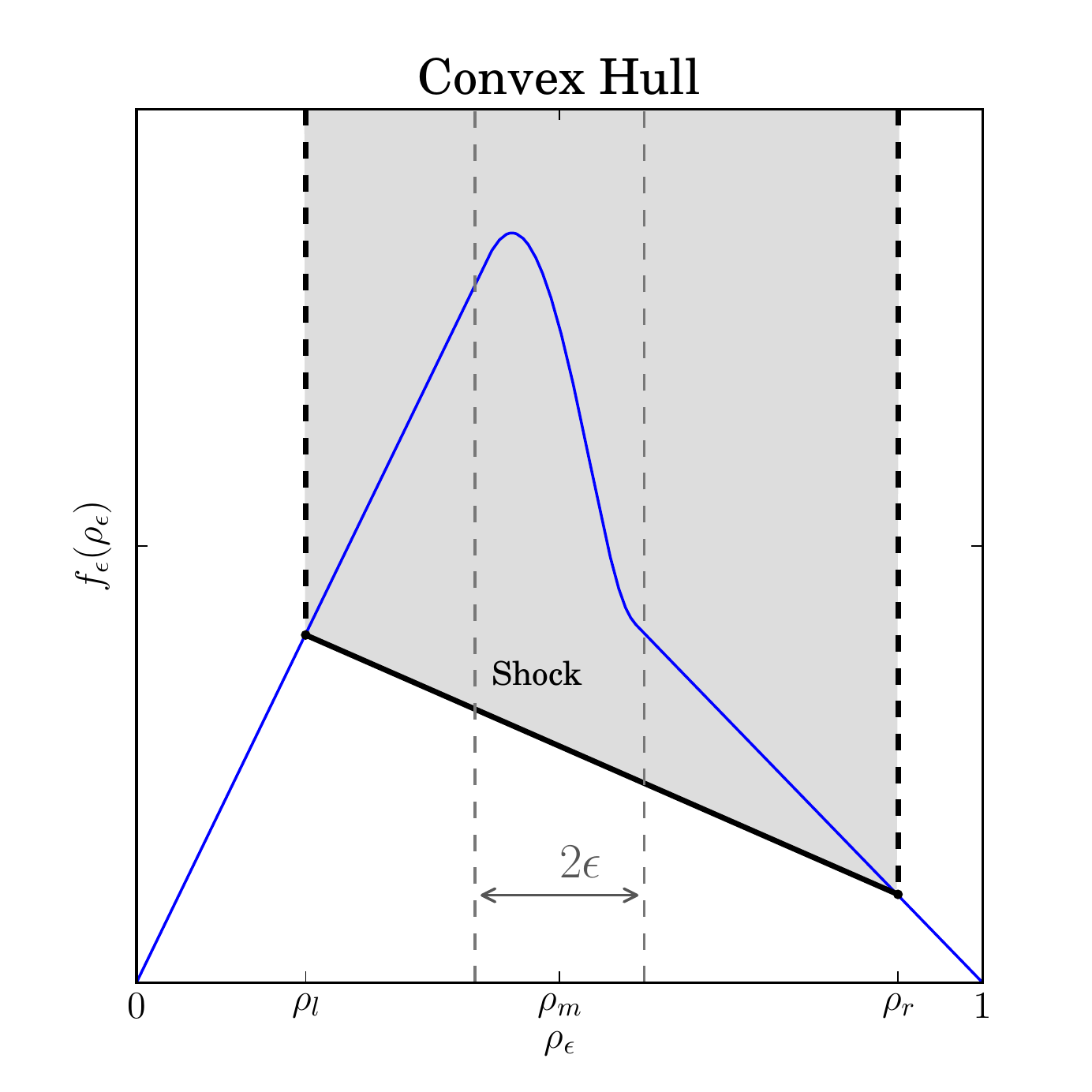}}
    \subfigure[]{\label{fig:Case3_ShockProfile}\includegraphics[width=0.45\textwidth]{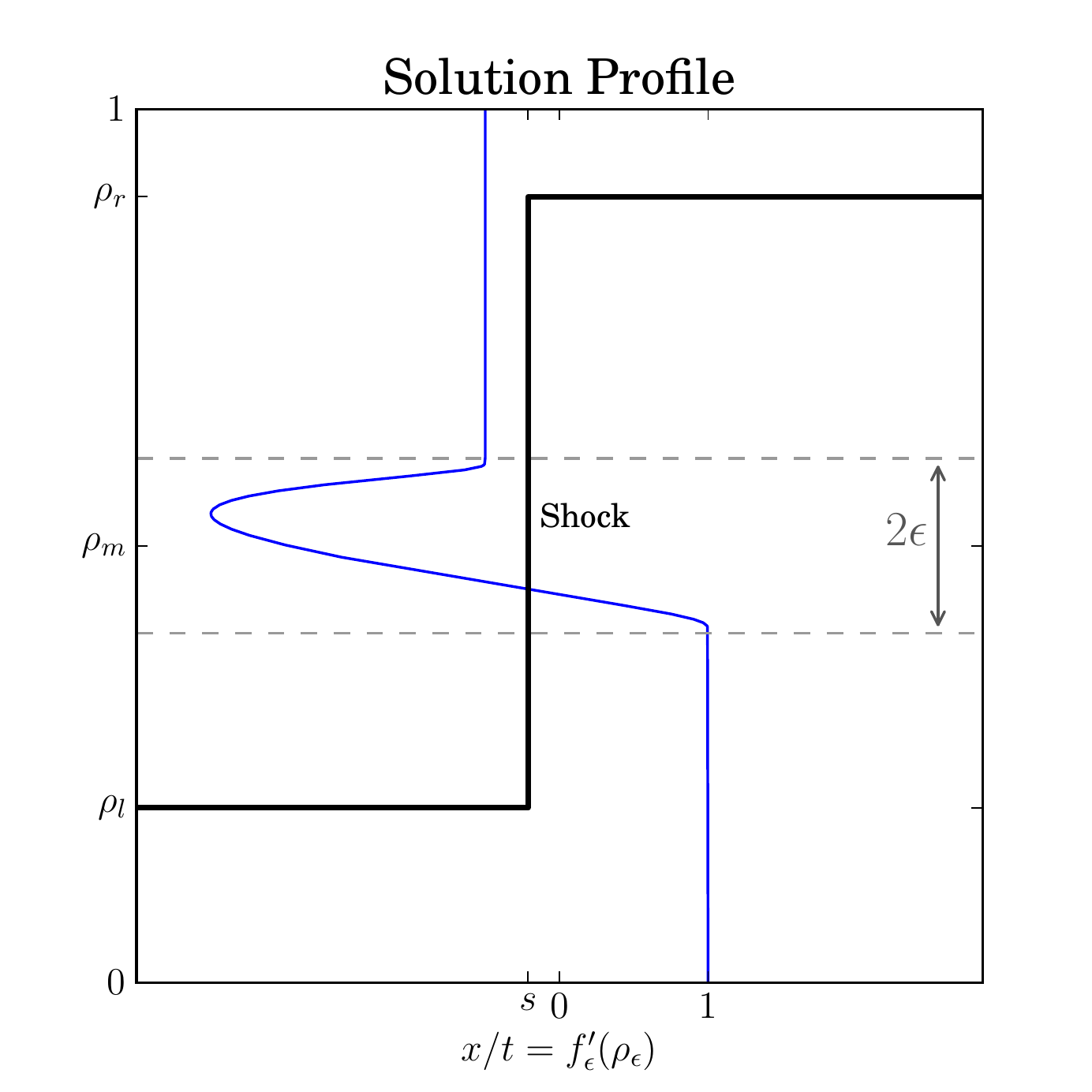}}
    \caption{The convex hull construction for Case~C.
      \subref{fig:Case3_ConvexHull} The thin (blue) line is the
      mollified flux function $f_\epsilon(\rho_\epsilon)$ and the thick
      (black) line is the convex hull.  The shaded oval region
      highlights the rarefaction wave. \subref{fig:Case3_ShockProfile}
      The solution $\rho_\epsilon(x/t)$ is shown as a thick (black)
      line, along with the corresponding flux derivative
      $f^\prime(\rho)$ as a thin (blue) line.  Parameter values are
      $\rho_l=0.2$, $\rho_r=0.9$, $\rho_{m}=0.5$, $\epsilon = 0.1$ and
      $\gamma = 0.5$.}
    \label{fig:Case3_RiemannSolution}
  \end{center}
\end{figure}

The solution structure in this case is significantly simpler than the
previous two in that there is only a single shock connecting the states
$\rho_r$ and $\rho_l$, and hence no intermediate state.  The convex hull
is depicted in
Fig.~\ref{fig:Case3_RiemannSolution}\subref{fig:Case3_ConvexHull} and
the corresponding solution profile in
Fig.~\ref{fig:Case3_RiemannSolution}\subref{fig:Case3_ShockProfile}.  As
$\epsilon \rightarrow 0$, the solution reduces to 
\bq
    \rho_\epsilon(x, t) = \left\{
	\begin{array}{ll}
	  \rho_l, & \text{if}~x < s t,\\
	  \rho_r, & \text{if}~x \geqslant s t, 
	\end{array}
    \right.
    \label{eqn:Case3_Solution}
\eq
which corresponds to a shock with speed
\bq
  s = \frac{f(\rho_r) - f(\rho_l)}{\rho_r-\rho_l},
  \label{eqn:Case3_ShockSpeed}
\eq
which can be either positive or negative depending on the sign of the
numerator. 

Note that as $\rho_l$ or $\rho_r$ approaches the discontinuity $\rho_m$,
the Riemann solution becomes sensitive to the choice of initial
data. This sensitivity is not unique to our problem, but is also
observed in other analytical solutions for problems with discontinuous
flux, such as in \cite{DiasFigueiraRodrigues2005,Gimse1993,Lu2009}.  In
the context of traffic flow, this sensitivity occurs in the
neighbourhood of the transition point between free-flow and congested
traffic, the exact location of which is expected to be highly sensitive
to the state of individual drivers comprising the flow.  Therefore, the
sensitivity in our model is consistent with actual traffic.

%% file: Sec4.tex
\section{Analysis of Zero Waves}
\label{sec:ZeroWaves}

In this section, we consider the two special situations that were not
addressed in Section~\ref{sec:RiemannProblem}, namely where either
$\rho_l=\rho_m$ or $\rho_r=\rho_m$.  In both cases, the mollified
problem gives rise to a {wave} having speed $O(1/\epsilon)$ and strength
$O(\epsilon)$, which we refer to as a \emph{zero rarefaction wave}
because of its similarity to the zero shocks identified by Gimse
\cite{Gimse1993}.  Since the speed of these waves becomes infinite as
$\epsilon \rightarrow 0$, information can be exchanged instantaneously
between neighbouring Riemann problems in any Godunov-type method.  We
demonstrate in this section how these effects can be incorporated into
the local Riemann solver.

\begin{figure}[tbp]
  \begin{center}
    \subfigure[Case~1.]{\label{fig:TripleRiemann_IC_Case1}\includegraphics[width=0.23\textwidth]{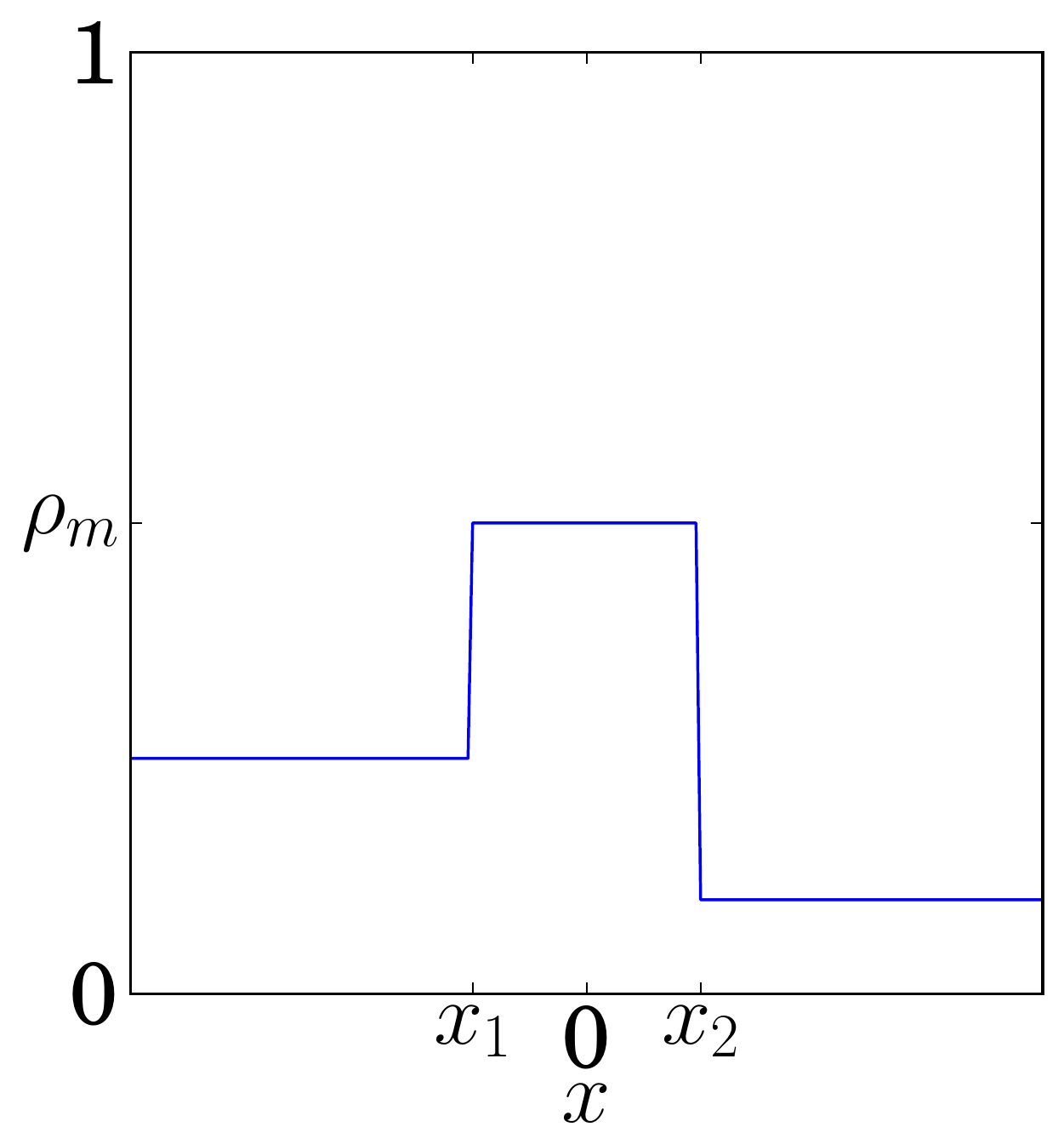}}
    \ 
    \subfigure[Case~2.]{\label{fig:TripleRiemann_IC_Case2}\includegraphics[width=0.23\textwidth]{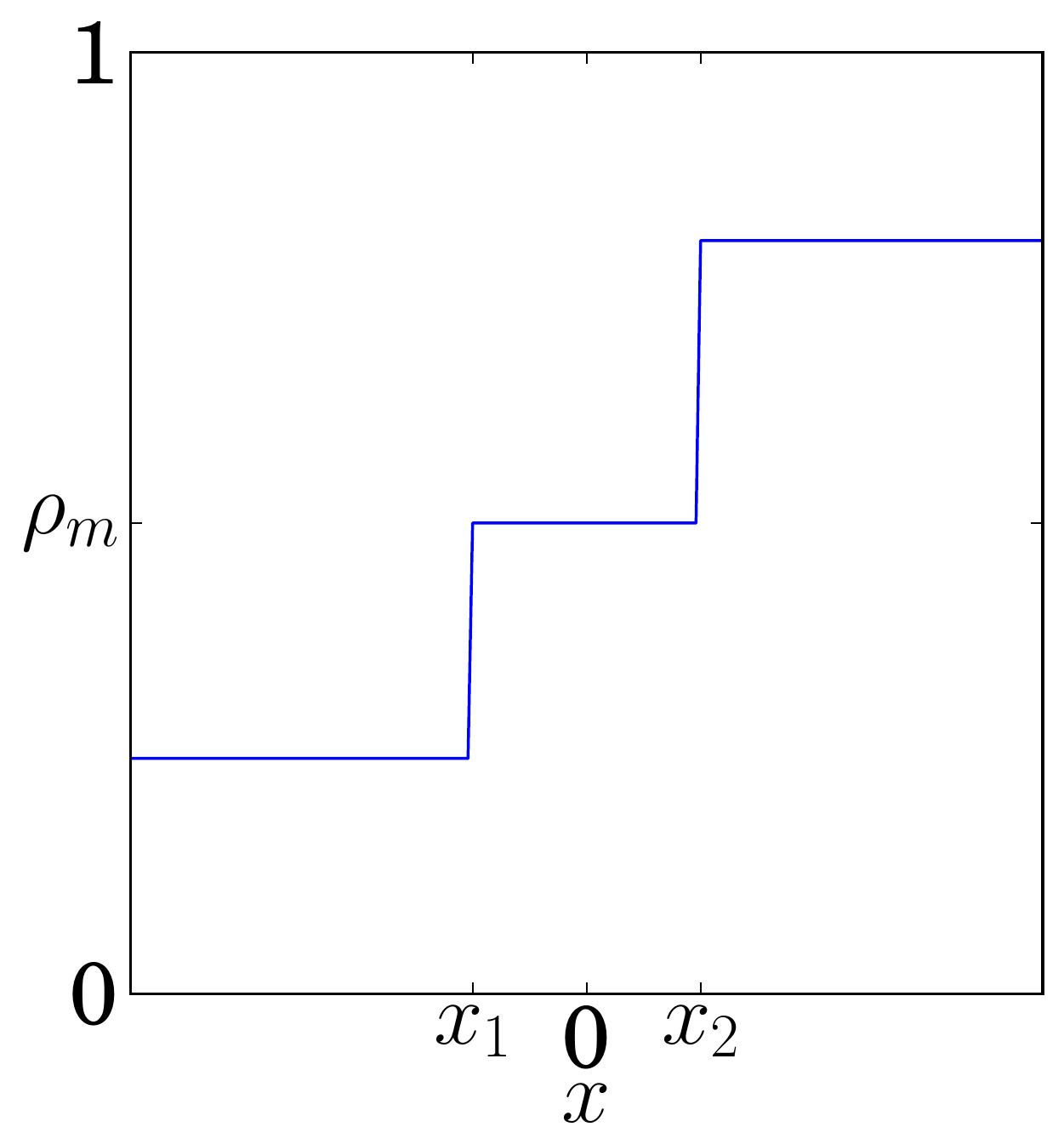}}
    \ 
    \subfigure[Case~3.]{\label{fig:TripleRiemann_IC_Case4}\includegraphics[width=0.23\textwidth]{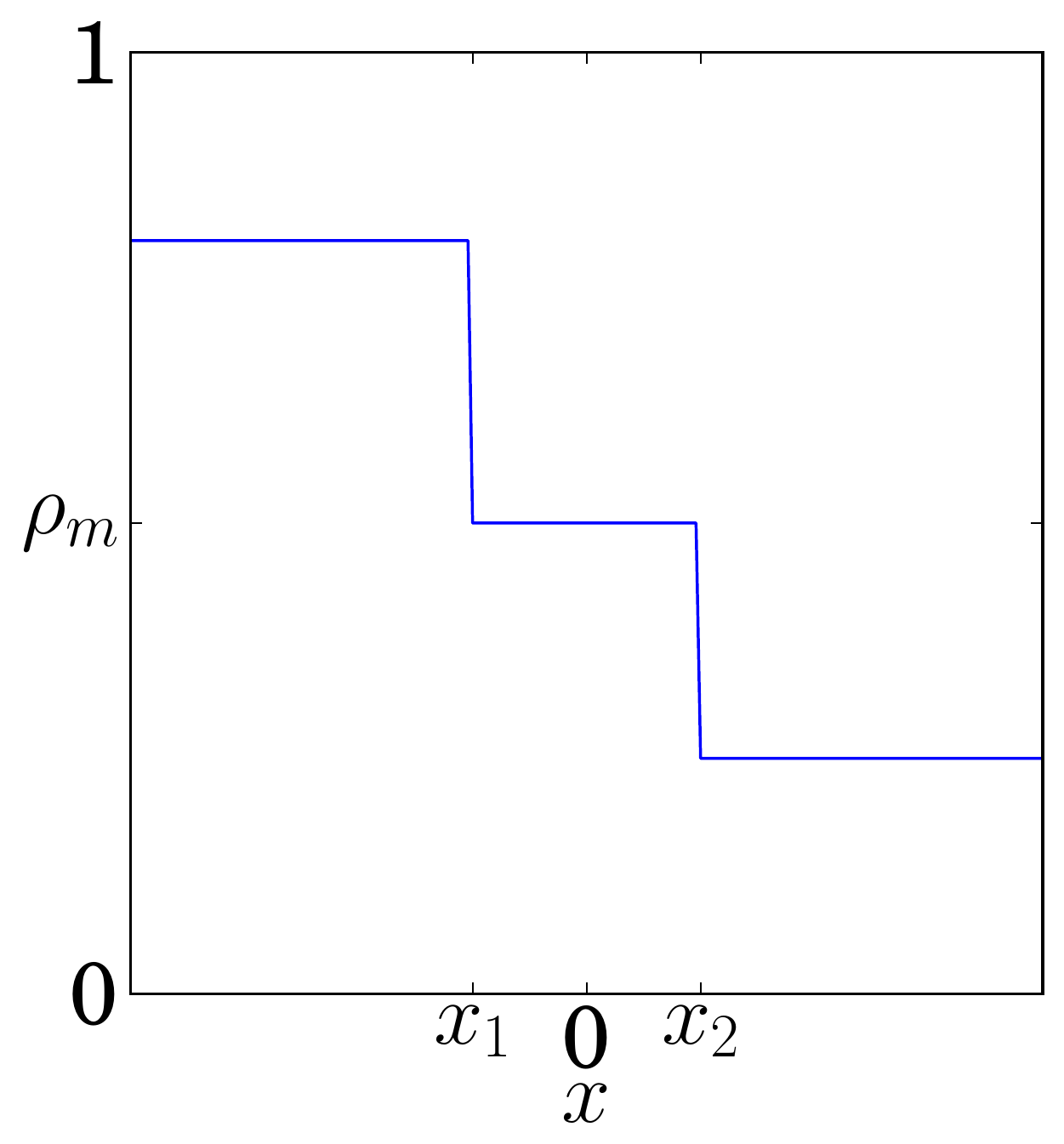}}
    \ 
    \subfigure[Case~4.]{\label{fig:TripleRiemann_IC_Case3}\includegraphics[width=0.23\textwidth]{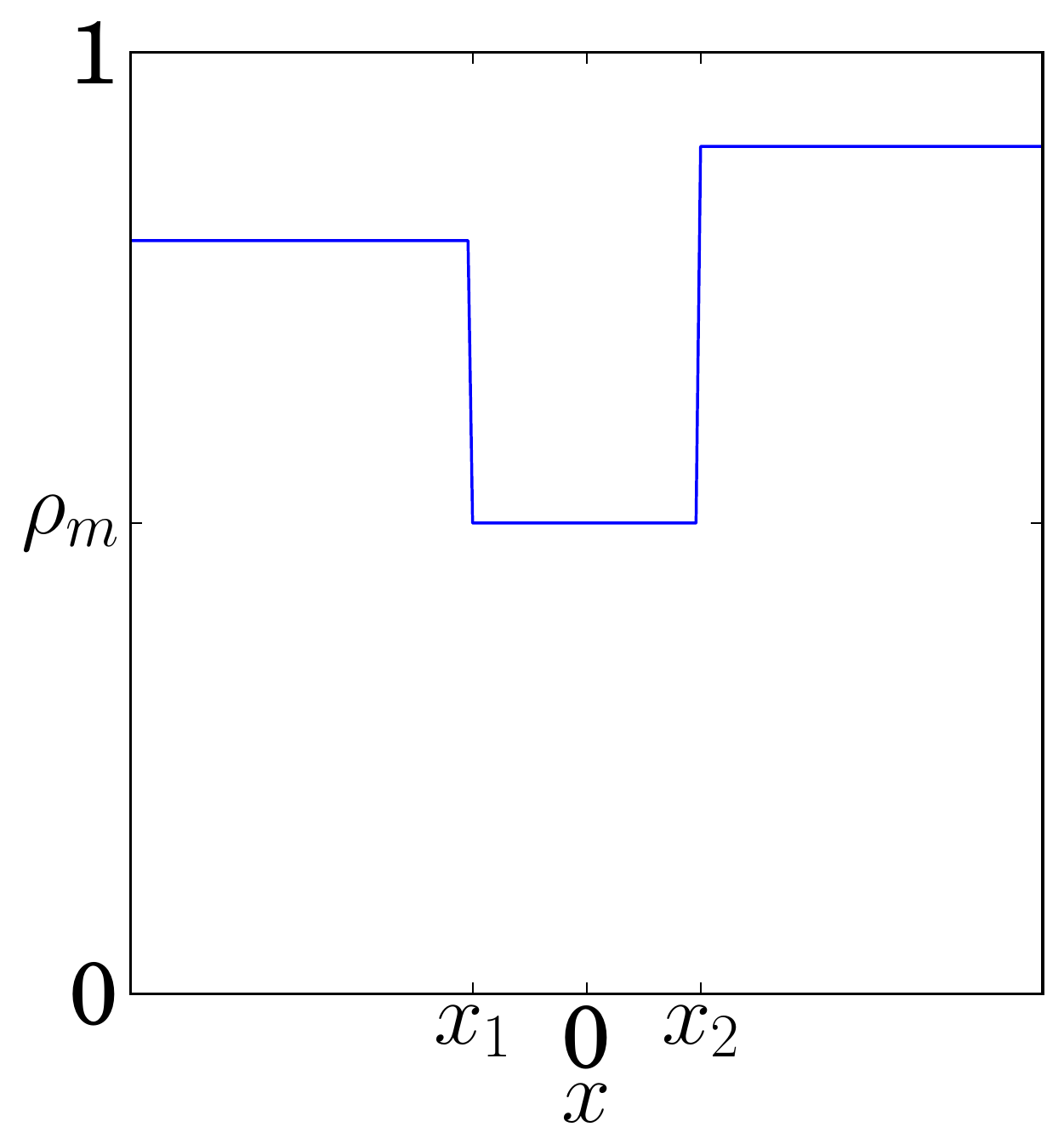}}
    \caption{Four possible cases for the initial conditions of the
      double Riemann problem \eqref{eqn:TripleRiemann_IC}.} 
    \label{fig:TripleRiemann_IC}
  \end{center}
\end{figure}

Motivated by the need to consider interactions between Riemann problems
arising from two pairs of piecewise constant states, we consider a
\emph{double Riemann problem} consisting of two ``usual'' Riemann
problems: one on the left with $(\rho_l,\rho_r) = (C_l,\rho_m)$, and a
second on the right with $(\rho_l,\rho_r) = (\rho_m,C_r)$.  As a result,
the mollified conservation law \eqref{eqn:MolConsLaw} is supplemented
with the following piecewise constant initial data pictured in
Fig.~\ref{fig:TripleRiemann_IC}
\bq
   \rho_\epsilon(x, 0) = \left\{
     \begin{array}{ll}
       C_l,    & \text{if}~x < x_1, \\
       \rho_m, & \text{if}~x_1 \leqslant x \leqslant x_2, \\
       C_r,    & \text{if}~x > x_2, 
     \end{array}
   \right.
\label{eqn:TripleRiemann_IC}
\eq 
where we have used the notation $C_l$ and $C_r$ for the left/right
states to emphasize the fact that we are solving a double Riemann
problem.

Since the solution can contain waves whose speed becomes unbounded, it
follows that we cannot solve the double Riemann problem with
intermediate state $\rho_m$ by simply splitting the solution into two
local Riemann problems and applying standard techniques.  Instead, the
origin and dynamics of zero waves need to be considered when
constructing the solution.  The double Riemann problem with a mollified
flux consists of up to three separate waves. Two of the waves -- which
we call the 1- and 2- wave using the standard terminology -- correspond
respectively to the left-most ($x_1$) and right-most ($x_2$) waves
arising from the initial discontinuities in the double Riemann problem.
We will see later that the 1-wave is a shock and the 2-wave is a contact
line.  The third wave corresponds to a left-moving {rarefaction wave}
having strength $O(\epsilon)$ and speed $O(1/\epsilon)$ that originates
from the $x_2$ interface and is located between the other two waves.  As
$\epsilon \rightarrow 0$, the rarefaction wave approaches a zero wave
that mediates an \emph{instantaneous} interaction between the 1- and
2-waves.

\subsection{Origin of Zero Waves}
\label{sec:CreateZeroRarefaction}

Since the local Riemann problem arising at the 2-wave (i.e., the contact
line at the $x_2$ interface) is the source of the {zero rarefaction
  wave}, we begin by focusing our attention on the right half of the
double Riemann problem.  The formation of a zero wave at the
discontinuity in the initial data located at $x_2$ can be divided into
two cases, corresponding to whether $C_r < \rho_m$ or $C_r > \rho_m$.
The specifics of the interaction between the 1-shock and the zero wave
will be treated separately Section~\ref{sec:ZeroRarefactionInteraction}.
\\
\\
\noindent{\bf Zero rarefaction with $\boldsymbol{C_r < \rho_m}$ and
  $\boldsymbol{C_l=\rho_m}$.} 

\begin{figure}[tbp]
  \begin{center}
    \includegraphics[width=0.40\textwidth]{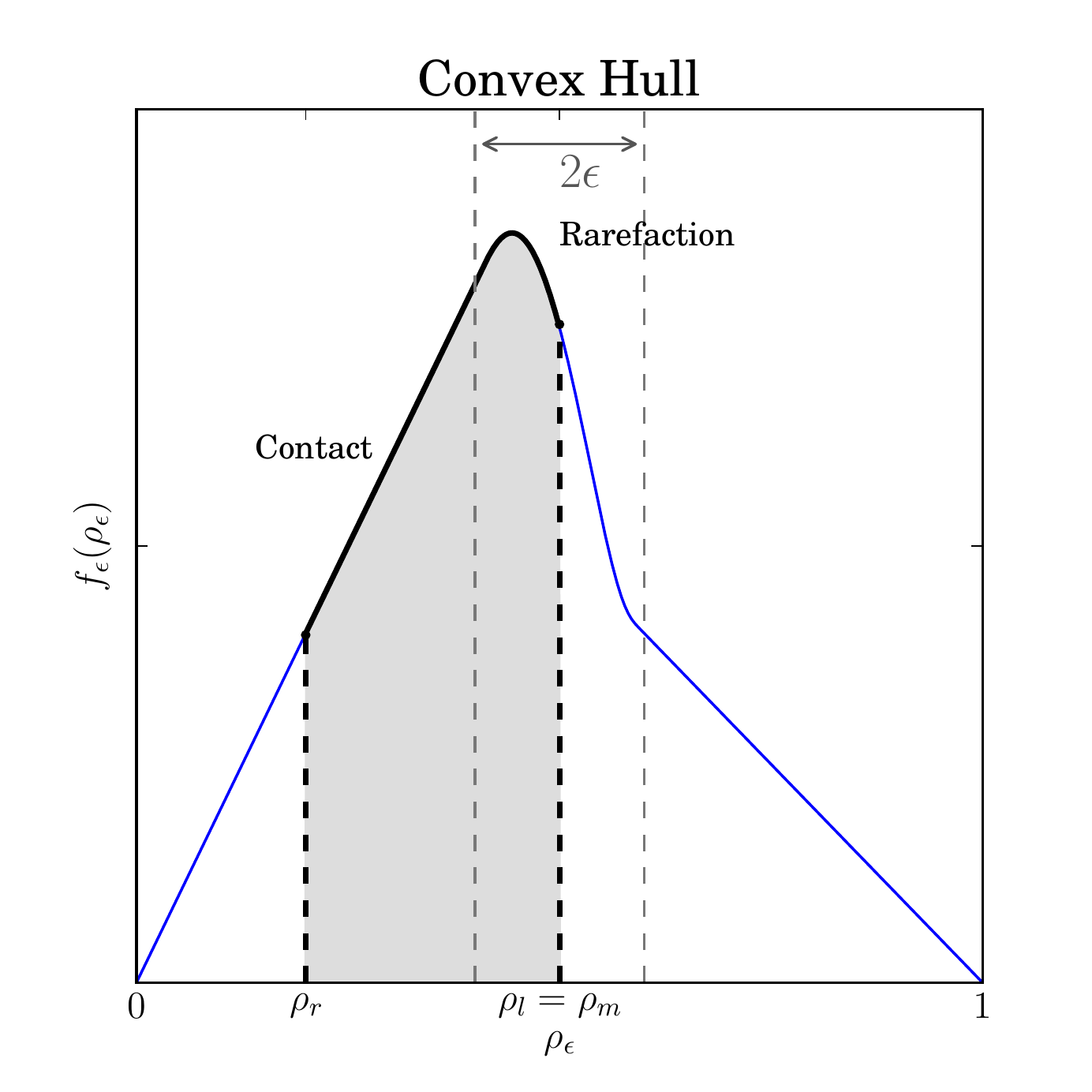}
    \includegraphics[width=0.40\textwidth]{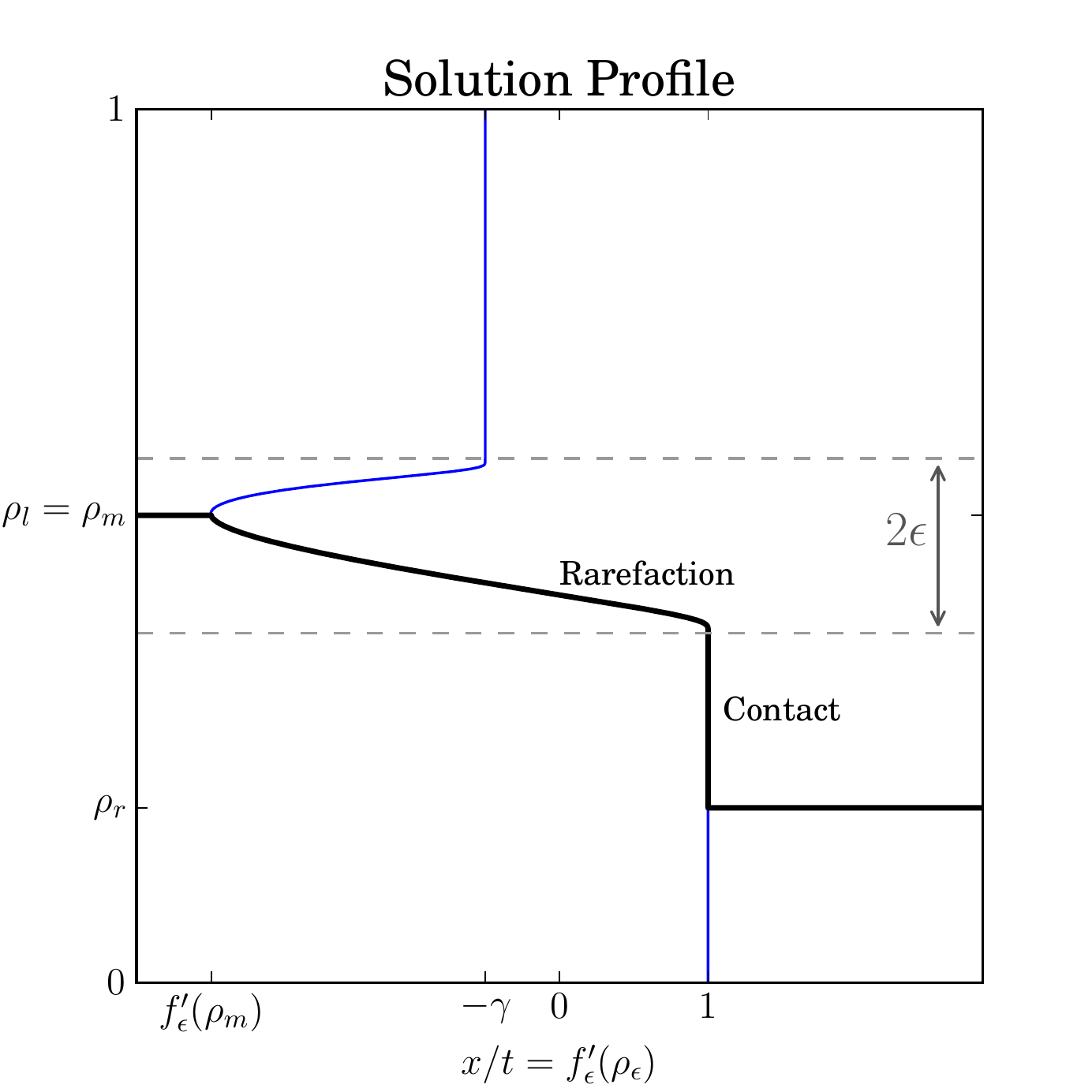}
    \caption{The convex hull (left, shaded) and solution profile (right)
      for the generation of the zero rarefaction wave when $\rho_l =
      \rho_m$ and $\rho_r < \rho_m$.}
    \label{fig:ZeroRarefaction_Case1}
  \end{center}
\end{figure}

We first consider the mollified Riemann problem with $C_l = \rho_l =
\rho_m$ and $C_r = \rho_r < \rho_m$.  The convex hull and solution
profile shown in Fig.~\ref{fig:ZeroRarefaction_Case1} exhibit a
right-moving contact line (the 2-wave) having speed $s=1$ and a
rarefaction wave of strength $O(\epsilon)$ travelling to the left with
speed $O(1/\epsilon)$.  For a simple isolated Riemann problem, the
solution would reduce to a lone contact line as $\epsilon \rightarrow 0$
and the zero rarefaction would have no impact.  However, when the zero
wave is allowed to interact with the solution of another neighbouring
Riemann problem -- such as when multiple Riemann problems are solved on
a sequence of grid cells in a Godunov-type numerical scheme -- the local
Riemann problems cannot be taken in isolation.
\\
\\
\noindent{\bf Zero rarefaction with $\boldsymbol{C_r > \rho_m}$ and
  $\boldsymbol{C_l=\rho_m}$.}  

\begin{figure}[tbp]
  \begin{center}
    \includegraphics[width=0.40\textwidth]{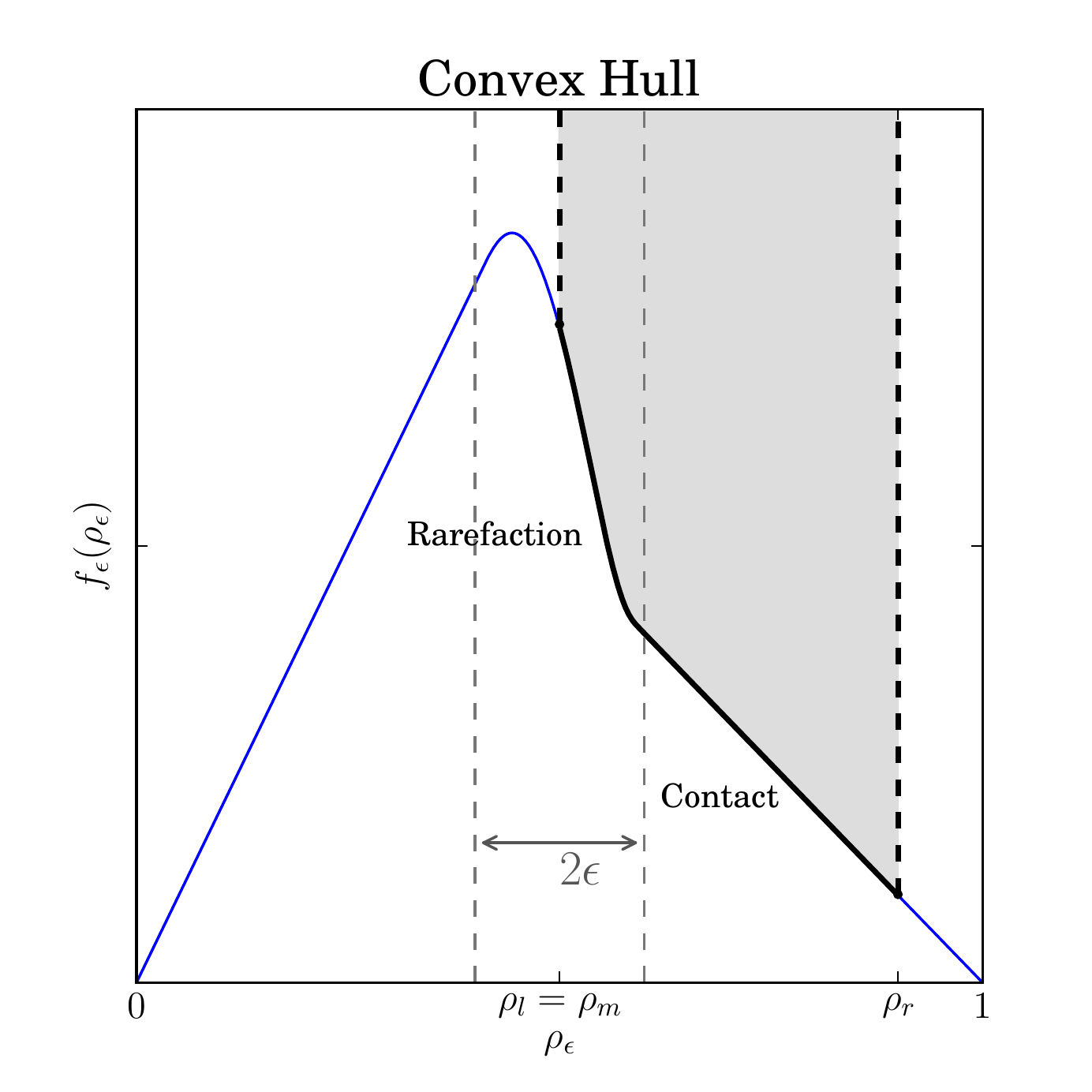}
    \includegraphics[width=0.40\textwidth]{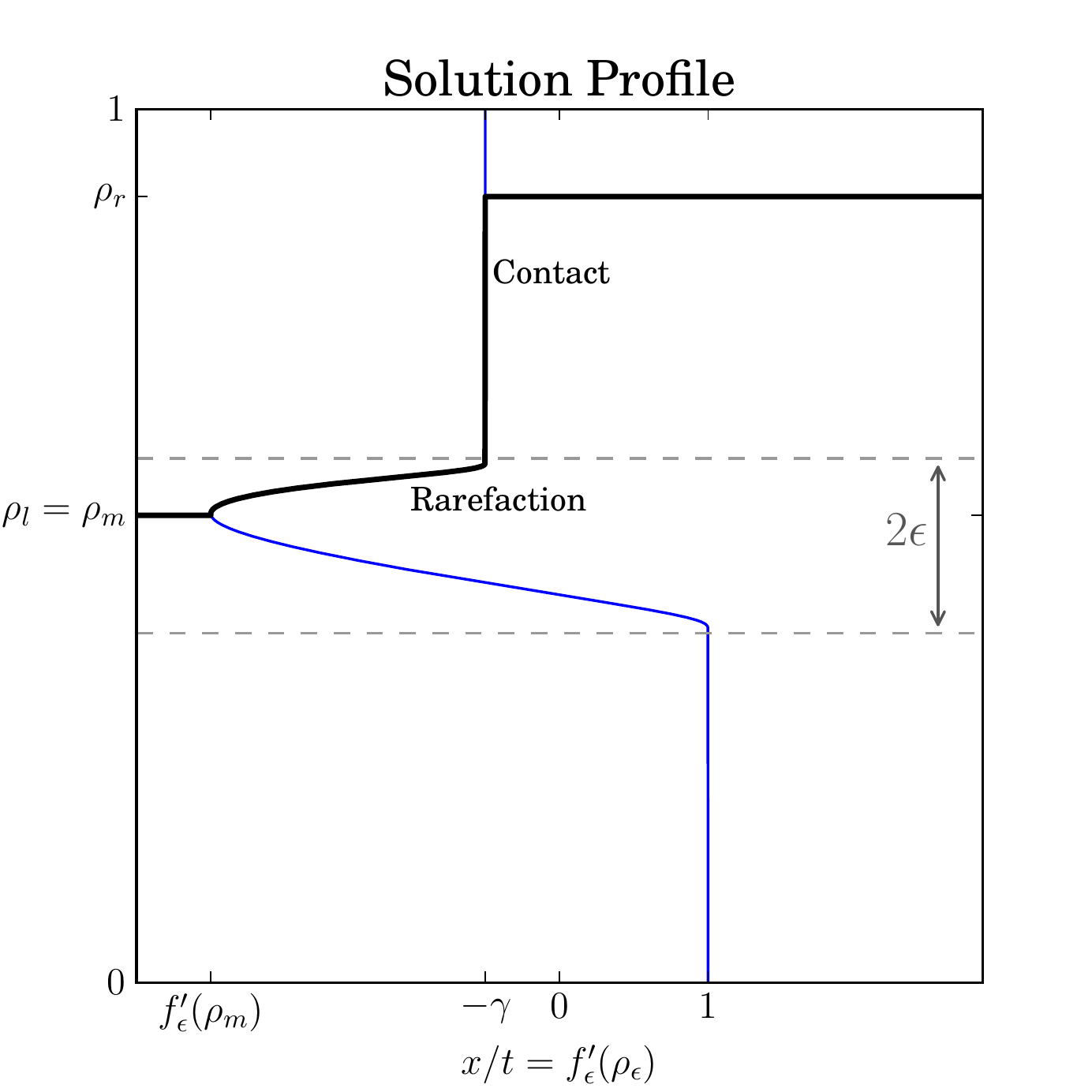}
    \caption{The convex hull (left, shaded) and solution profile (right)
      for the generation of the zero rarefaction wave when $\rho_l =
      \rho_m$ and $\rho_r > \rho_m$.}
    \label{fig:ZeroRarefaction_Case2}
  \end{center}
\end{figure}

Next we consider the mollified Riemann problem with $C_l = \rho_l =
\rho_m$ and $C_r = \rho_r > \rho_m$.  The convex hull and solution are
depicted in Fig.~\ref{fig:ZeroRarefaction_Case2}, and the solution again
consists of a contact line and {zero rarefaction} wave.  The main
difference from the previous case is that the contact line travels to
the left with speed $-\gamma$ instead of to the right.

Note that since the location of the inflection point of
$f_\epsilon(\rho_\epsilon)$ approaches $\rho_m$ as $\epsilon \rightarrow
0$, an additional {zero shock} of strength $O(\epsilon^2)$ and speed
$O(1/\epsilon)$ is generated.  For the sake of clarity, we have
not included this wave in Fig.~\ref{fig:ZeroRarefaction_Case2} since it
is of higher order than the $O(\epsilon)$ zero rarefaction and so has
negligible impact on the solution in the limit as $\epsilon \rightarrow
0$.

\subsection{Interaction Between 1-Shock and Zero Rarefaction Wave} 
\label{sec:ZeroRarefactionInteraction}

When solving the double Riemann problem, we need to determine how the
{zero rarefaction} wave produced at $x_2$ interacts with the 1-wave,
which we will see shortly must be a shock.  The details of the
interaction can be studied using the method of characteristics for the
four cases shown in Fig.~\ref{fig:TripleRiemann_IC}.  Since the {zero
  wave} has speed $O(1/\epsilon)$, we must examine the shock-{zero
  rarefaction} interaction on two different time scales of length
$O(\epsilon)$ and $O(1)$.  Since $\lambda_1 = O(1)$ as $\epsilon
\rightarrow 0$ in each case, the problem can be simplified significantly
by neglecting the shock dynamics on the $O(\epsilon)$ time scale.  We
will then show that the 1-wave approaches a constant speed as $\epsilon
\rightarrow 0$ when $t=O(1)$.
\\
\\
\noindent{\bf Case~1: $\boldsymbol{C_r,C_l < \rho_m}$.}

\begin{figure}[tbp]
  \begin{center}
    \subfigure[]{\label{fig:ZeroRarefactionInteraction_Plots_Case1_Plot1}\includegraphics[width=0.37\textwidth]{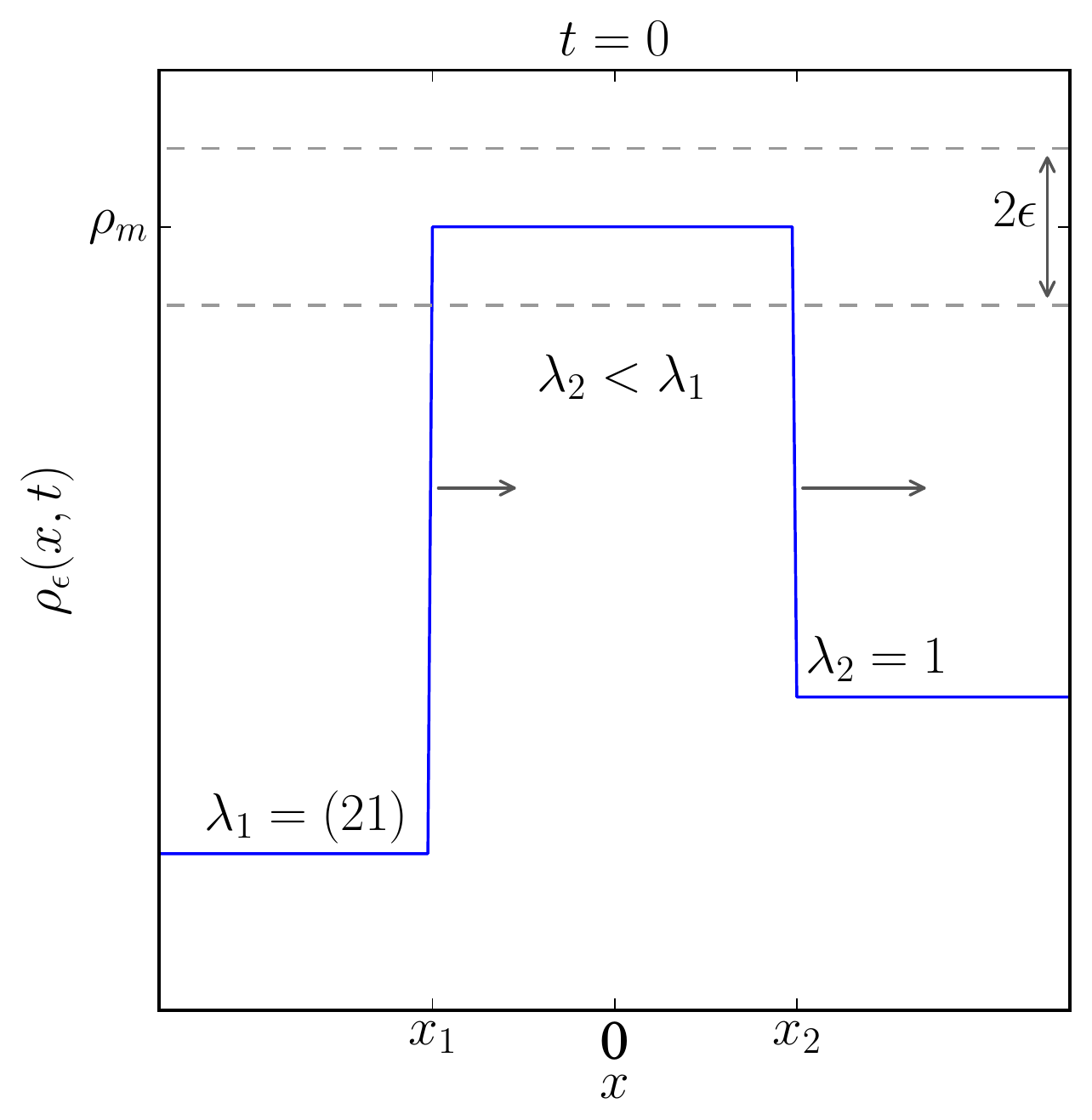}}
    \qquad
    \subfigure[]{\label{fig:ZeroRarefactionInteraction_Plots_Case1_Plot2}\includegraphics[width=0.37\textwidth]{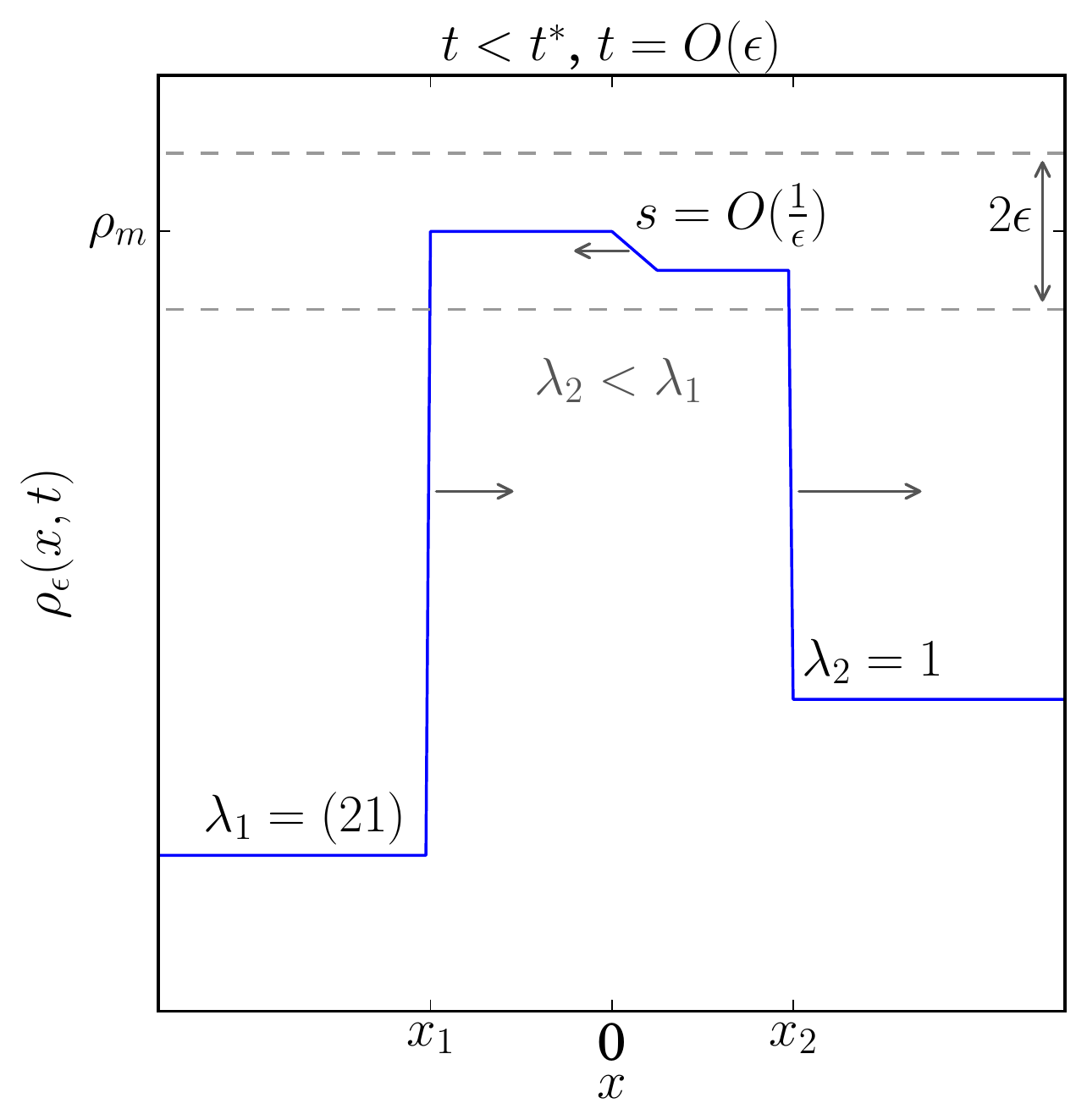}}\\
    \subfigure[]{\label{fig:ZeroRarefactionInteraction_Plots_Case1_Plot3}\includegraphics[width=0.37\textwidth]{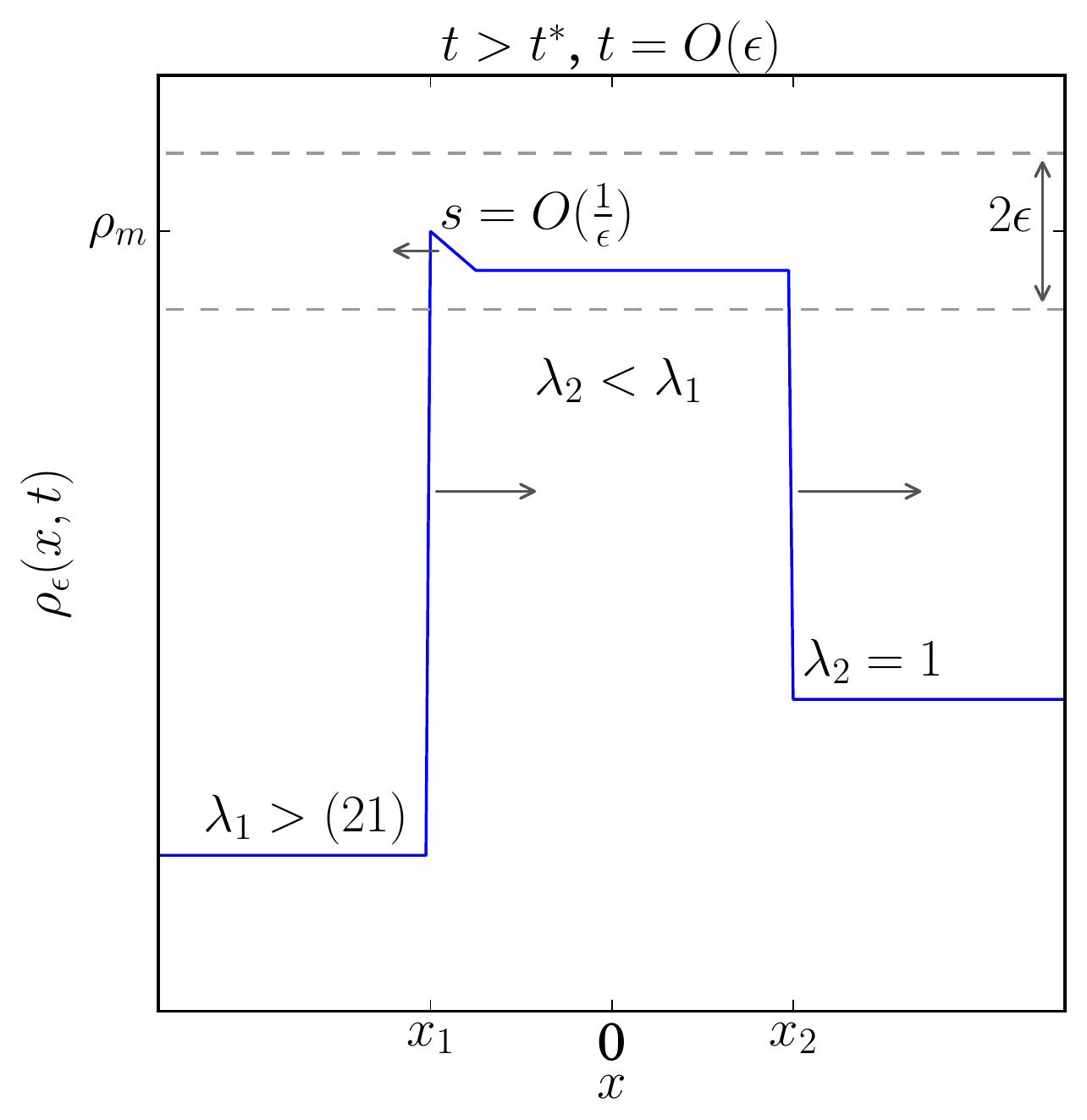}}
    \qquad
    \subfigure[]{\label{fig:ZeroRarefactionInteraction_Plots_Case1_Plot4}\includegraphics[width=0.37\textwidth]{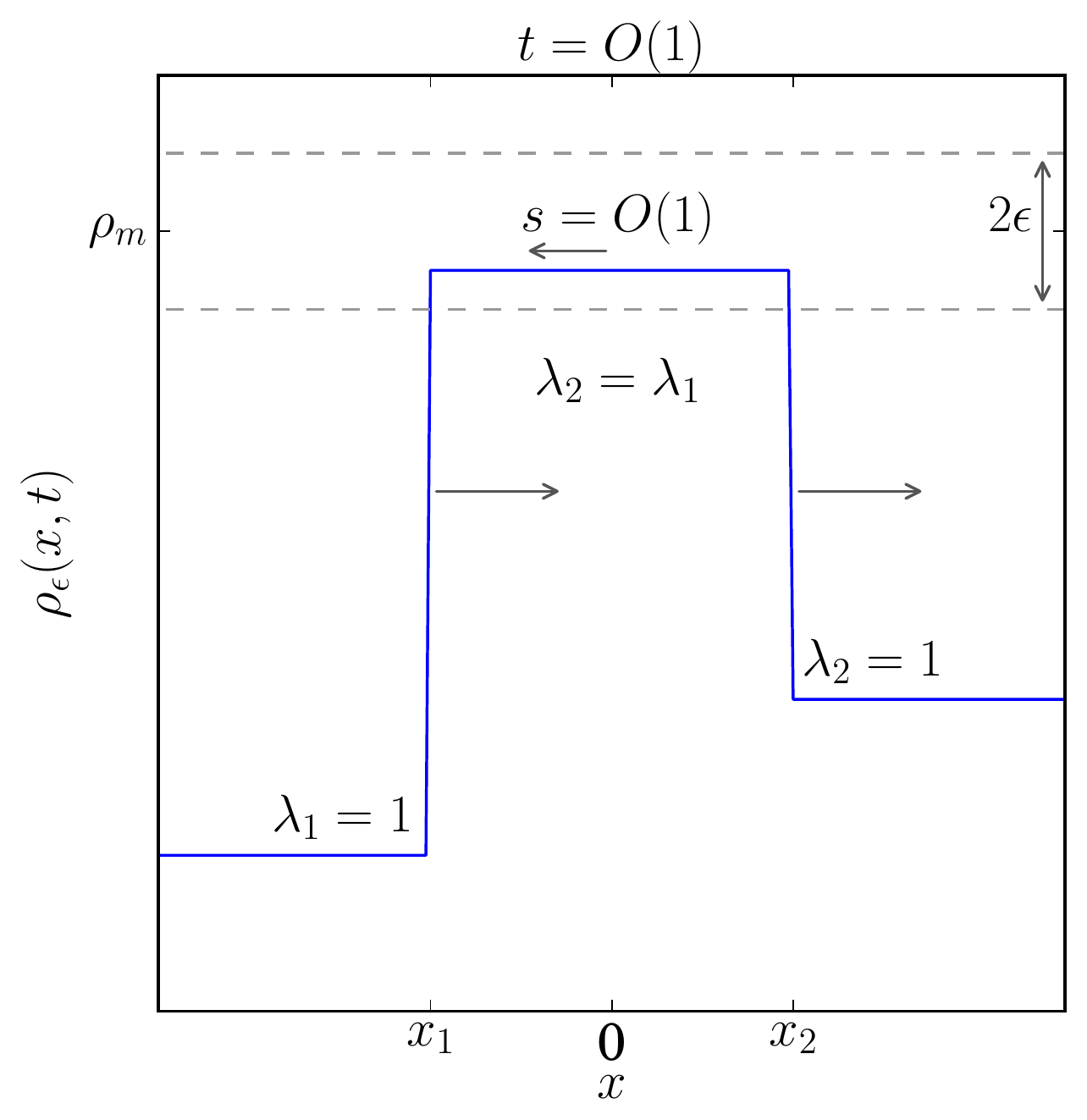}}
    \caption{Time evolution of the double Riemann solution in
      Case~1 when $C_r,C_l < \rho_m$.  A similar picture holds in
      Case~4 when $C_r,C_l > \rho_m$, except that the wave directions
      are reversed.}
    \label{fig:ZeroRarefactionInteraction_Plots_Case1}
  \end{center}
\end{figure}

\begin{figure}[tbp]
  \begin{center}
    \footnotesize
    \begin{tabular}{ccc}
      Case 1: $C_r,C_l < \rho_m$ 
      & & 
      Case 2a: $C_r>\rho_m$ and $\gamma/(\gamma+1) \leqslant C_l < \rho_m$
      \\  
      \includegraphics[width=0.45\textwidth]{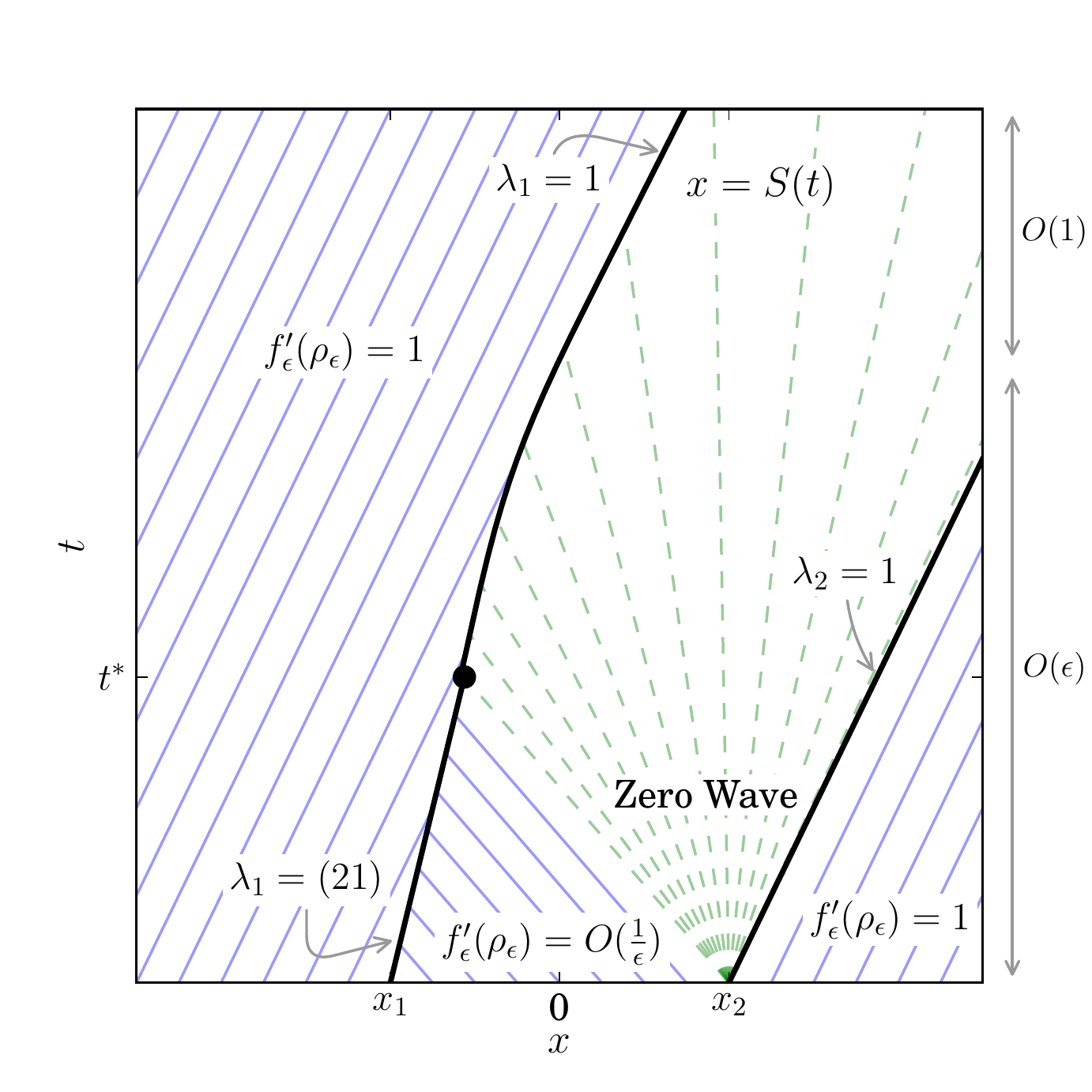}
      & &
      \includegraphics[width=0.45\textwidth]{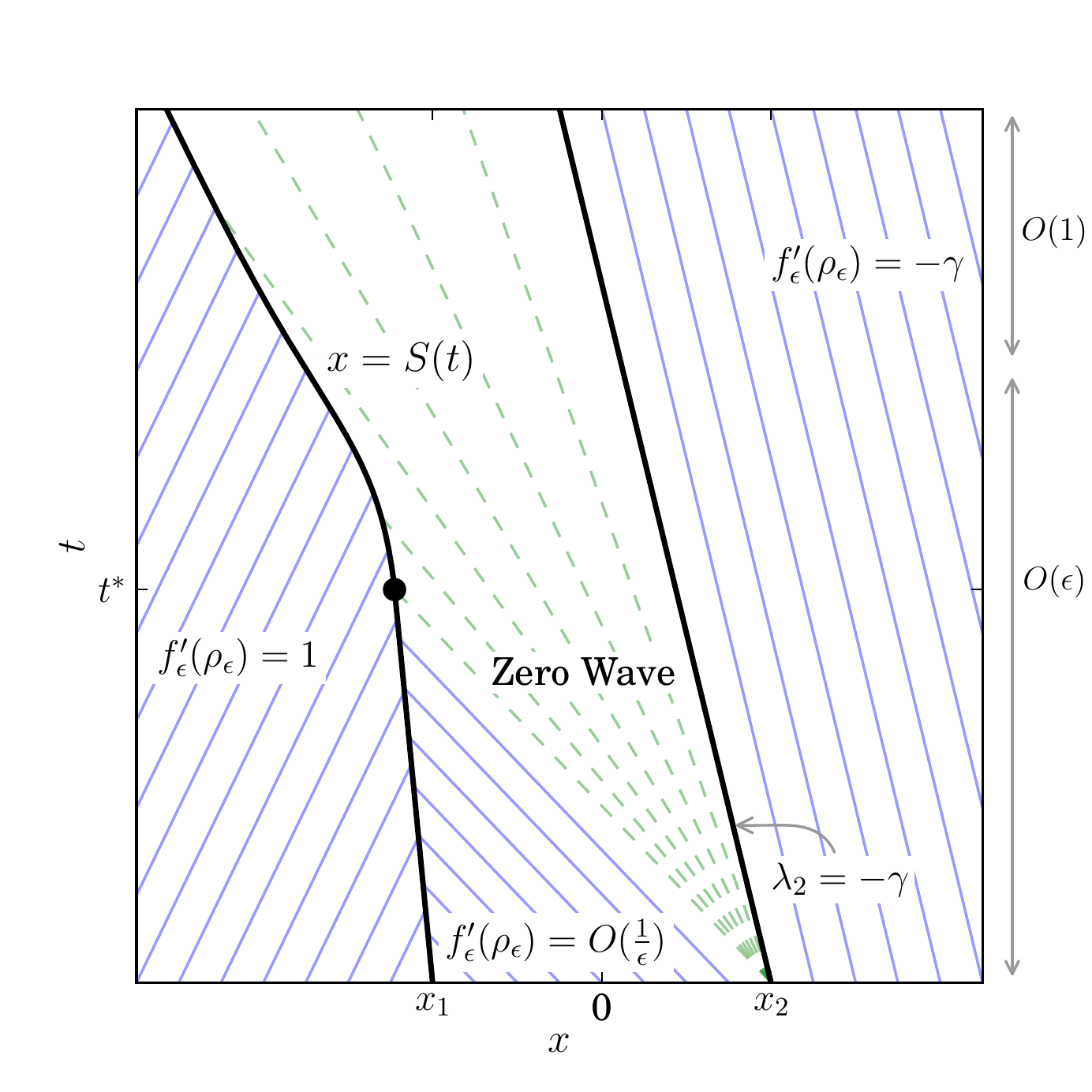}
      \\
      Case 2b: $C_r>\rho_m$ and $C_l < \gamma/(\gamma+1) < \rho_m$
      & & 
      Case 3: $C_r < \rho_m < C_l$
      \\
      \includegraphics[width=0.45\textwidth]{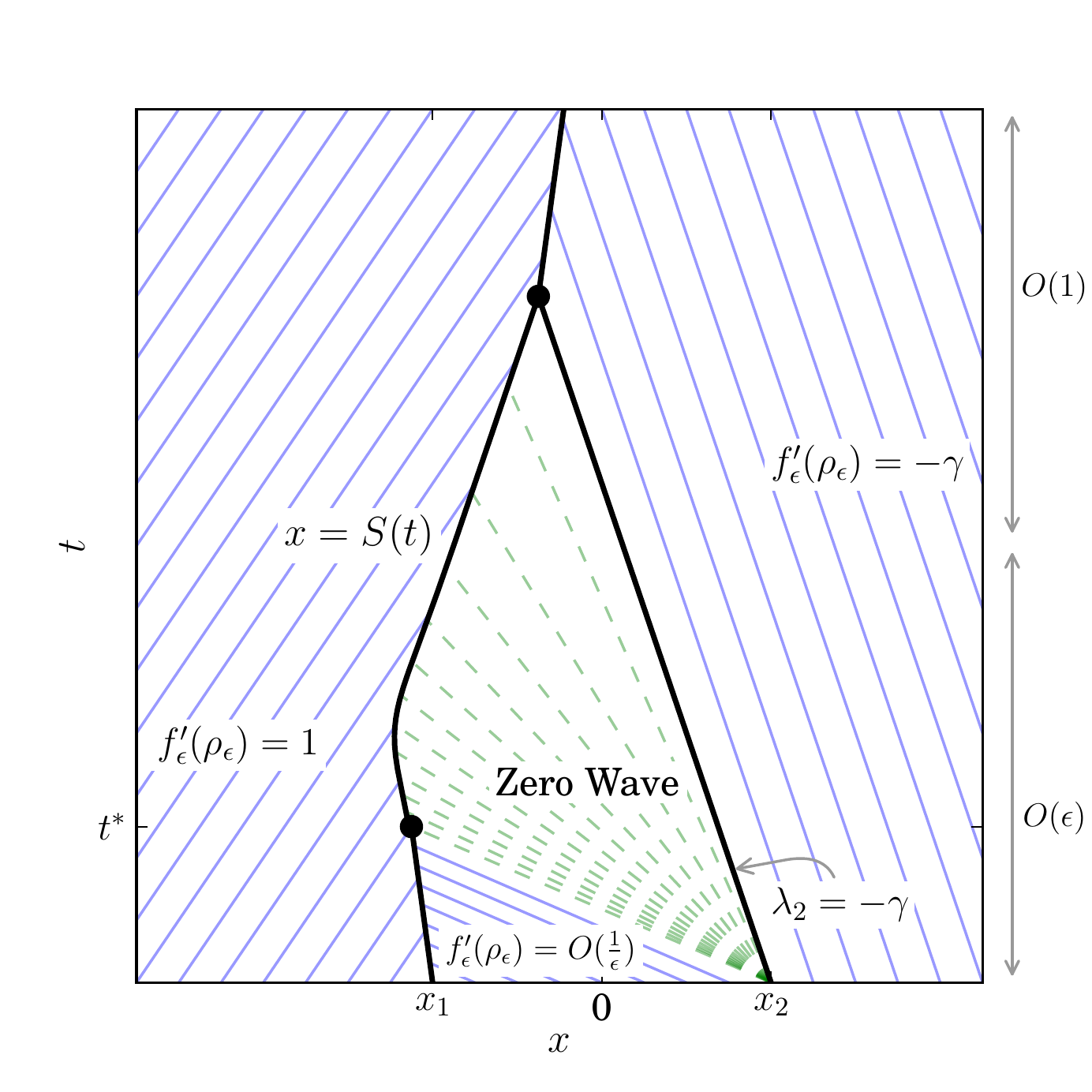}
      & & 
      \includegraphics[width=0.45\textwidth]{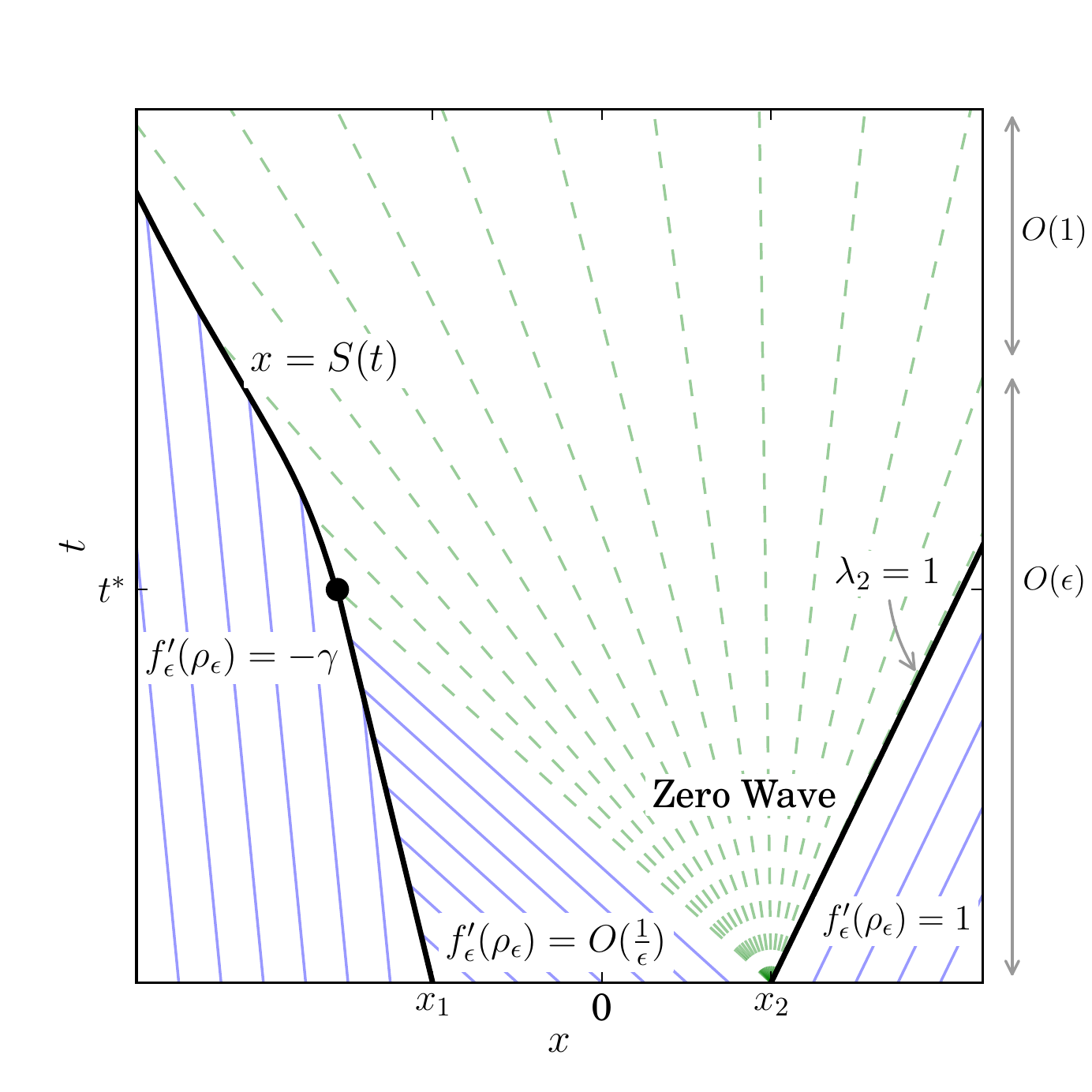}
    \end{tabular}
    \caption{Characteristic lines and elementary waves in the
      $x,t$--plane for various cases for the double Riemann
      solution.  The thick lines (black) represent shocks while thin
      lines (blue and green) represent characteristics.} 
    \label{fig:ZeroRarefactionInteraction_PhasePlot_Case1-4}
  \end{center}
\end{figure}

We first consider solutions to the mollified conservation law
\eqref{eqn:MolConsLaw} having initial conditions
\eqref{eqn:TripleRiemann_IC} that satisfy $C_r,C_l < \rho_m$ as shown in
Fig.~\ref{fig:TripleRiemann_IC}\subref{fig:TripleRiemann_IC_Case1}.  The
solution in this case consists of three waves: a shock, a {zero
  rarefaction} wave, and a contact line. The initial discontinuity at
$x_2$ generates a contact line (the 2-wave) that travels at speed
$\lambda_2=1$, along with a {zero rarefaction} wave.  A 1-shock
originates from location $x_1$ and travels along the trajectory
$x=S(t)$, where $S(t)$ has yet to be determined.  Suppose that the {zero
  rarefaction} wave intersects the 1-wave at time $t^*$; then the speed
of the 1-wave satisfies the Rankine-Hugoniot condition
\bq
  \frac{d S(t)}{dt} = \lambda_1(t) = \frac{f_\epsilon(C_l) -
    f_\epsilon(\rho_m)}{C_l - \rho_m}, \qquad \text{for\;} t \leqslant t^* \simeq
  {\mathcal O}(\epsilon) \mbox{.} 
  \label{eqn:BeforeTStarLambda1ShockSpeed}
\eq
The time evolution of the solution for $t<t^*$ is illustrated in
Fig.~\ref{fig:ZeroRarefactionInteraction_Plots_Case1}\subref{fig:ZeroRarefactionInteraction_Plots_Case1_Plot2},
and the corresponding plot of characteristics in the $x,t$--plane
is shown in Fig.~\ref{fig:ZeroRarefactionInteraction_PhasePlot_Case1-4}.
The characteristics that intersect with the 1-wave for $t<t^*$ have
speed equal to $1$ to the left of the 1-wave, and speed $O(1/\epsilon)$
to the right.

At time $t^*$, the {zero rarefaction} wave starts to interact with the
shock $x=S(t)$ which decreases the density $\rho_\epsilon$ to the right
of the 1-wave as illustrated in
Fig.~\ref{fig:ZeroRarefactionInteraction_Plots_Case1}\subref{fig:ZeroRarefactionInteraction_Plots_Case1_Plot3}.
As a result, the shock wave attenuates leading to an increase in the
shock speed $\lambda_1$.
Since the {zero wave} contains $\rho_\epsilon$ values lying within the
interval $[\rho_m-\epsilon, \rho_m]$, we can determine the portion of
the shock--{zero wave} interaction that occurs on the $\epsilon$ time
scale by finding the range of $\rho_\epsilon$ for which
$f^\prime_\epsilon(\rho_\epsilon) = O(1/\epsilon)$.
Equation~\eqref{eqn:FirstDerivativeFlux} implies that
\bq
  f^\prime_\epsilon(\rho_\epsilon) =
  O\left(\frac{1}{\epsilon}\eta\left(\frac{\rho_\epsilon-\rho_m}{\epsilon}\right)\right)
  \mbox{,}  
  \label{eqn:MollifiedFluxDerivativeRelationship}
\eq
so that we only need to determine the range of $\rho_\epsilon$ where
\bq
  \eta\left(\frac{\rho_\epsilon-\rho_m}{\epsilon}\right) = O(1)\mbox{.}
  \label{eqn:OrderOfEta}
\eq
Using the formula for the mollifier \eqref{eqn:BumpMollifier}, it is
easy to verify that \eqref{eqn:OrderOfEta} holds when
\bq
  \rho_\epsilon = \rho_m + \alpha_\epsilon \epsilon, 
  \qquad
  \lim_{\epsilon \rightarrow 0}\left( \alpha_\epsilon^2 - 1 \right) \ne
  0
  \qquad \text{and} \qquad 
  \alpha_\epsilon\in(-1,1) \mbox{.} 
  \label{eqn:RhoWithInfiniteSpeed}
\eq
It is only as $\rho_\epsilon$ approaches the boundary of the
$\epsilon$-region and $\lim_{\epsilon \rightarrow 0}\left(
  \alpha_\epsilon^2 - 1 \right) = 0$, that
$f^\prime_\epsilon(\rho_\epsilon) = o(1/\epsilon)$. Therefore, the
values of $\rho_\epsilon$ within the {zero rarefaction wave} that
satisfy \eqref{eqn:RhoWithInfiniteSpeed} will interact with the 1-shock
on the $\epsilon$ time scale.

Because $\lambda_1(t) = O(1)$, we can ignore the influence of the shock
over the $\epsilon$ time scale when $\epsilon \rightarrow 0$.
Therefore, we only need to determine the interaction between the shock
and the {zero rarefaction} on the $O(1)$ time scale. By finding the
range of $\rho_\epsilon$ within the {zero wave} where
$f^\prime_\epsilon(\rho_\epsilon) = O(1)$, we can show that the shock
speed $\lambda_1$ approaches a constant as $\epsilon \rightarrow 0$.  By
using the relationship \eqref{eqn:MollifiedFluxDerivativeRelationship},
we know that there exists a $\phi$ such that
$f^\prime_\epsilon(\rho_\epsilon) = O(1)$ for all
\bq
  \rho_\epsilon = \rho_m - \epsilon + \phi(\epsilon), 
  \label{eqn:RhoWithOrderOneSpeed}
\eq
where $\phi(\epsilon)=o(\epsilon)$ and $\phi(\epsilon)>0$.  Next, by bounding the
integral
\begin{align*}
  \displaystyle \int_{\rho_m-\epsilon}^{\rho_m - \epsilon +
    \phi(\epsilon)} \eta_\epsilon(s-\rho_m) \, ds \leqslant~& [\rho_m -
  \epsilon + \phi(\epsilon) - (\rho_m-\epsilon) ] \eta_\epsilon(\rho_m -
  \epsilon + \phi(\epsilon)-\rho_m), \\ 
  =~& \phi(\epsilon) \eta_\epsilon(-\epsilon + \phi(\epsilon)), \\
  =~& O(\phi(\epsilon))\mbox{,}
\end{align*}
we know from Eq.~\eqref{eqn:MollifiedFlux} that
\bqs
  f_\epsilon(\rho_\epsilon) = \rho_\epsilon + O(\phi(\epsilon)),
\eqs
for the range of $\rho_\epsilon$ defined by
Eqs.~\eqref{eqn:RhoWithOrderOneSpeed}.  Therefore, as $\epsilon
\rightarrow 0$, we have that $f_\epsilon(\rho_\epsilon) \rightarrow
g_f(\rho_m)$ and the shock speed $\lambda_1 \rightarrow 1$.
This results in a solution of the form
\bq
  \rho(x,t) = \left\{
    \begin{array}{ll}
      C_l,    & \text{if}~x < x_1 + \lambda_1 t, \\
      \rho_m, & \text{if}~x_1 + \lambda_1 t \leqslant x \leqslant x_2 + \lambda_2 t, \\
      C_r,    & \text{if}~x > x_2 + \lambda_2 t,
    \end{array}
  \right.
  \label{eqn:TripleRiemann_Ansatz}
\eq
where $\lambda_1 = \lambda_2 = 1$.  Therefore, at longer times the
solution takes the form of a ``square wave'' propagating to the
right at constant speed as pictured in
Fig.~\ref{fig:ZeroRarefactionInteraction_Plots_Case1}\subref{fig:ZeroRarefactionInteraction_Plots_Case1_Plot4}.
\\
\\
\noindent{\bf Case~2: $\boldsymbol{C_l < \rho_m < C_r}$.}

\begin{figure}[tbp]
  \begin{center}
    \subfigure[]{\label{fig:ZeroRarefactionInteraction_Plots_Case3_Plot1}\includegraphics[width=0.37\textwidth]{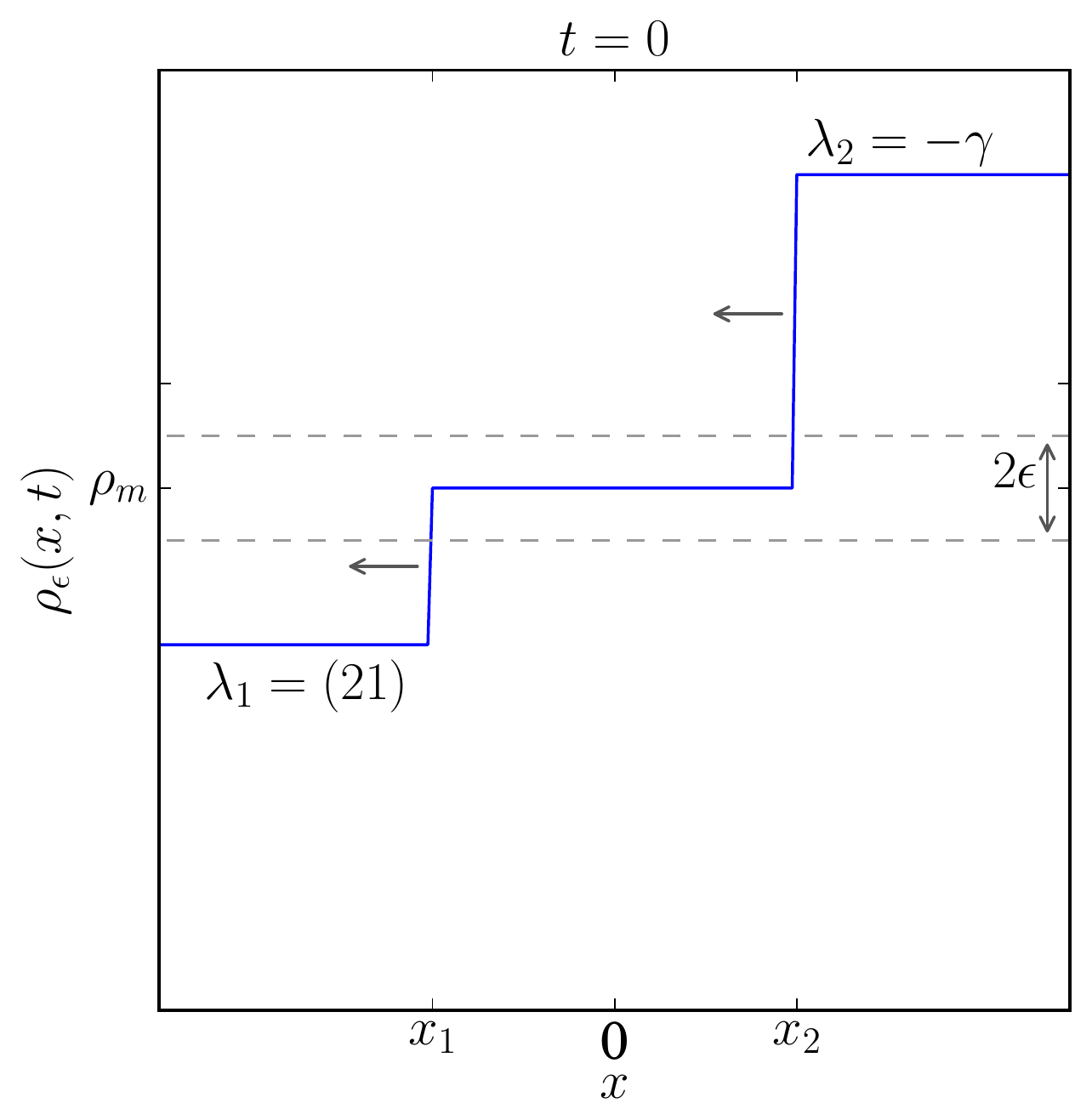}}
    \qquad
    \subfigure[]{\label{fig:ZeroRarefactionInteraction_Plots_Case3_Plot2}\includegraphics[width=0.37\textwidth]{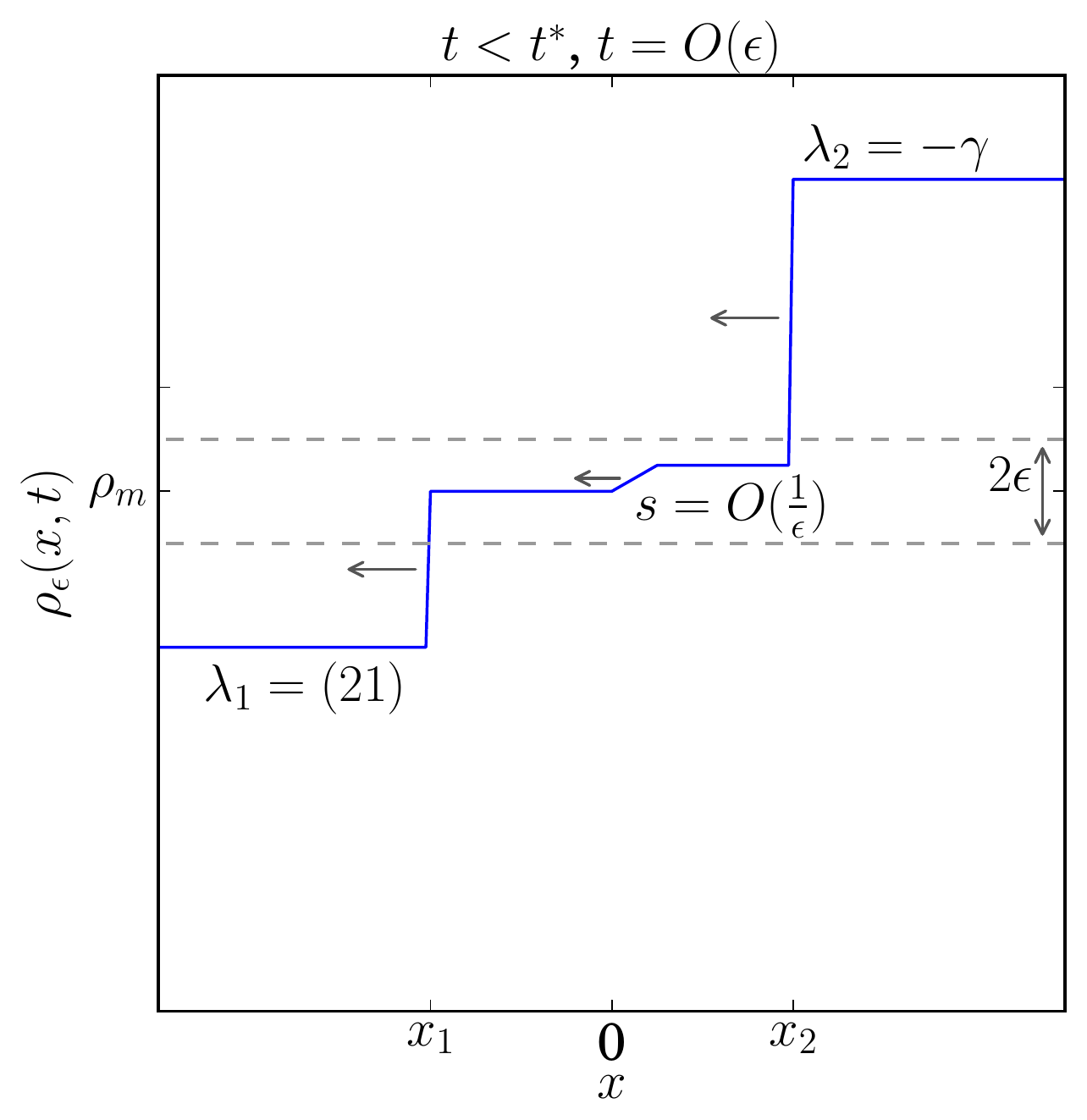}}\\
    \subfigure[]{\label{fig:ZeroRarefactionInteraction_Plots_Case3_Plot3}\includegraphics[width=0.37\textwidth]{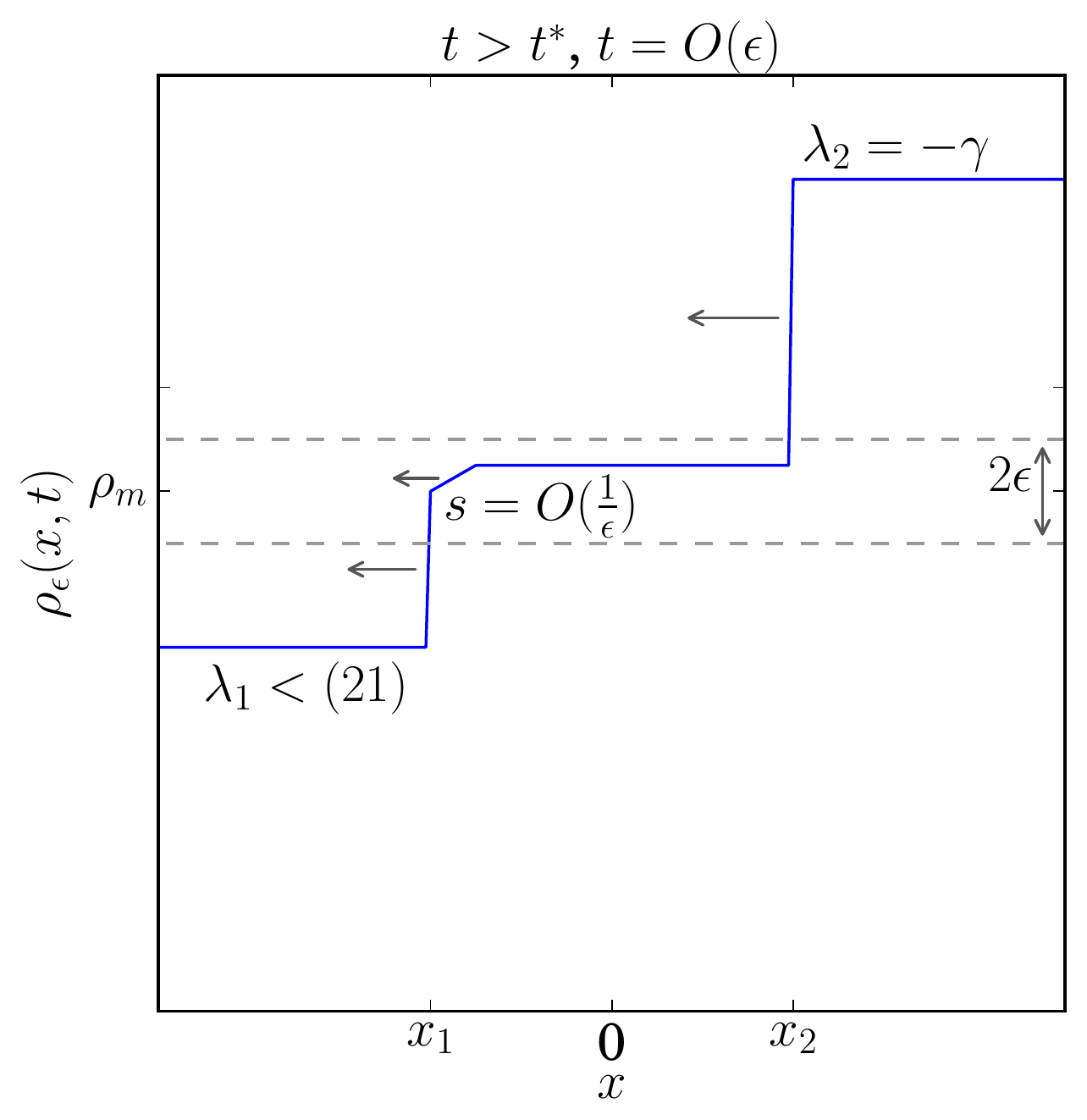}}
    \qquad
    \subfigure[]{\label{fig:ZeroRarefactionInteraction_Plots_Case3_Plot4}\includegraphics[width=0.37\textwidth]{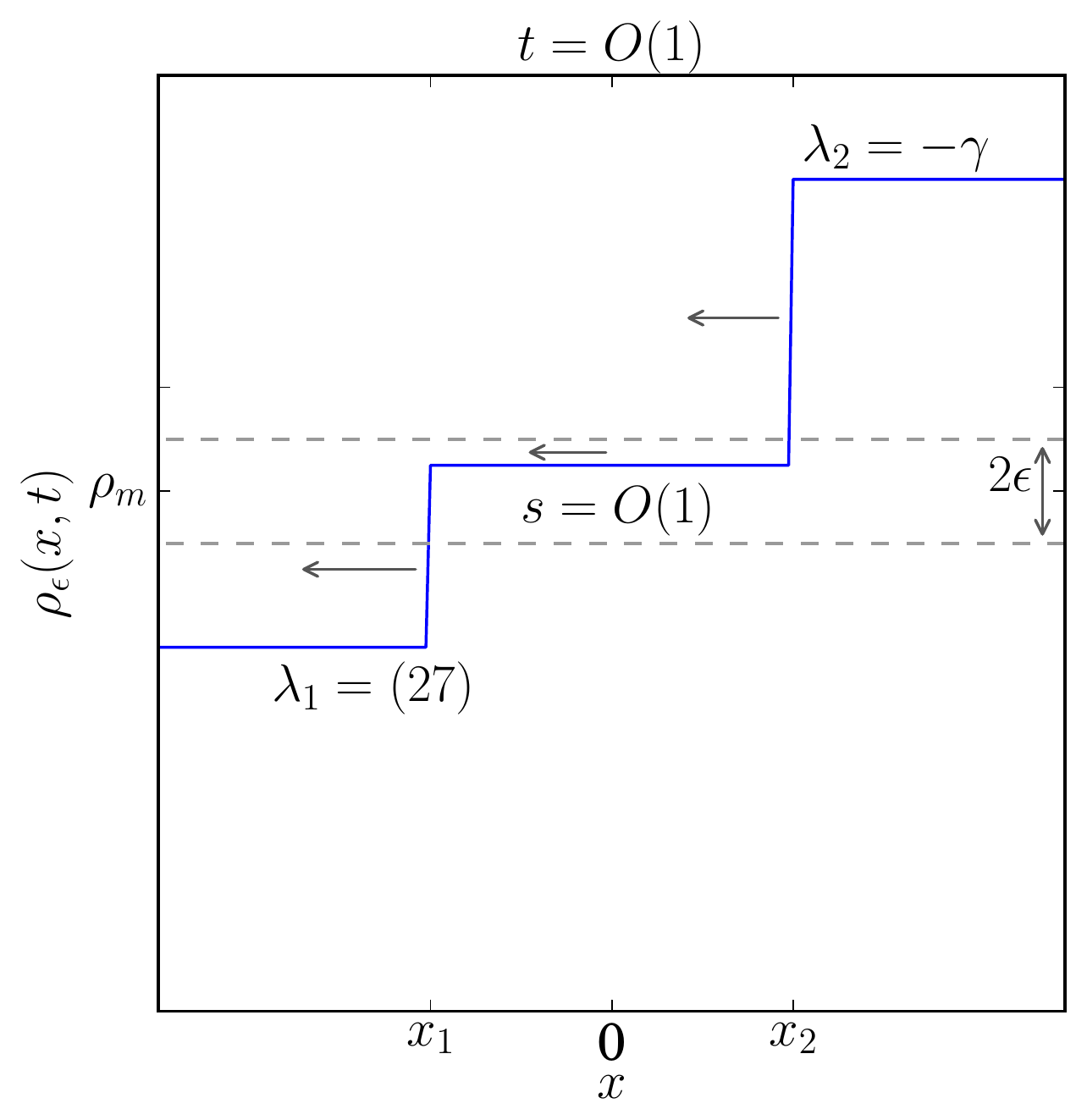}}
    \caption{Time evolution of the double Riemann solution in
      Case~2a when $C_r>\rho_m$ and $\gamma/(\gamma+1) \leqslant C_l < \rho_m$.}
    \label{fig:ZeroRarefactionInteraction_Plots_Case3}
  \end{center}
\end{figure}

Next, we consider the double Riemann problem when $C_r>\rho_m$ and $C_l
< \rho_m$ which generates a {zero rarefaction wave} and contact line at
$x_2$, and a shock at $x_1$. Before the 1-shock and {zero rarefaction
  wave} interact at time $t^*=O(\epsilon)$, the shock speed $\lambda_1$
satisfies the Rankine-Hugoniot condition
\eqref{eqn:BeforeTStarLambda1ShockSpeed}.  For $t>t^*$, the shock and
{zero rarefaction} interact, thereby causing the value of
$\rho_\epsilon$ to the right of the 1-wave to increase.  Using a similar
argument as in Case~1, we can deduce that when $t=O(1)$,
\bq
  \lambda_1 = \frac{g_c(\rho_m) - f(C_l)}{\rho_m - C_l} 
  \qquad \text{and} \qquad 
  \lambda_2 = -\gamma,
\eq
as $\epsilon \rightarrow 0$.  Note that these wave speeds
are consistent with the Riemann problem in
Eq.~\eqref{eqn:Case2_ShockSpeed}.  The time evolution of the solution is
illustrated in Fig.~\ref{fig:ZeroRarefactionInteraction_Plots_Case3}.

There are actually two distinct sub-cases that need to be considered
here, corresponding to whether $C_l\geqslant\gamma/(\gamma+1)$ (which we
call Case~2a) or $C_l<\gamma/(\gamma+1)$ (Case~2b).  In Case~2a, the two
elementary waves (1--shock and 2--contact) do not interact, while in
Case~2b we have $\lambda_1 > \lambda_2$ and so the elementary waves
collide to form a single shock that has speed given by
Eq.~\eqref{eqn:Case3_ShockSpeed}.  This distinction is illustrated in
the characteristic plots in
Fig.~\ref{fig:ZeroRarefactionInteraction_PhasePlot_Case1-4}.
\\
\\
\noindent{\bf Case~3: $\boldsymbol{C_r < \rho_m < C_l}$.}

\begin{figure}[tbp]
  \begin{center}
    \subfigure[]{\label{fig:ZeroRarefactionInteraction_Plots_Case2_Plot1}\includegraphics[width=0.37\textwidth]{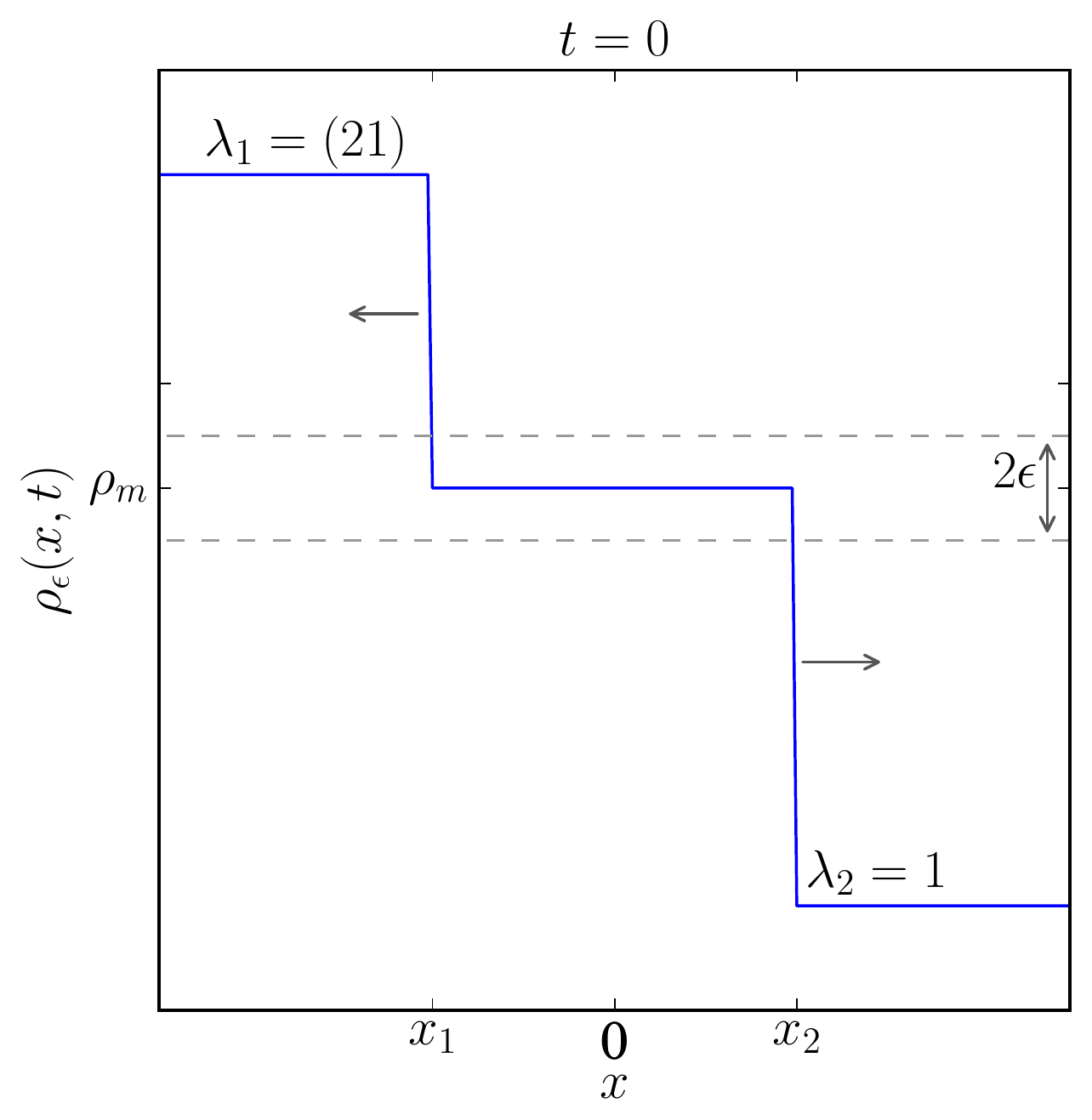}}
    \qquad
    \subfigure[]{\label{fig:ZeroRarefactionInteraction_Plots_Case2_Plot2}\includegraphics[width=0.37\textwidth]{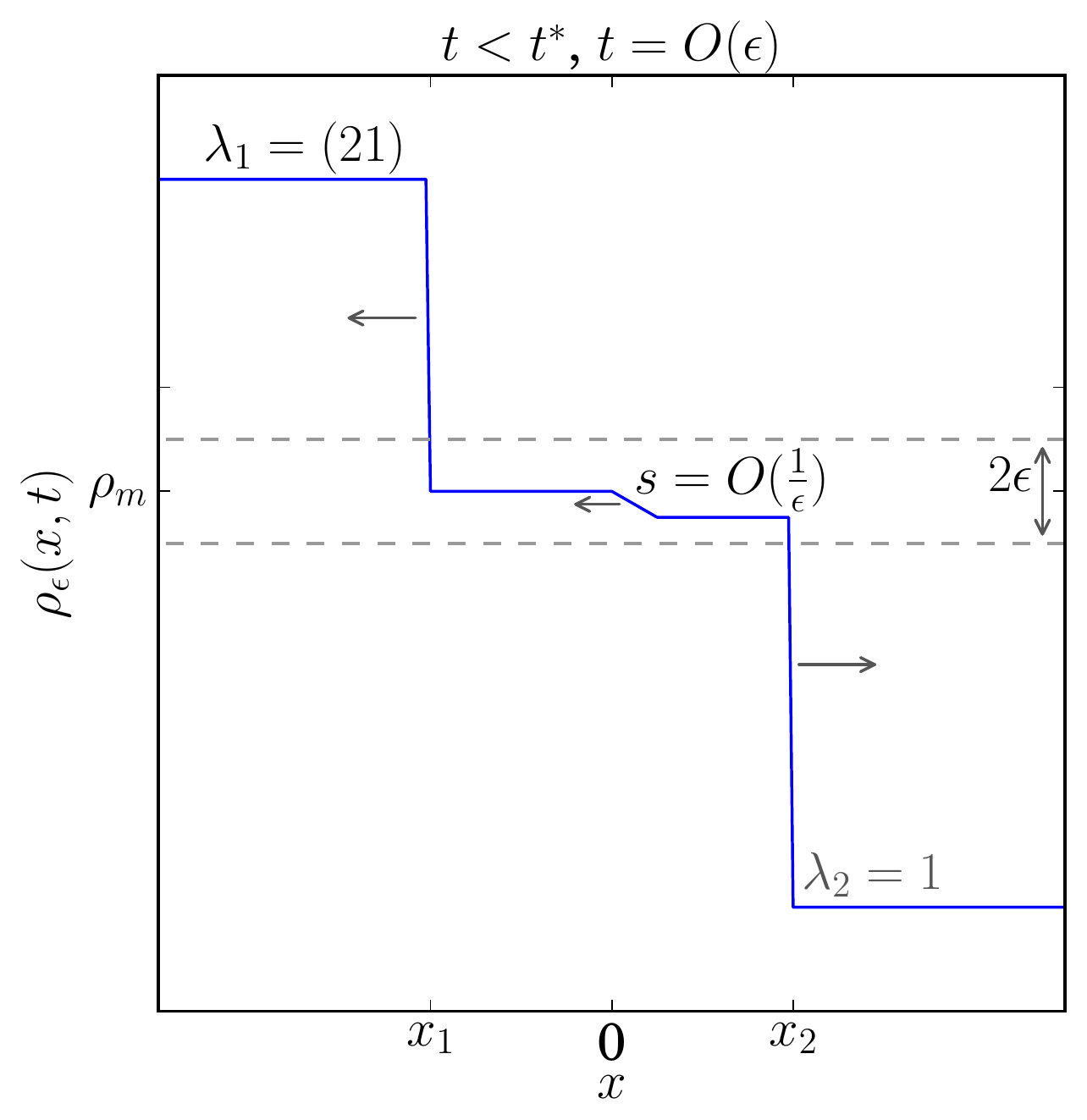}}\\
    \subfigure[]{\label{fig:ZeroRarefactionInteraction_Plots_Case2_Plot3}\includegraphics[width=0.37\textwidth]{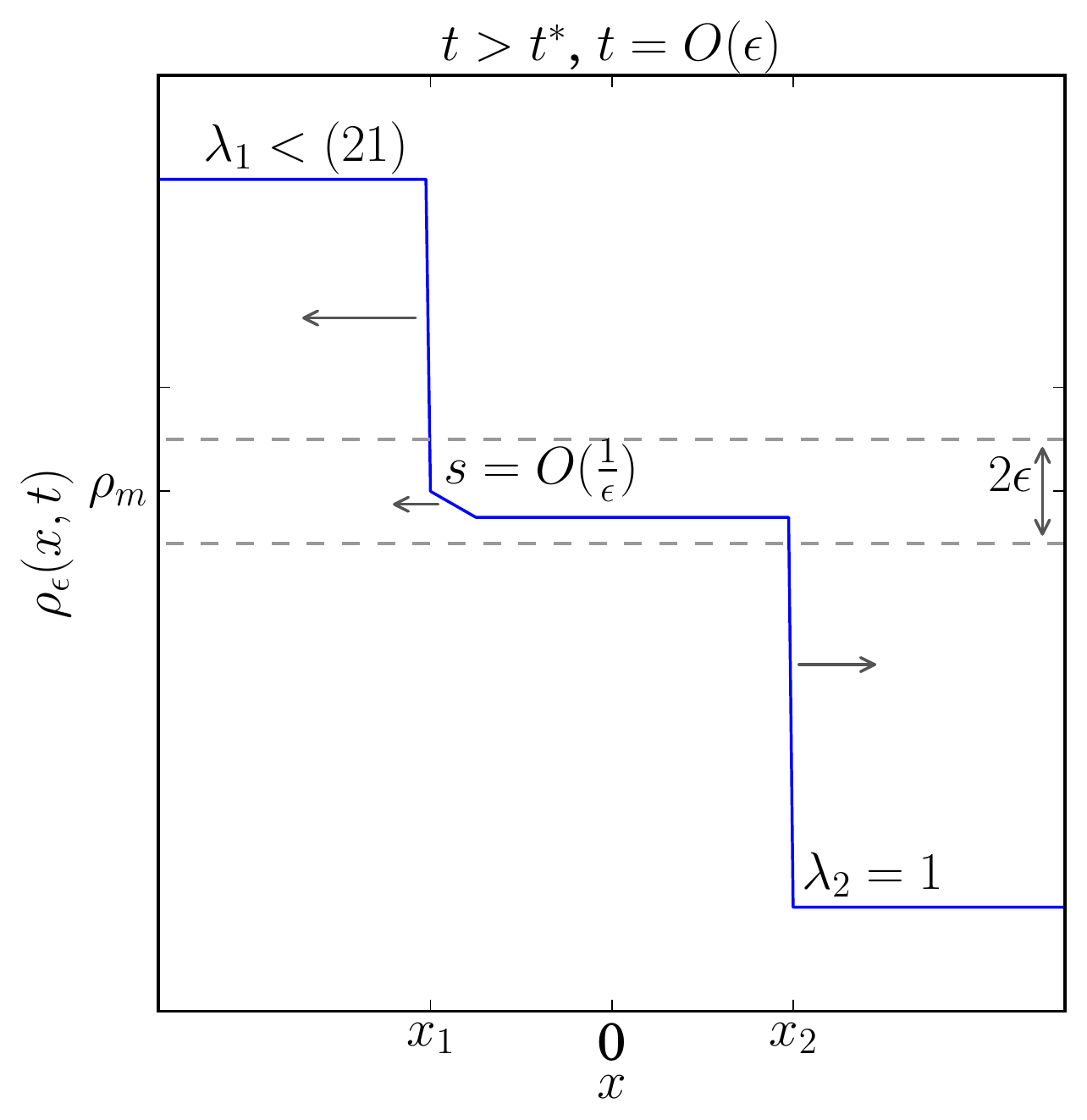}}
    \qquad
    \subfigure[]{\label{fig:ZeroRarefactionInteraction_Plots_Case2_Plot4}\includegraphics[width=0.37\textwidth]{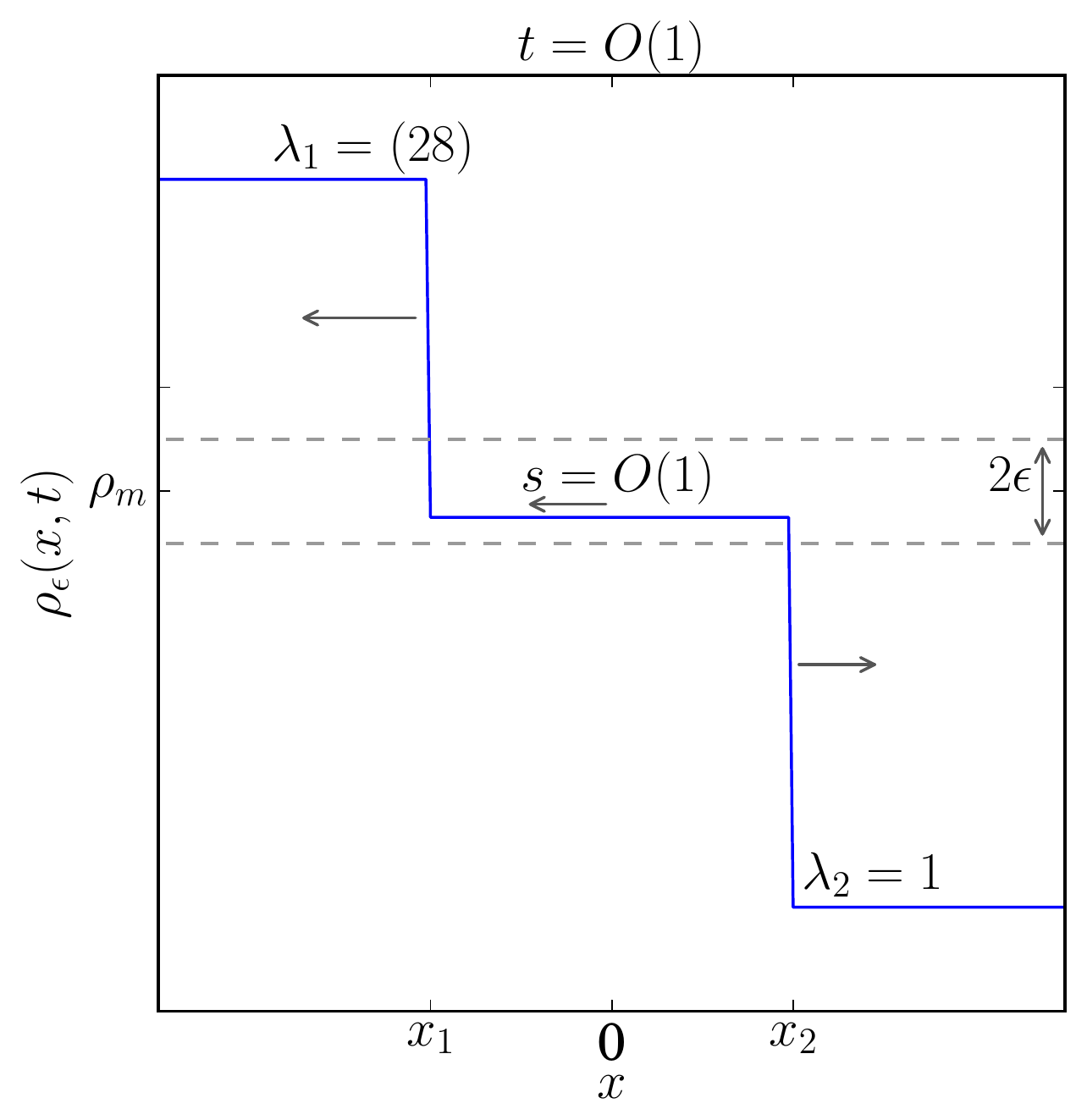}}
    \caption{Time evolution of the double Riemann solution in
      Case~3 when $C_r < \rho_m < C_l$.}
    \label{fig:ZeroRarefactionInteraction_Plots_Case2}
  \end{center}
\end{figure}

We next consider the situation when $C_r < \rho_m < C_l$, where a {zero
  rarefaction wave} and contact line are produced at $x_2$ and a shock
wave is generated at $x_1$.  Before the 1-shock and {zero rarefaction
  wave} interact at time $t^*=O(\epsilon)$, the shock speed $\lambda_1$
satisfies the Rankine-Hugoniot condition
\eqref{eqn:BeforeTStarLambda1ShockSpeed} as illustrated in
Figs.~\ref{fig:ZeroRarefactionInteraction_Plots_Case2}\subref{fig:ZeroRarefactionInteraction_Plots_Case2_Plot1}
and \subref{fig:ZeroRarefactionInteraction_Plots_Case2_Plot2}.  Once the
zero wave collides with the shock, the shock speed $\lambda_1$ increases
according to the Rankine-Hugoniot condition since the value of
$\rho_\epsilon$ to the right of the 1-wave decreases (shown in
Fig.~\ref{fig:ZeroRarefactionInteraction_Plots_Case2}\subref{fig:ZeroRarefactionInteraction_Plots_Case2_Plot3}). By
ignoring interactions on the $O(\epsilon)$ time scale, we find that
\bq
  \lambda_1 = \frac{g_f(\rho_m) - f(C_l)}{\rho_m - C_l},
  \qquad \text{and} \qquad 
  \lambda_2 = 1,
  \label{eqn:lambda_case3}
\eq
for $t=O(1)$ as $\epsilon \rightarrow 0$. Notice that the wave speeds
$\lambda_1$ and $\lambda_2$ are identical to those for the Riemann
problem considered in Case~C from Section~\ref{sec:RiemannProblem}.
\\
\\
\noindent{\bf Case~4: $\boldsymbol{C_r,C_l>\rho_m}$.}

The final remaining case corresponds to $C_r,C_l>\rho_m$ and leads to a
solution satisfying Eq.~\eqref{eqn:TripleRiemann_Ansatz} with
\bq
  \lambda_1 = \frac{g_c(\rho_m) - f(C_l)}{\rho_m - C_l}
  \qquad \text{and} \qquad 
  \lambda_2 = -\gamma.
  \label{eqn:lambda_case4}
\eq
The speed of the 1-shock is determined simply by the
Rankine-Hugoniot condition where $f(\rho_m)$ is evaluated on the
congested flow branch.  This should be compared to the 1-wave in Case~3
where the shock speed is determined by Eq.~\eqref{eqn:lambda_case3} in
which $f(\rho_m)$ is evaluated on the free flow branch.  No plots of the
solution have been provided for Case~4 since the structure is
essentially the same as in Case~1 except that all waves propagate in the
opposite direction.

%% file: Sec5.tex
\section{Finite Volume Scheme of Godunov Type}
\label{sec:NumericalScheme}

In this section, we construct a finite volume scheme of Godunov type for
our discontinuous flux problem based on ideas originally developed by
Godunov~\cite{Godunov1959} for solving the Euler equations of gas
dynamics.  In Godunov's method, the spatial domain is divided into cells
$[x_{j-1/2}, x_{j+1/2}]$ of constant width, $\Delta
x=x_{j+1/2}-x_{j-1/2}$, and the solution is assumed to be piecewise
constant on each grid cell.  Cell-averaged solution values $Q_j^n$
(located at cell centers) are updated using the exact solution from
local Riemann problems evaluated at interfaces between adjacent cells.
In each time step, the following three-stage algorithm, referred to as
the \emph{Reconstruct--Evolve--Average} or \emph{REA algorithm} in
\cite{LeVequeRedBook}, is used to update the $Q_j^n$:
\begin{itemize}
\item \emph{Reconstruct} a piecewise constant function $\tilde{\rho}(x,
  t_n) = Q_j^n$ for all $x \in [x_{j+1/2}, x_{j-1/2}]$ from the cell
  average $Q_j^n$ at time $t_n$.
\item \emph{Evolve} the conservation law exactly using initial
  data $\tilde{\rho}(x, t_n)$, thereby obtaining $\tilde{\rho}(x,
  t_{n+1})$ at time $t_{n+1}=t_n+\Delta t$.
\item \emph{Average} the solution $\tilde{\rho}(x, t_{n+1})$ to
  obtain new cell average values
  \bqs 
    Q_j^{n+1} = \frac{1}{\Delta x} \int_{x_{j-1/2}}^{x_{j+1/2}}
      \tilde{\rho}(x, t_{n+1}) \, dx \mbox{.} 
  \eqs
\end{itemize}
For problems with smooth flux, the evolution step can be performed by
solving a local Riemann problem at each cell interface having left state
$\rho_l = Q_j^n$ and right state $\rho_r = Q_{j+1}^n$.  As long as the
time step is chosen small enough, the elementary waves produced at each
interface do not interact (remembering that waves travel at finite speed
in the smooth case) and hence the evolution step yields an appropriate
approximation of the solution. However, as we have already shown, when
the flux is discontinuous the presence of {zero waves} travelling at
infinite speed gives rise to long-range interactions between local
Riemann problems.  Consequently, we make use of solutions to the {double
  Riemann solution} derived in the previous section that incorporate the
effects of zero waves.  We note that the resulting algorithm has some
similarities to the method of Gimse \cite{Gimse1993}.

We next provide details of our implementation using LeVeque's high
resolution wave propagation formulation \cite{LeVequeRedBook}, in which
the Riemann solver returns a set of wave strengths
$\mathcal{W}^p_{j+1/2}$ and speeds $s^p_{j+1/2}$ generated at each
interface between states $Q_j$ and $Q_{j+1}$ (in what follows, we will
omit the superscript $n$ when it is clear that time index $n$ is
assumed).  The \emph{evolution} step of the algorithm above can then be
written as
\leavethisout{
  \begin{gather}
    Q_j^{n+1} = Q_j^n - \frac{\Delta t}{\Delta x} \left( F_{j+1/2}^n -
      F_{j-1/2}^n \right), 
    \label{eqn:evolve-step}
  \end{gather}
  where the numerical flux function is $F^n_{j-1/2} =
  f\left(q^\ast\left(Q_{j-1}^n,Q_j^n\right)\right)$ and
  $q^\ast(Q_L,Q_R)$ is the exact solution to the Riemann problem joining
  constant states $Q_L$ and $Q_R$.  
}
\bq
  Q_j^{n+1} = Q_j^n - \frac{\Delta t}{\Delta x} \left[ 
    \sum_p \left( s_{j-1/2}^p \right)^+ \mathcal{W}_{j-1/2}^p + 
    \sum_p \left( s_{j+1/2}^p \right)^- \mathcal{W}_{j+1/2}^p 
    \right],
  \label{eqn:evolve-step}
\eq
where $(s)^+=\max(s,0)$ and $(s)^-=\min(s,0)$.  When the flux is smooth,
the wave speed $s_{j+1/2}$ depends only on the states $Q_j$ and
$Q_{j+1}$, whereas for a discontinuous flux function this is no longer
the case.  When $Q_j,Q_{j+1} \ne \rho_m$, the Riemann solver yields the
``standard'' elementary waves whose strength and speed are given in
Section \ref{sec:RiemannProblem}.  Because our flux is strictly
non-convex, we also observe compound waves which are described within
Cases 1 and 2 in Section \ref{sec:RiemannProblem}. 

Our Riemann problem solution diverges from the standard one when either
$Q_j = \rho_m$ or $Q_{j+1} = \rho_m$, in which case we construct the
solution to a local double Riemann problem that requires the appropriate
1- or 2-wave given in Section~\ref{sec:ZeroWaves}.  Note that the double
Riemann solution should be viewed as {two} separate local Riemann
problems that each produce {one} elementary wave: the left Riemann
problem corresponds to the 1-wave and the right Riemann problem
corresponds to the 2-wave.

For example, if $Q_j = \rho_m$, then we construct the double
Riemann problem with $Q_j = \rho_m$ and $Q_{j+1} = C_r$, thereby
obtaining the strength and wave speed corresponding to the 2-wave in
Section~\ref{sec:ZeroWaves}.  Combining together all four cases in
Section~\ref{sec:ZeroRarefactionInteraction}, the speed can be written
as
\bq
    s^1_{j+1/2} = \left\{
	\begin{array}{ll}
	  1      ~, & ~ Q_{j+1} < \rho_{m}\\
	  -\gamma~, & ~ Q_{j+1} > \rho_{m}
	\end{array}
    \right. \mbox{,}
\label{eqn:Lambda2_Shockspeed}
\eq
which we note is independent of $C_l$. 

\begin{figure}[tbp]
  \begin{center}
    \includegraphics[width=0.8\textwidth]{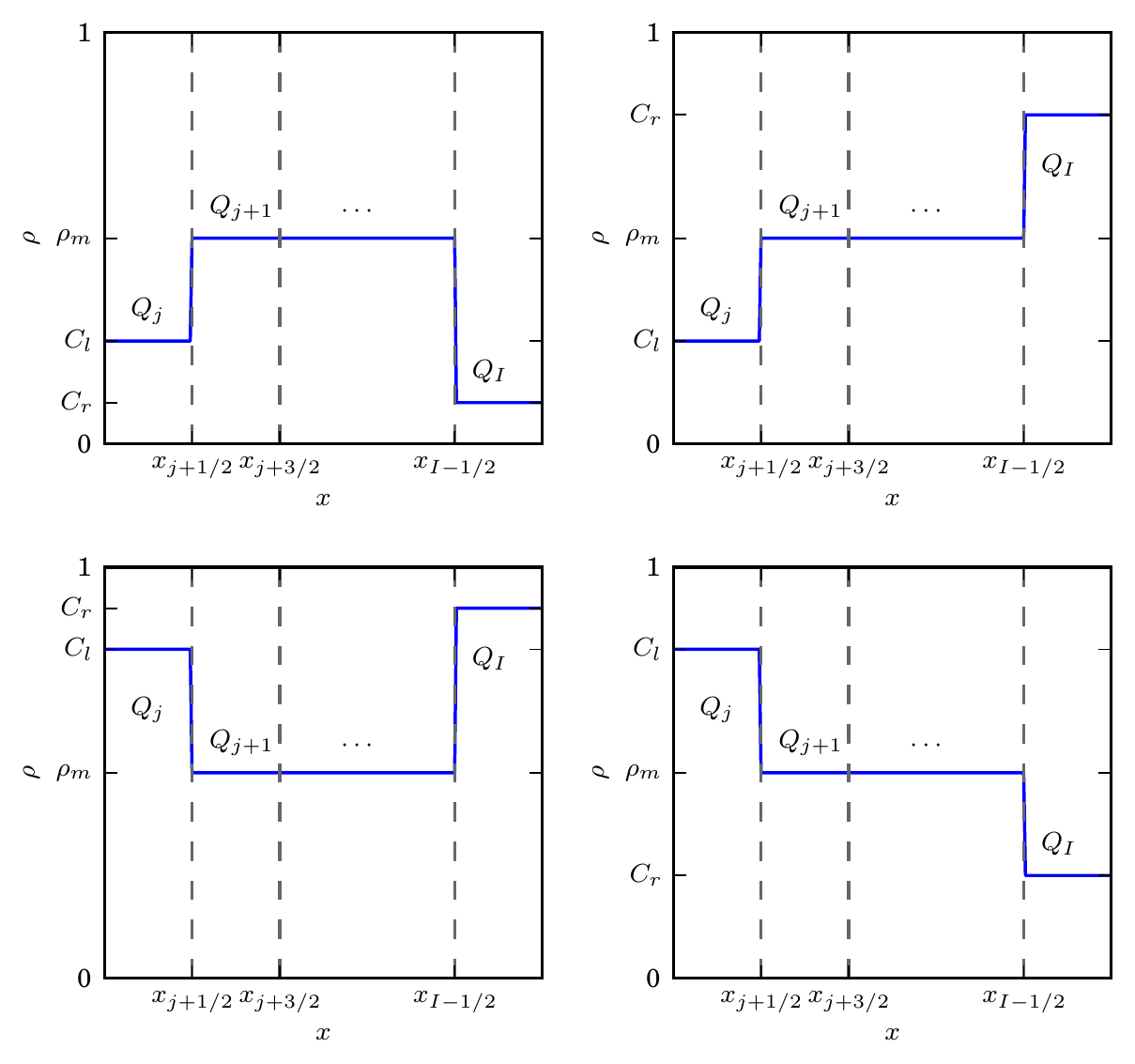}
    \caption{The four possible cases that can arise when calculating the
      shock speed in Eq.~\eqref{eqn:Lambda1_Shockspeed}.}
    \label{fig:Lambda1_Shockspeed}
  \end{center}
\end{figure}

Alternatively, if $Q_{j+1} = \rho_m$ then we construct the double
Riemann problem with $Q_{j} = C_l$, $Q_{j+1} = \rho_m$ and unknown $Q_I
= C_r$, thereby obtaining the strength and wave speed for the 1-wave in
Section~\ref{sec:ZeroWaves}.  In contrast with the case $Q_j=\rho_m$
just considered, the 1-wave's speed depends on values of $C_l$,
$\rho_m$, and $C_r$.  Therefore, when determining the speed of the wave
at the interface between $Q_{j} = C_l$ and $Q_{j+1} = \rho_m$, we must
\emph{look ahead} to find the value $Q_I = C_r$ corresponding to the
first value of $Q_j$ not equal to $\rho_m$; that is,
\bq
    I = \min \left \{~k~|~j+1 < k \leqslant N \mbox{ and } Q_{k} \ne \rho_{m} 
    \right \}. 
    \label{eqn:LookAheadIndex}
\eq
In summary, when $Q^n_{j+1}=\rho_{m}$, the wave speed reduces to 
\bq
  s^1_{j+1/2} = \left\{
    \begin{array}{ll}
      \displaystyle \frac{g_f(\rho_m) - f(Q_j)}{\rho_m - Q_j} ~,& ~ Q_I < \rho_{m}\\
      \displaystyle \frac{g_c(\rho_m) - f(Q_j)}{\rho_m - Q_j} ~,& ~ Q_I > \rho_m
    \end{array}
  \right.
  \label{eqn:Lambda1_Shockspeed}
\eq
which is visualized in Fig.~\ref{fig:Lambda1_Shockspeed}. 

Note that we have not yet discussed the two simple cases when
$Q_j,Q_{j+1} < \rho_m$ and $Q_j,Q_{j+1} > \rho_m$, both of which reduce
to the linear advection equation and can be trivially solved.  A summary
of wave strengths and speeds for all cases is presented in
Table~\ref{Table:Wave:StrengthAndSpeed}.
\begin{table*}[htbp]\centering\small
  \ra{1.6}
  \caption{Wave speed and strength for the 1- and 2-waves in all possible cases.}
  \begin{tabular}{@{}l@{}l|cc|cc@{}}\toprule
    \multicolumn{2}{l|}{Case} & $s^1_{j+1/2}$ & $\mathcal{W}^1_{j+1/2}$ & $s^2_{j+1/2}$ & $\mathcal{W}^2_{j+1/2}$  \\
    \midrule
    & $Q_{j} = Q_{j+1} = \rho_{m}$                & $0$ & $0$ & $0$ & $0$ \\
    & $Q_{j} = \rho_{m}$ and $Q_{j} \ne Q_{j+1}$  & Eq.~\eqref{eqn:Lambda2_Shockspeed} & $Q_{j+1}-Q_{j}$ & $0$ & $0$\\
    & $Q_{j+1} = \rho_{m}$ and $Q_{j} \ne Q_{j+1}$& Eq.~\eqref{eqn:Lambda1_Shockspeed} & $Q_{j+1}-Q_{j}$ & $0$ & $0$  \\
    & $\rho_{m} < Q_{j}$ and $\rho_{m} < Q_{j+1}$ & $-\gamma$ & $Q_{j+1}-Q_{j}$ & $0$ & $0$\\
    & $\rho_{m} > Q_{j}$ and $\rho_{m} > Q_{j+1}$ & $1$ & $Q_{j+1}-Q_{j}$ & $0$ & $0$\\
    & $Q_{j} < \rho_{m} < Q_{j+1}$ and $\displaystyle \frac{\gamma}{\gamma+1} \geqslant Q_{j}$ & $\displaystyle\frac{f(Q_{j+1})-f(Q_{j})}{Q_{j+1} - Q_{j}}$ & $Q_{j+1}-Q_{j}$ & $0$ & $0$ \\
    & ${\displaystyle \frac{\gamma}{\gamma+1}} < Q_{j} < \rho_{m} < Q_{j+1}$ & $-\gamma$ & $Q_{j+1}-\rho_{m}$ & $\displaystyle\frac{g_c(\rho_{m}) - f(Q_{j})}{\rho_{m}-Q_{j}}$ & $\rho_{m}-Q_{j}$\\
    & $Q_{j+1} < \rho_{m} < Q_{j}$                & $\displaystyle \frac{f(Q_{j})-g_f(\rho_{m})}{Q_{j} - \rho_{m}}$ & $\rho_{m}-Q_{j}$ & $1$ & $Q_{j+1} - \rho_{m}$ \\
    \bottomrule
  \end{tabular}
  \label{Table:Wave:StrengthAndSpeed}
\end{table*}

A slight modification to the local Riemann solver is required to deal
with the fact that algebraic operations are actually performed in
floating-point arithmetic.  It is highly unlikely that the numerical
value of $Q_{j}$ ever exactly equals $\rho_m$, and yet we find that it
is necessary to employ the solution of the double Riemann problem when
$Q_{j}$ is \emph{close to} $\rho_m$.  Therefore, we need to relax the
requirement slightly for the zero-wave cases by replacing the condition
$Q_{j} = \rho_m$ with
\bq 
  | Q_{j} - \rho_m | \leqslant \delta
  \mbox{,} \label{eqn:FloatingPointFix}  
\eq
where $\delta$ is a small parameter that is typically assigned values on
the order of $10^{-5}$.  The choice of $\delta$ is a balance between
accuracy and efficiency in that taking a larger $\delta$ value 
gives rise to significant deviations in the height of the $\rho_m$
plateau regions, and hence also errors in mass conservation.  Taking
$\delta$ values much smaller than $10^{-5}$ does not improve the
solution significantly but does require a smaller time step for
stability.  This modification influences the accuracy of the simulations
and also creates an artificial upper bound on the maximum wave speed,
not including {zero waves}. For example, consider the shock solution in
Eq.~\eqref{eqn:Case2_ShockSpeed}, where as $\rho_l$ approaches $\rho_m$
the shock speed becomes unbounded.  By enforcing
\eqref{eqn:FloatingPointFix}, the parameter $\delta$ determines how
close $\rho_l$ can be to $\rho_m$ before the algorithm switches to the
double Riemann solution.

The time step $\Delta t$ is chosen adaptively to enforce stability of
our explicit update scheme, using a restriction based on the wave speeds
from all local Riemann problems.  In particular, we take
\bqs
  \Delta t = \mbox{{CFL}} \cdot \min_{j,p} \left(\frac{\Delta
      x}{|s^p_{j+1/2}|}\right), 
\eqs
where $0<\mbox{{CFL}}<1$ is a constant chosen to be around 0.9 in
practice.  As long as the parameter $\delta$ is not taken too small,
this condition is sufficient to guarantee stability.  Because the
effects of the zero waves have been incorporated directly into the
Riemann problem, they have no direct influence on the stability
restriction.

The Riemann solver described above forms the basis for the first-order
Godunov scheme.  We have also implemented a high resolution variant
using \emph{wave limiters} which limit the waves $\mathcal{W}^p_{j+1/2}$
in a manner similar to the limiting of fluxes in flux-based finite
volume schemes.  The details of this modification are described
in~\cite{LeVequeRedBook}.

\section{Numerical Results}
\label{sec:NumericalResults}

\begin{figure}[tbp]
  \centering
  \includegraphics[width=0.6\textwidth]{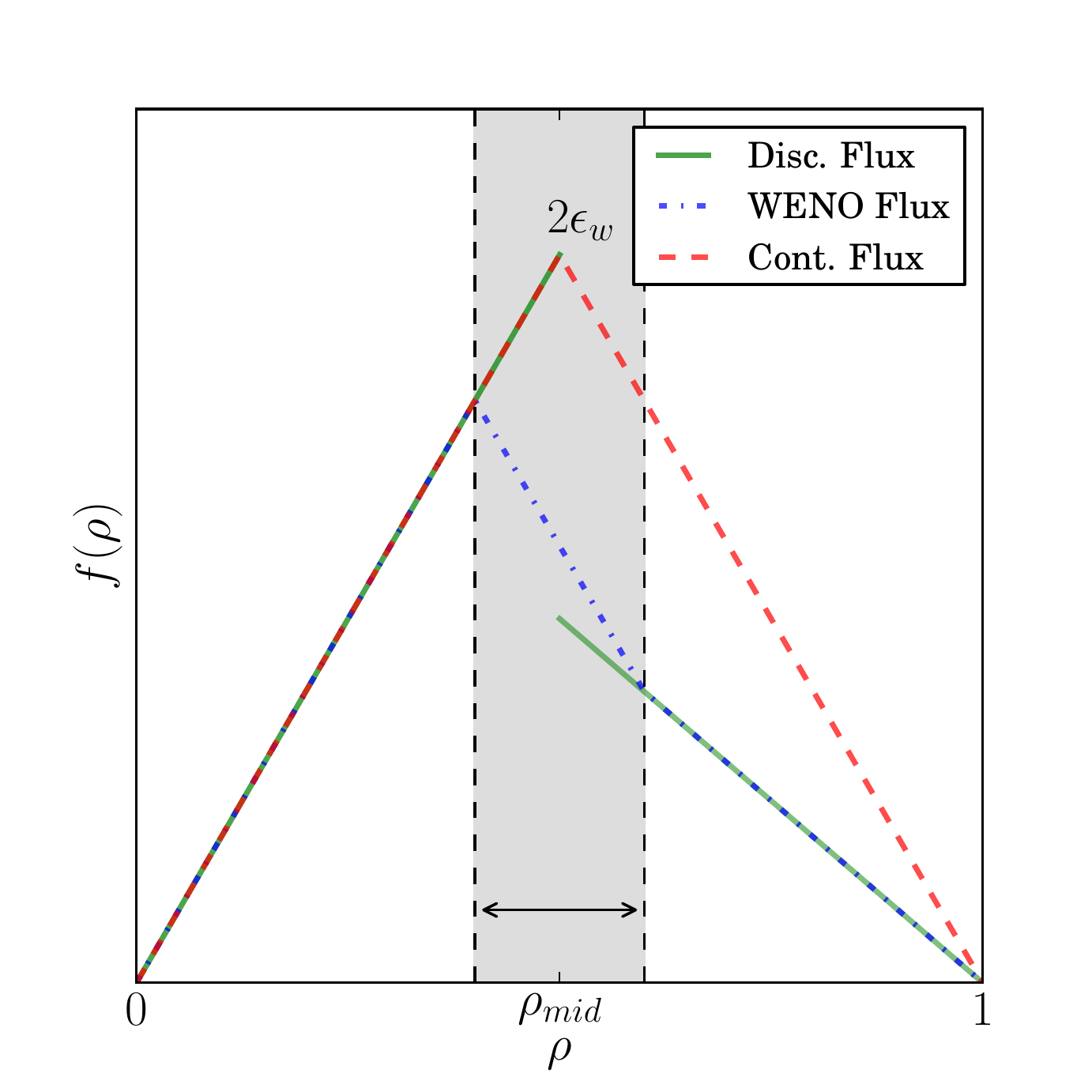}
  \caption{Plots of the discontinuous flux function used in our high
    resolution simulations (with $\rho_m=0.5$, $\gamma=0.5$), the
    regularized flux for WENO calculations in
    Section~\ref{sec:results:smooth}, and the continuous flux used in
    Section~\ref{sec:results:continuous}.}
  \label{fig:myflux}
\end{figure}

We now apply the method described in the previous section to a number of
test problems.  For the high resolution scheme, we employ the superbee
limiter function.   In all cases, we use the discontinuous, piecewise
linear flux function \eqref{eqn:DiscConsLaw}--\eqref{eqn:DiscFlux} with
parameters $\rho_m=0.5$ and $\gamma=0.5$ that is pictured in
Fig.~\ref{fig:myflux}.  

\subsection{Riemann Initial Data}
\label{sec:results:riemann}

As a first illustration of our numerical method, we use three sets of
Riemann initial data:
\begin{enumerate}
  \renewcommand{\theenumi}{\alph{enumi}}
  \renewcommand{\labelenumi}{(\theenumi)}
  \renewcommand{\theenumi}{\Alph{enumi}}
  \renewcommand{\labelenumi}{\theenumi.}
\item $\rho_l = 0.9$ and $\rho_r=0.2$,
\item $\rho_l = 0.4$ and $\rho_r=0.9$,
\item $\rho_l = 0.3$ and $\rho_r=0.98$,
\end{enumerate}
for which the exact solution can be determined using the methods in
Sections \ref{sec:RiemannProblem} and \ref{sec:ZeroWaves}.  These
initial data labeled A, B and C correspond to the three ``Cases'' with
the same labels in Section~\ref{sec:RiemannProblem} and pictured in
Figs.~\ref{fig:Case1_RiemannSolution}, \ref{fig:Case2_RiemannSolution}
and \ref{fig:Case3_RiemannSolution} respectively.  We note that the
numerical values for the left and right states in Case~C differ slightly
from the initial data used in Fig.~\ref{fig:Case3_RiemannSolution} in
order to generate larger wave speeds.

\begin{figure}[!tbp]
  \begin{center}
    \subfigure[]{\label{fig:RP_Sim_Case1}\includegraphics[width=0.40\textwidth]{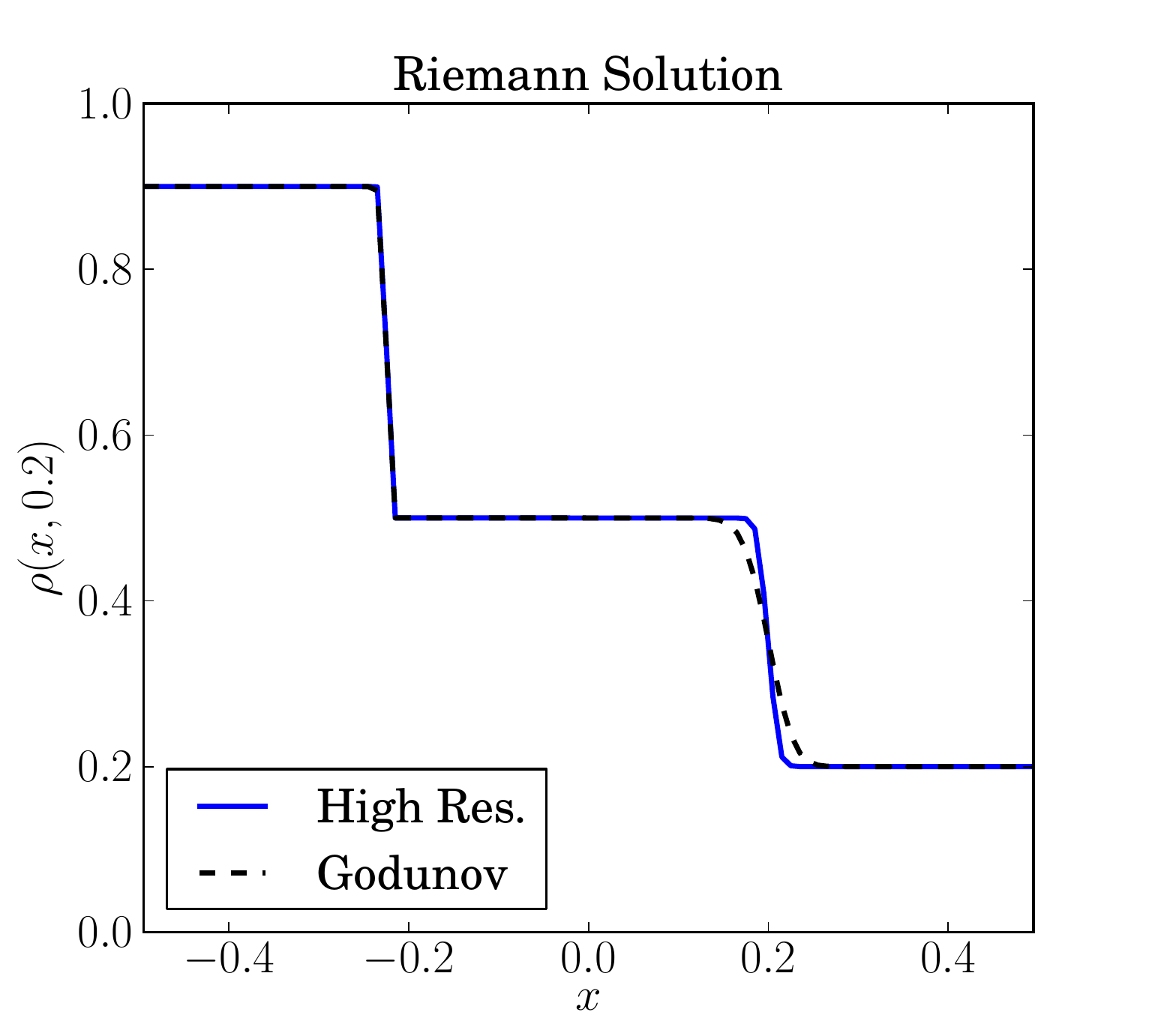} 
      \qquad
      \includegraphics[width=0.40\textwidth]{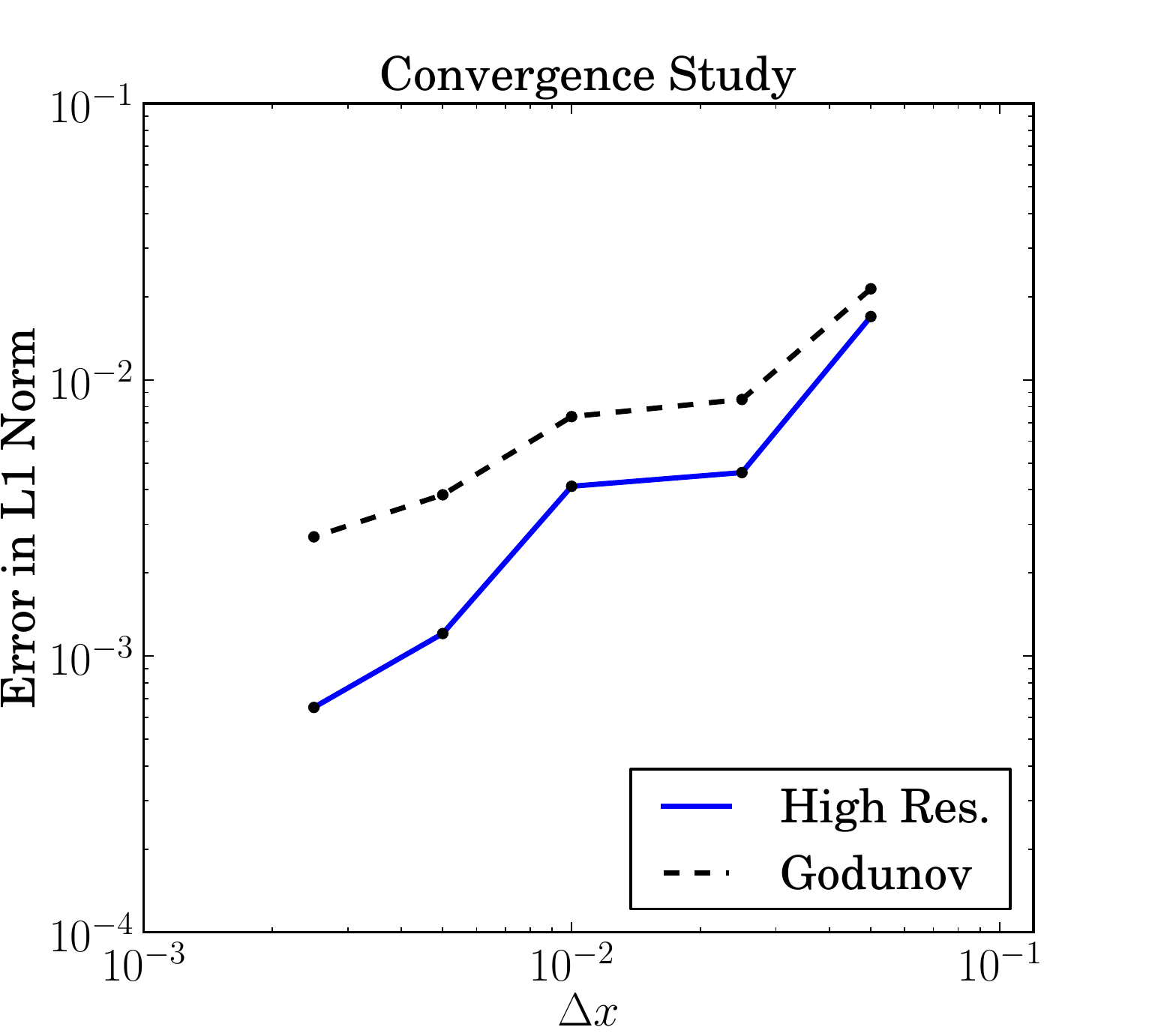}}
    \subfigure[]{\label{fig:RP_Sim_Case2}\includegraphics[width=0.40\textwidth]{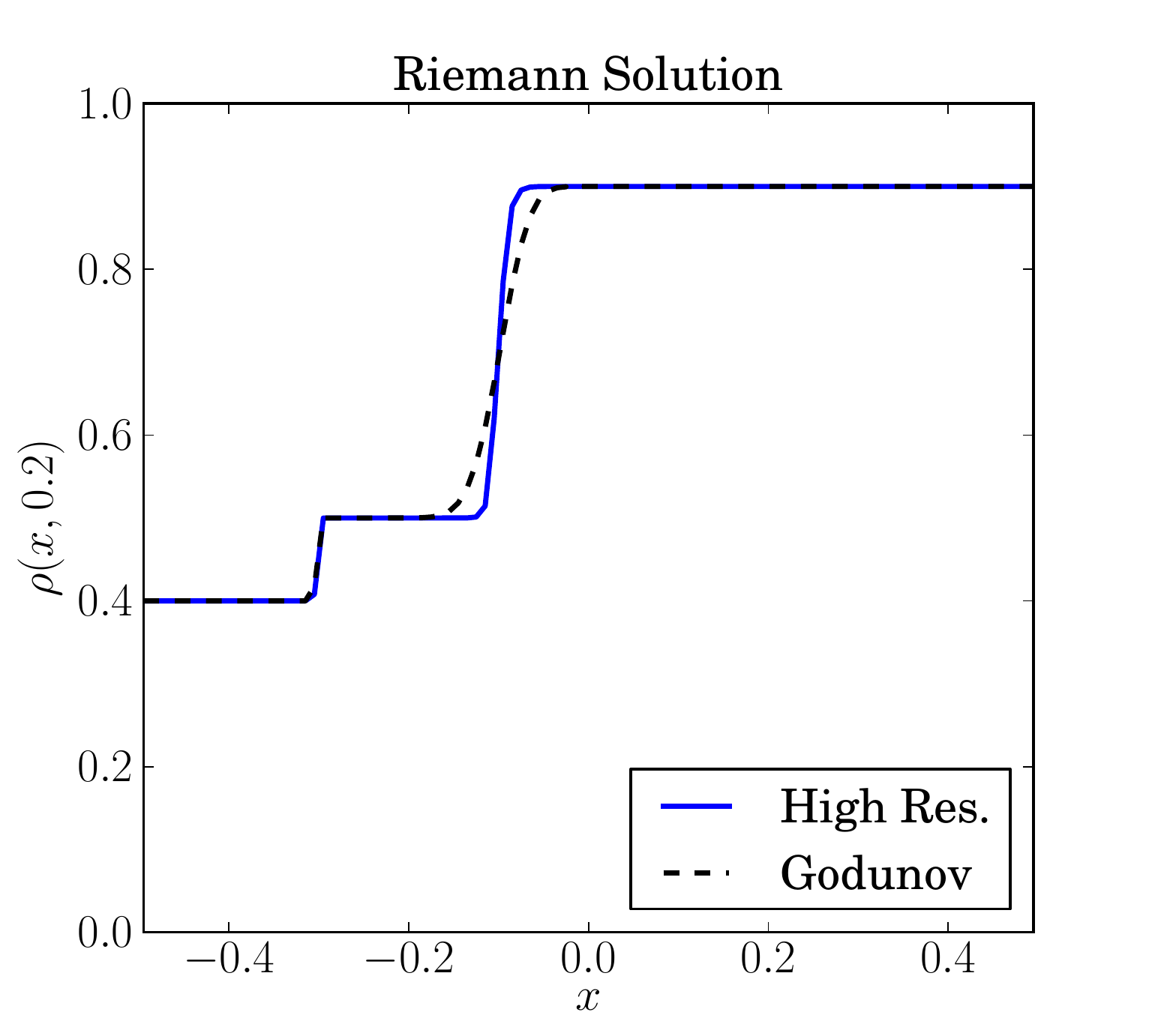} 
      \qquad
      \includegraphics[width=0.40\textwidth]{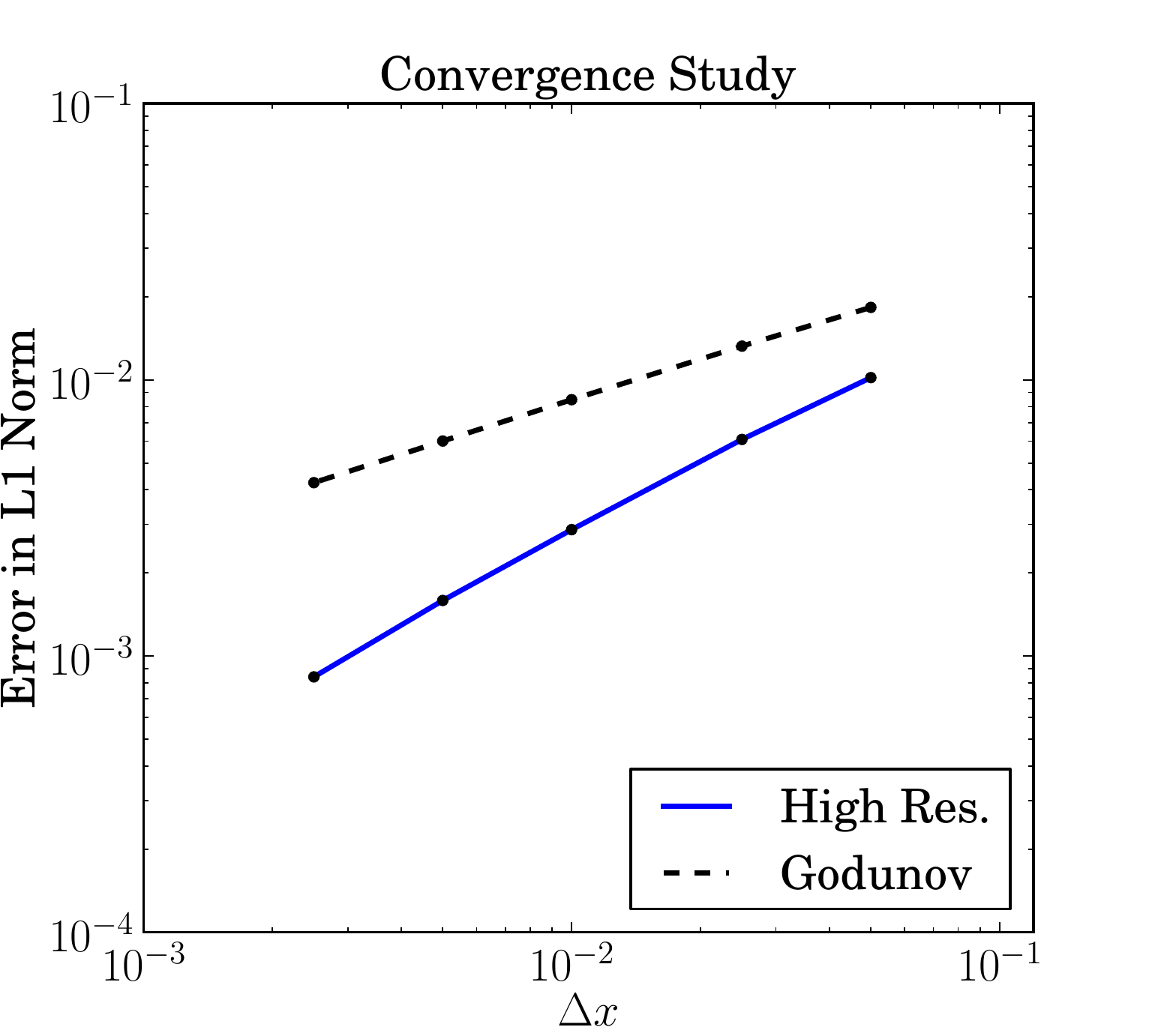}}
    \subfigure[]{\label{fig:RP_Sim_Case3}\includegraphics[width=0.40\textwidth]{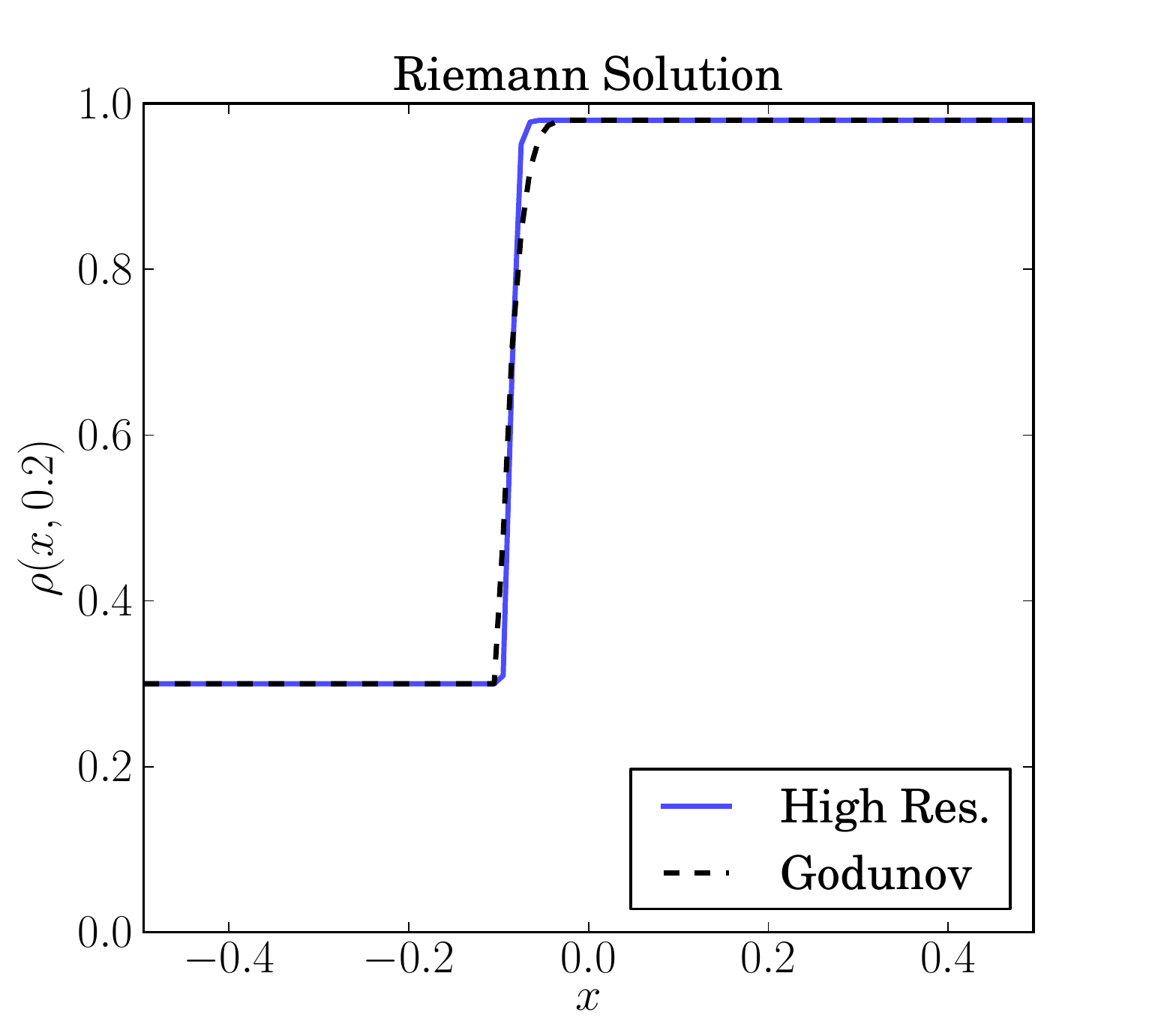} 
      \qquad
      \includegraphics[width=0.40\textwidth]{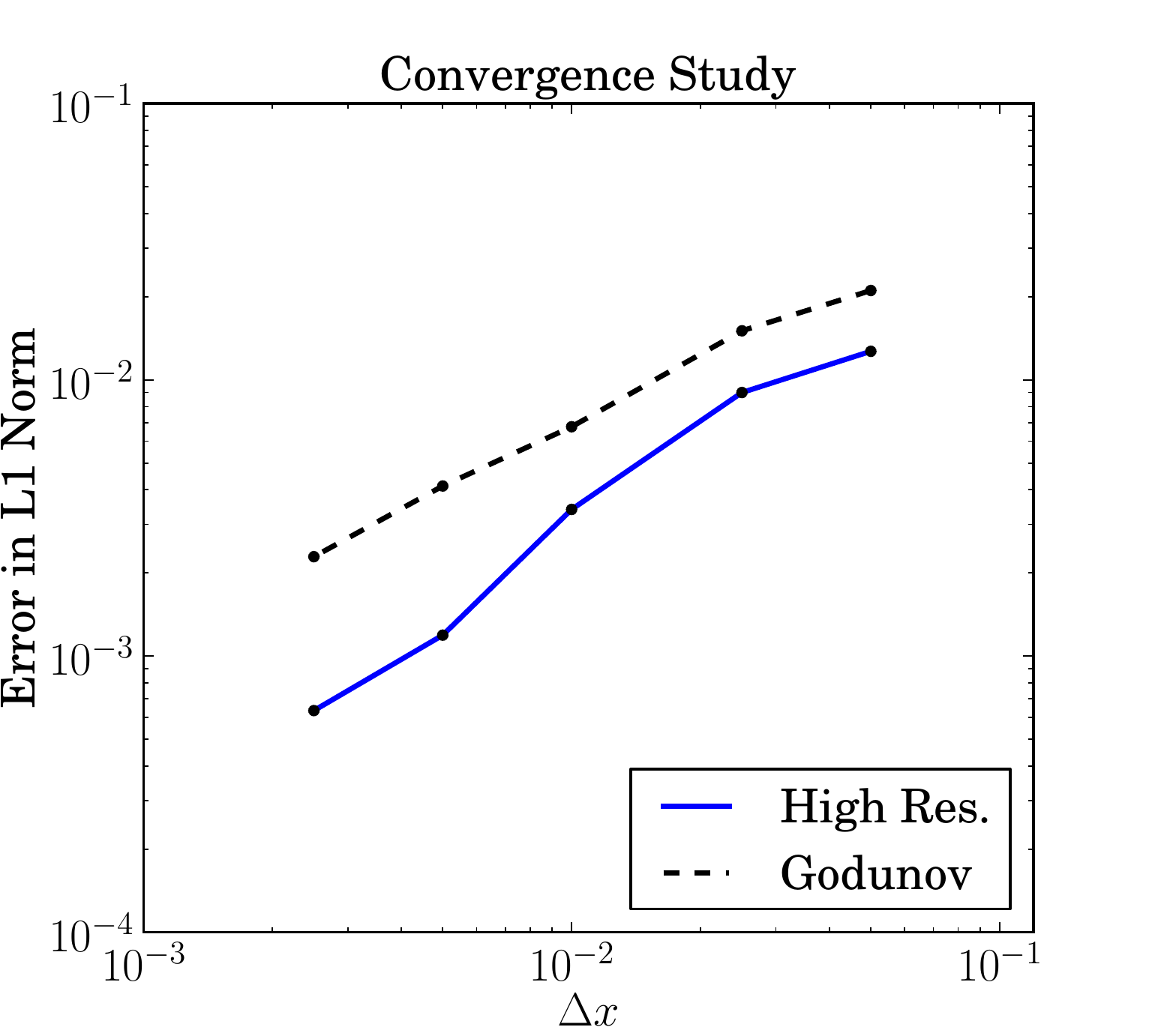}} 
    \caption{Left: Computed solutions for the Riemann initial data at
      $t=0.2$ where $\delta = 10^{-7}$, $\Delta x = 0.01$,
      $\mbox{CFL}=0.95$, $\rho_{m} = 0.5$, and $\gamma = 0.5$.  Right:
      Convergence study in the L1 norm.}
    \label{fig:RP_Sim}
  \end{center}
\end{figure}

In Fig.~\ref{fig:RP_Sim} (left), we compare the results from the first
order Godunov method and the high resolution scheme with wave limiters.
In each case, we perform a convergence study of error in the discrete L1
norm for grid resolutions ranging from $\Delta x=0.05$ to $0.0025$.
Convergence rates estimated using a linear least squares fit are
summarized in Table~\ref{Table:RP_ConvRates}.  As expected, the
numerical scheme converges to the exact solution for all three test
problems.\ \
Godunov's method captures the correct speed for both the shock and
contact discontinuities, although there is a more noticeable smearing of
the contact line which is typical for this type of problem.  The L1
convergence rates are consistent with the order $\sqrt{\Delta x}$
spatial error estimate established analytically for discontinuous
solutions of hyperbolic conservation laws having a smooth
flux~\cite{kuznetsov-1984,lucier-1985}.  The convergence rates in the L2
norm are also provided for comparison purposes and are significantly
smaller than the corresponding L1 rates, as expected.

\begin{table*}[ht!bp]\centering\small
  \caption{Convergence rates for the three Riemann problems
    pictured in Fig.~\ref{fig:RP_Sim}.}  
  \ra{1.6}
  \begin{tabular}{@{}lllcccc@{}}\toprule
    & & & \multicolumn{2}{c}{Godunov} & \multicolumn{2}{c}{High Resolution} \\
    \cmidrule(r){4-5} \cmidrule(r){6-7}
    $\rho_l$ & $\rho_r$ & \phantom{a} & L1 & L2 & L1 & L2 \\
    \midrule
    ${0.9}$ & ${0.2}$  & & $0.643$ & $0.367$ & $1.022$ & $0.569$  \\
    ${0.4}$ & ${0.9}$  & & $0.488$ & $0.232$ & $0.832$ & $0.375$  \\
    ${0.3}$ & ${0.98}$ & & $0.754$ & $0.373$ & $1.053$ & $0.627$  \\
    ${0.1}$ & ${0.4}$  & & $0.487$ & $0.145$ & $0.700$ & $0.238$  \\
    \bottomrule
  \end{tabular}
  \label{Table:RP_ConvRates}
\end{table*}

\subsection{Smooth Initial Data, With WENO Comparison}
\label{sec:results:smooth}

\begin{figure}[!tbp]
  \begin{center}
    \subfigure[]{\includegraphics[width=0.30\textwidth]{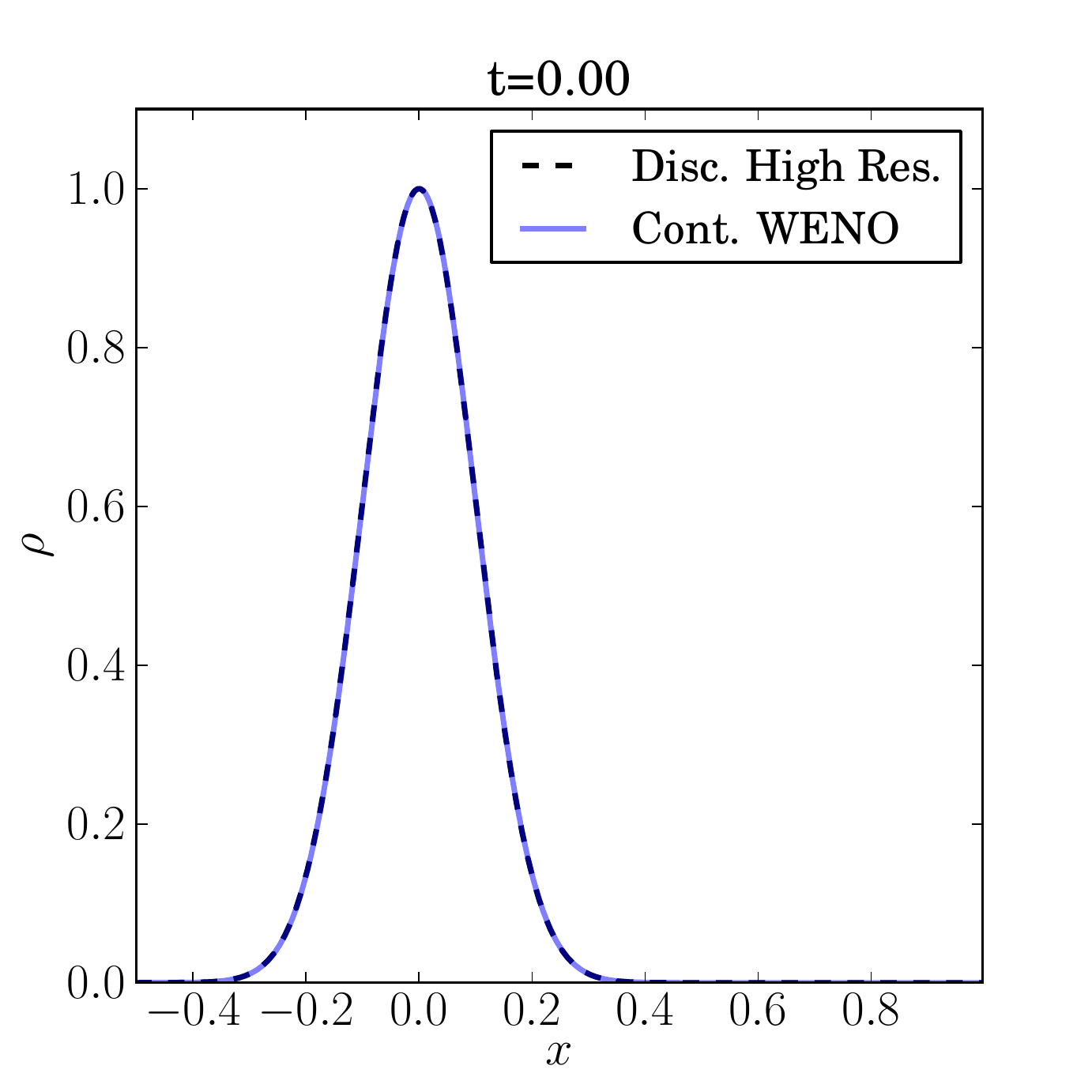}}
    \quad
    \subfigure[]{\includegraphics[width=0.30\textwidth]{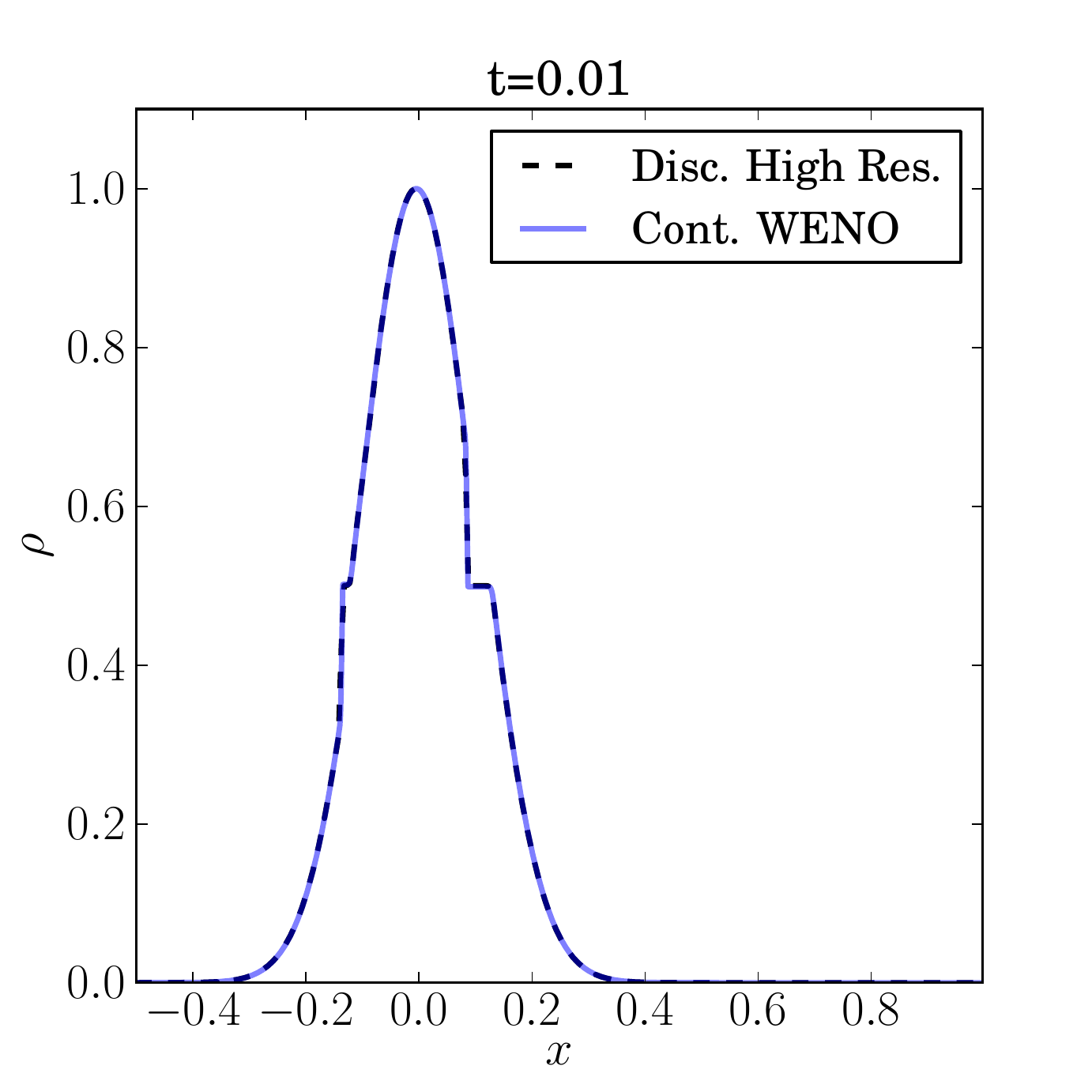}}
    \quad
    \subfigure[]{\includegraphics[width=0.30\textwidth]{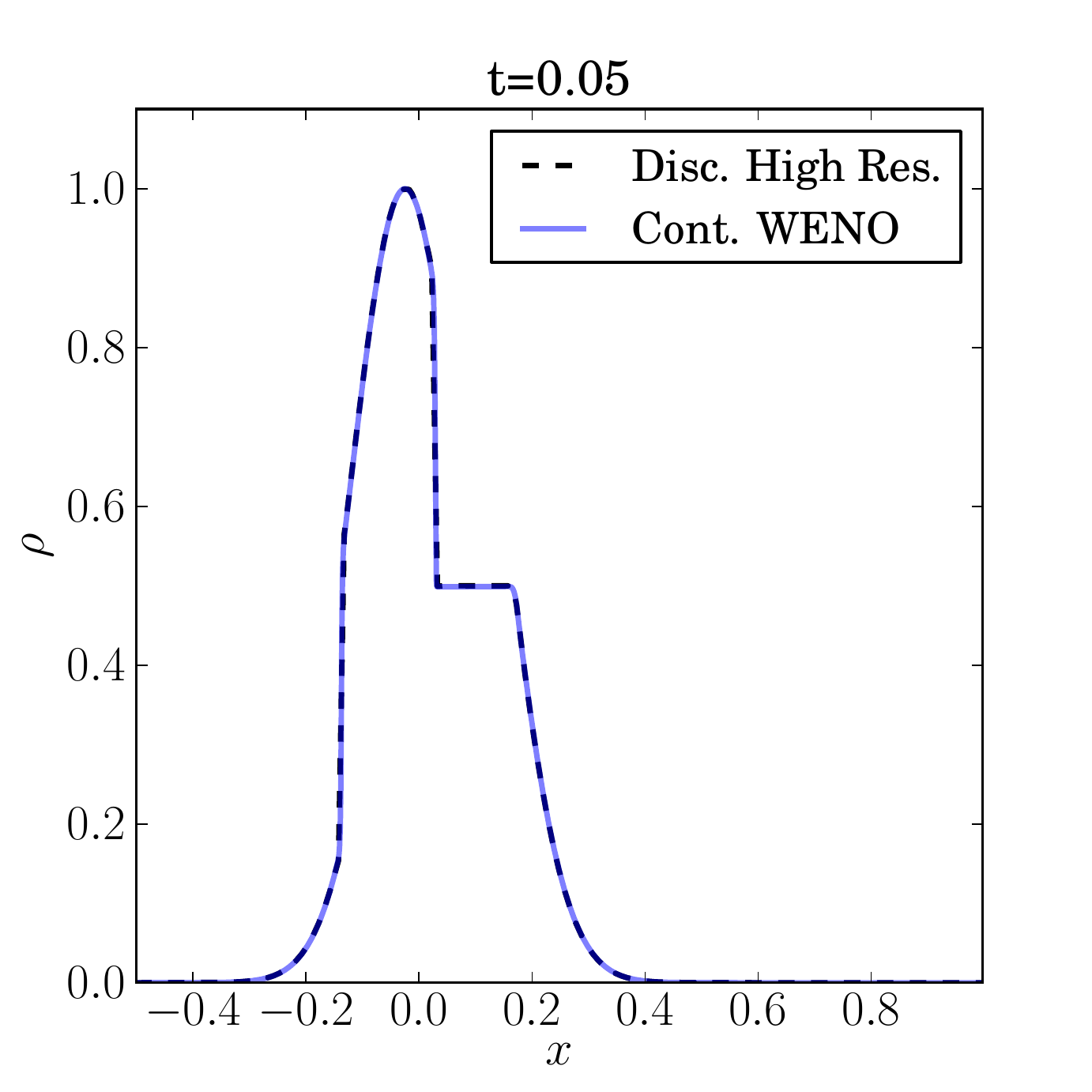}}
    \\
    \subfigure[]{\includegraphics[width=0.30\textwidth]{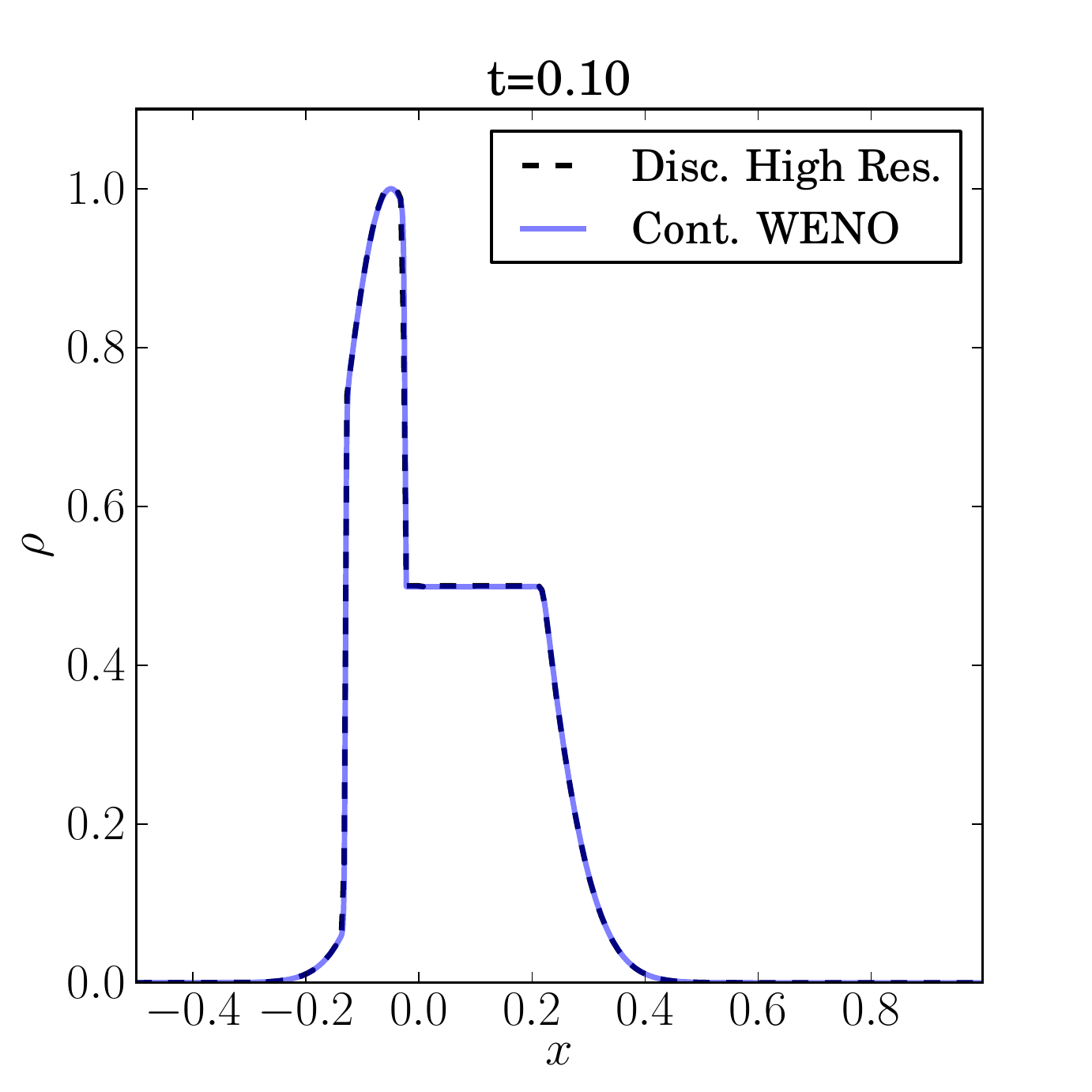}}
    \quad
    \subfigure[]{\includegraphics[width=0.30\textwidth]{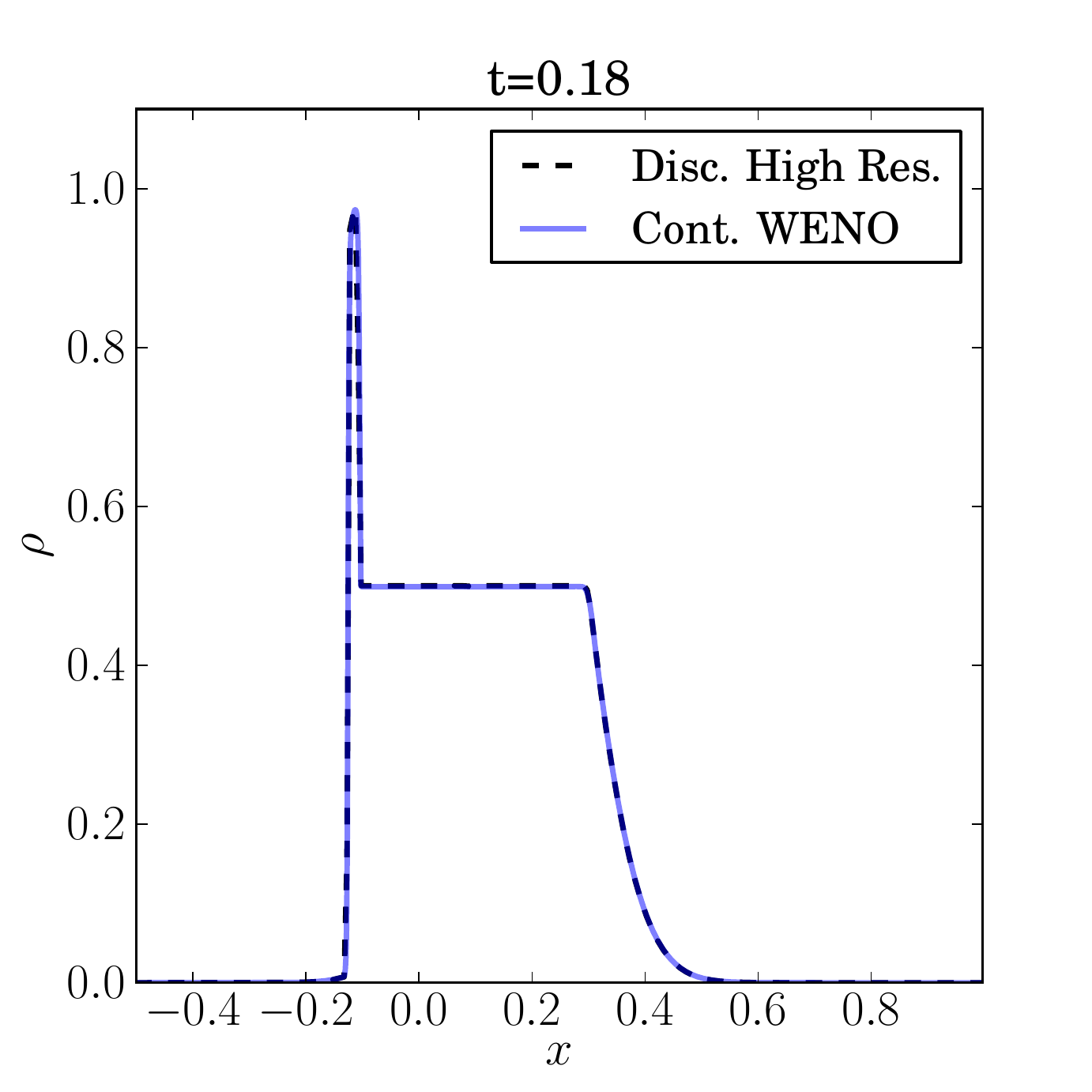}}
    \quad
    \subfigure[]{\includegraphics[width=0.30\textwidth]{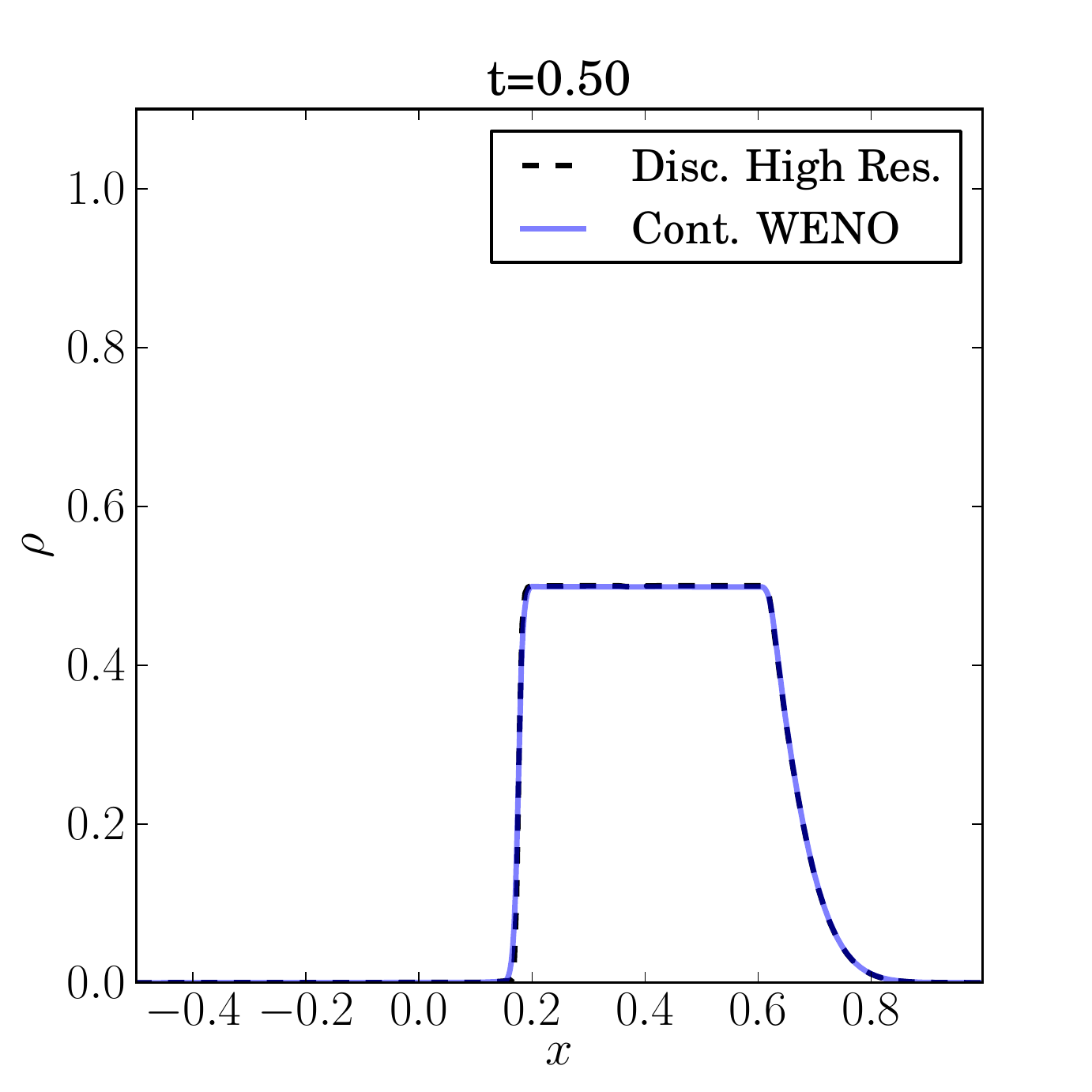}}
    \caption{Time evolution of the solution with smooth initial
      conditions using a high resolution scheme ($\delta = 10^{-5}$,
      $\Delta x = 0.005$, $\mbox{CFL}=0.9$) and a third-order WENO
      scheme \cite{CentPack} ($\Delta x = 0.0004$, $\epsilon_w =
      10^{-3}$, $\mbox{CFL}=0.9$) where $\rho_{m} = 0.5$ and $\gamma =
      0.5$.}
    \label{fig:Gaussian_Sim}
  \end{center}
\end{figure}

For the next series of simulations, we use the smooth initial data
\bq
    \rho(x, 0) = \exp\left(-\frac{x^2}{2\sigma^2}\right), 
\eq
where $\sigma = 0.1$ and the domain is the interval $[-1,1]$ with
periodic boundary conditions.  This corresponds in the traffic flow
context to a single platoon of cars travelling on a ring road, where
the initial density peak has a maximum value of 1 and decreases smoothly
to zero on either side.  Fig.~\ref{fig:Gaussian_Sim} depicts the time
evolution of the solution computed using our high resolution Godunov
scheme.  Initially, the platoon begins to spread out and a horizontal
plateau appears on the right side at a value of $\rho=\rho_m$.  This
plateau value corresponds to the ``optimal'' traffic density that is the
maximum value of $\rho$ for which cars can still propagate at the
free-flow speed.  The plateau lengthens as a shock propagates to the
left into the upper half of the density profile, reducing the width of
the peak.  At the same time, the cars in the dense region spread out to
the right as the plateau also extends in the same direction, while the
left edge of the platoon remains essentially stationary until the dense
peak is entirely gone.  When the peak finally disappears (near time
$t=0.18$), the remaining platoon of cars propagates to the right with
constant speed and unchanged shape.  Note that the traffic density
evolves such that the area under the solution curve remains
approximately constant, since the total number of cars must be
conserved.
 
For comparison purposes, we have also performed simulations using the
WENO scheme in CentPack \cite{CentPack}, which employs a third-order
CWENO reconstruction in space \cite{KurganovLevy2000,LevyEtAl2001} and a
third-order SSP Runge-Kutta time integrator
\cite{GottliebEtAl2001}. Because this algorithm requires the flux
function to be continuous (although not smooth), we have used a
regularized version of the flux that is piecewise linear and continuous,
replacing the discontinuity by a steep line segment connecting the two
linear pieces over a narrow interval of width $2\epsilon_w$ (instead of
using the function $f_\epsilon$ in \eqref{eqn:MollifiedFlux} because
that would require an integral to be evaluated for every flux function
evaluation).  The regularized flux is shown in Fig.~\ref{fig:myflux}.

From Fig.~\ref{fig:Gaussian_Sim}, we observe that the WENO
simulation requires a substantially smaller time step and grid spacing
in order to obtain results that are comparable to our method.  In
particular, the WENO scheme requires a time step of $\Delta t =
O(10^{-6})$ and a spatial resolution of $\Delta x = O(10^{-4})$ when
$\epsilon_w = 10^{-3}$; this can be compared with our high resolution
Godunov scheme for which we used a time step of $\Delta t = O(10^{-4})$
when $\Delta x = O(10^{-3})$ and $\delta = 10^{-5}$.  This performance
difference is magnified further as $\epsilon_w$ decreases due to the
ill-conditioning of the regularized-flux problem.

\begin{figure}[!tbp]
  \begin{center}
    \subfigure[]{\includegraphics[width=0.44\textwidth]{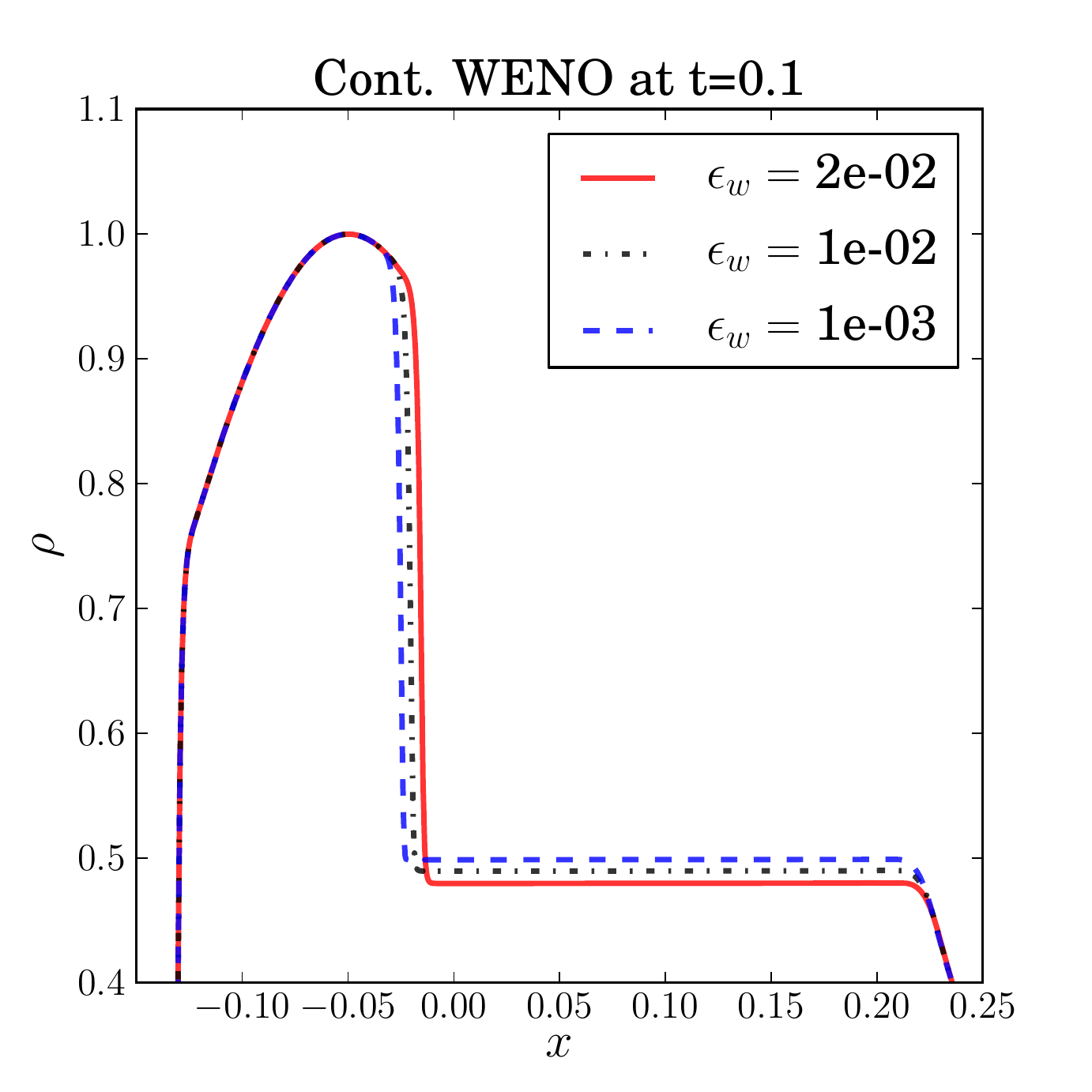} \label{fig:Gaussian_Weno_Eps}}
    \qquad 
    \subfigure[]{\includegraphics[width=0.44\textwidth]{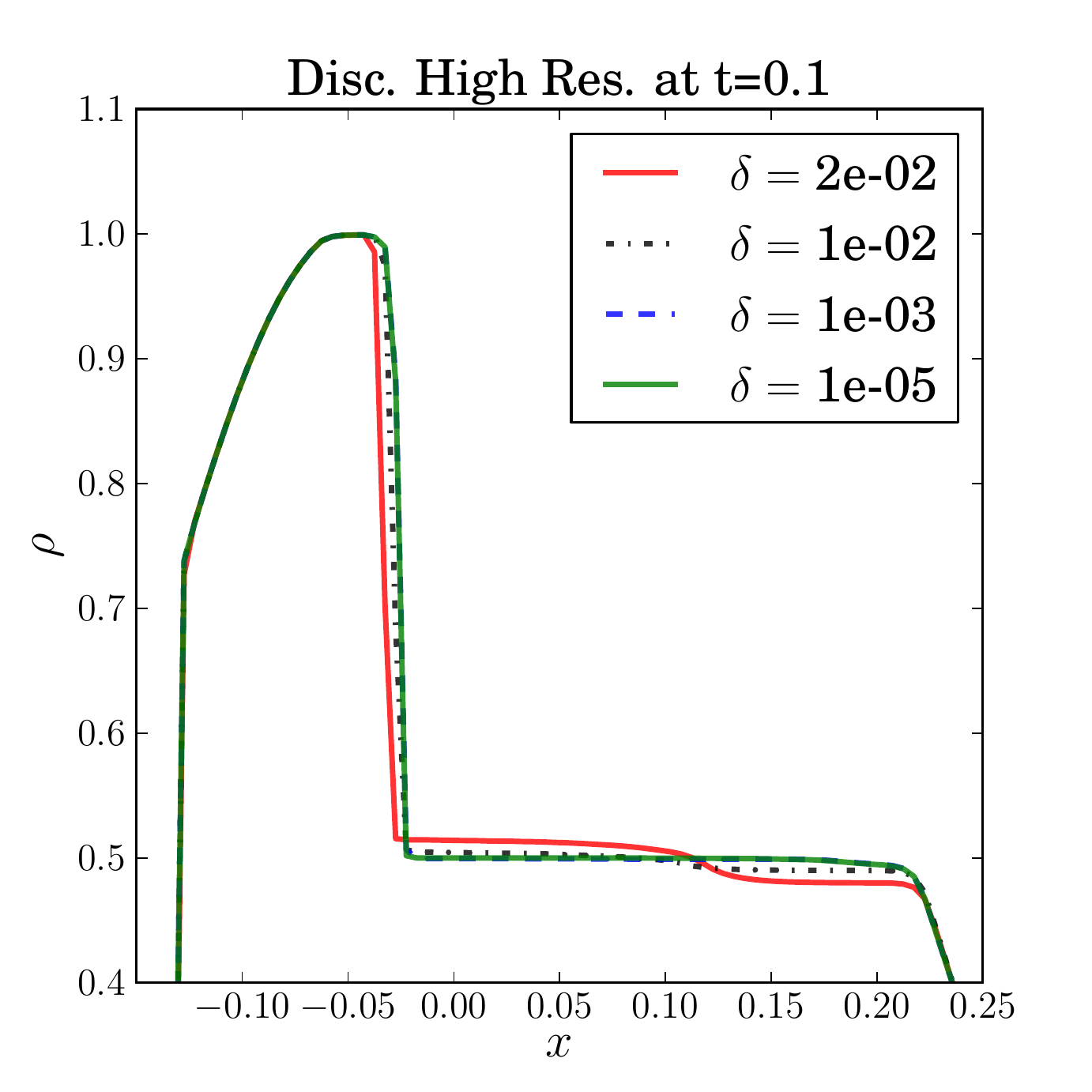} \label{fig:Gaussian_Superbee_Eps}}
    \caption{\subref{fig:Gaussian_Weno_Eps} WENO simulation for different $\epsilon_w$ for problem pictured in Fig.~\ref{fig:Gaussian_Sim}
      where $\Delta x = 0.0004$ and $\Delta t= O(10^{-6})$). 
      \subref{fig:Gaussian_Superbee_Eps} High resolution simulation for different $\delta$ for problem pictured in Fig.~\ref{fig:Gaussian_Sim}
      where $\Delta x = 0.005$ and $\Delta t= O(10^{-4})$). }
    \label{fig:Gaussian_Eps} 
  \end{center}
\end{figure}

In Fig.~\ref{fig:Gaussian_Eps}, we magnify the region containing the
density peak at time $t=0.1$ to more easily visualize the difference
between the two sets of results for different values of $\epsilon_w$ and
$\delta$. From these plots we observe that the two solutions approach
one another as both $\epsilon_w$ and $\delta$ are reduced, which
provides further evidence that our high resolution Godunov scheme
computes the correct solution.  Note that the WENO scheme does yield a
slightly sharper resolution of the shock than our method, but on the
other hand it significantly underestimates the height of the plateau
region when $\epsilon_w$ is too large.  

\begin{figure}[!tbp]
  \begin{center}
    \subfigure[]{\includegraphics[width=0.44\textwidth]{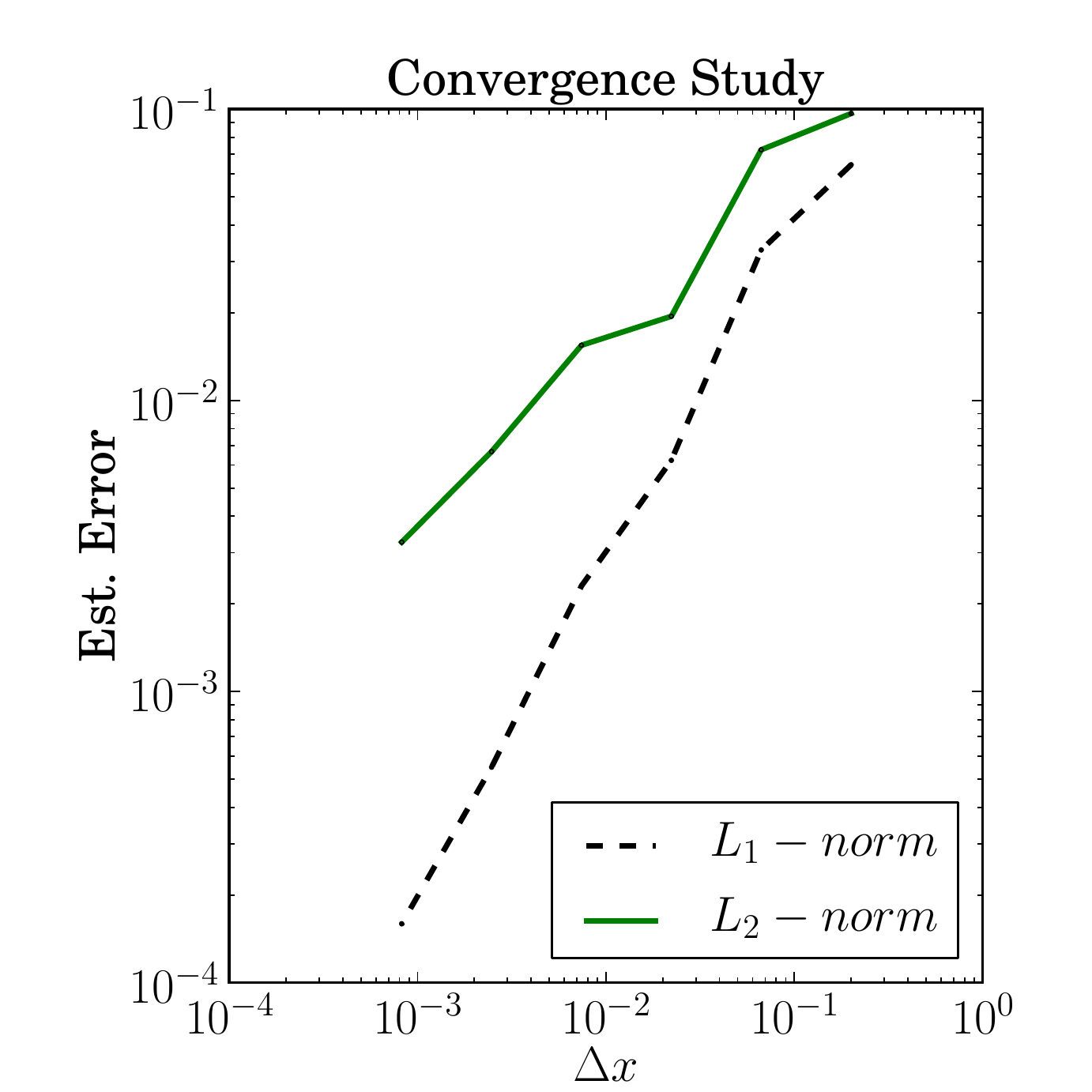} \label{fig:Gaussian_ConvStudy}}
    \qquad 
    \subfigure[]{\includegraphics[width=0.44\textwidth]{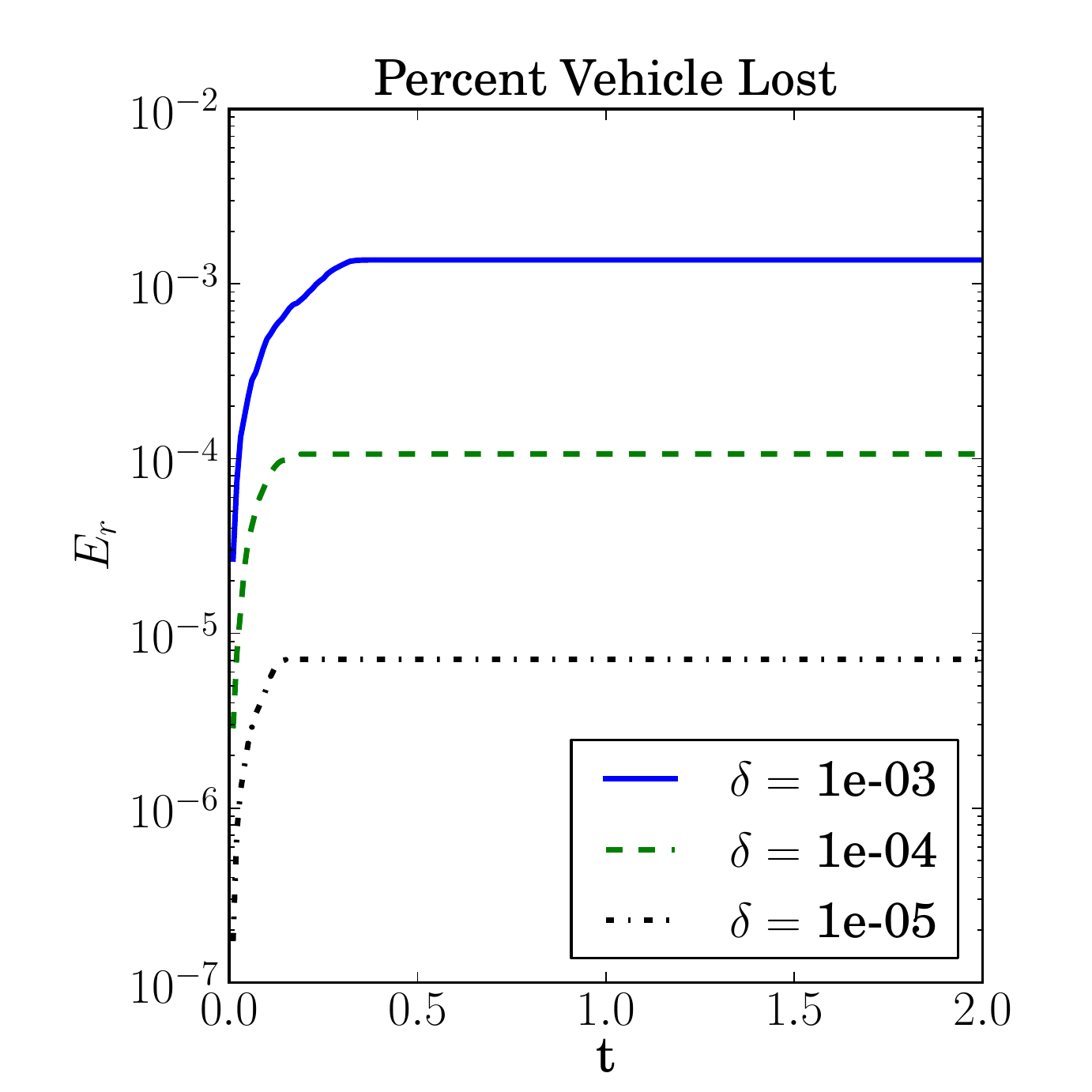} \label{fig:Gaussian_CarLost}}
    \caption{\subref{fig:Gaussian_ConvStudy} Convergence study for the
      problem with smooth initial data at $t=0.05$ using $\Delta x_0
      =0.2$ and $P=6$ levels of refinement (refer to
      Fig.~\ref{fig:Gaussian_Sim}).  \subref{fig:Gaussian_CarLost}
      Percentage of the initial vehicles lost during the simulation
      using different $\delta$ values.}
    \label{fig:Gaussian_Errors} 
  \end{center}
\end{figure}

Next, we estimate the error in our high resolution Godunov scheme by
comparing the computed solutions on a sequence of successively finer
grids with $\Delta x_p = \Delta x_0/3^p$ for $p = 0,1,\ldots,P$, with
$P=6$ levels of refinement.  The finest grid solution with $\Delta x =
\Delta x_P$ is treated as the ``exact'' solution for the purposes of
this convergence study.  The errors at $t=0.05$ shown in
Fig.~\ref{fig:Gaussian_Errors}\subref{fig:Gaussian_ConvStudy} exhibit
convergence rates of approximately 1.125 in the L1 norm and 0.632 in the
L2 norm which are consistent with the results from
Section~\ref{sec:results:riemann}.  We note that even though the initial
data are smooth, our high resolution Godunov scheme does not obtain
second order accuracy because of the shock that appears immediately on
the right side of the plateau.

Because convergence rates provide only a rough measure of solution
error, we can gain additional insight into the accuracy of the method by
measuring \emph{conservation error}, which is expressed in terms of the
variation $V$ in the total number of vehicles via
\bqs
    E_r = \frac{V^n - V^0}{V^0}
    \quad \text{where} \quad
    V^n = \Delta x \sum_{j=1}^{N} Q_j^n. 
\eqs
Since we are using periodic boundary conditions and have no external
sources or sinks, $V^n$ should remain approximately constant for all
$n$.  Indeed,
Fig.~\ref{fig:Gaussian_Errors}\subref{fig:Gaussian_CarLost} shows that
the numerical scheme accumulates only a small conservation error over
time, and that the rate of vehicles lost can be controlled by reducing
$\delta$.

\subsection{Smooth Initial Data, With Continuous Flux Comparison}
\label{sec:results:continuous}

In this final set of simulations, we compare solutions of the
conservation law \eqref{eqn:DiscConsLaw} with the discontinuous
piecewise linear flux \eqref{eqn:DiscFlux} and the continuous piecewise
linear flux
\bq 
  f_c(\rho) =  \left\{
    \begin{array}{ll}
      \displaystyle \rho, & \text{if}~0 \leqslant \rho < \rho_{m}, \\
      \displaystyle \frac{\rho_m}{1-\rho_m}(1-\rho), & \text{if}~\rho_{m} \leqslant
      \rho \leqslant 1.
    \end{array}
  \right.
  \label{eqn:ContFlux}
\eq
Both fluxes are illustrated in Fig.~\ref{fig:myflux}.

We begin by emphasizing that these two seemingly very different flux
functions can still give rise to similar solutions to the Riemann
problem with suitably chosen piecewise constant initial data.  For
Cases~A and~C from Section \ref{sec:RiemannProblem}, the convex hulls
and the corresponding solutions are identical for the two
fluxes. However, the convex hulls are different in Case~B, where the
discontinuous flux gives rise to the compound wave illustrated in
Fig.~\ref{fig:Case2_RiemannSolution} while the continuous flux generates
a single shock (analogous to the solution in
Fig.~\ref{fig:Case3_RiemannSolution}).  Furthermore, the continuous flux
\eqref{eqn:ContFlux} does not give rise to any zero waves; instead, when
either $\rho_l=\rho_m$ or $\rho_r=\rho_m$, a single contact line is
produced that satisfies the Rankine-Hugoniot condition.

A clear illustration of the difference between these two fluxes is
provided by comparing simulations on the periodic domain $x \in [-1,1]$
for smooth initial data
\bq
    \rho(x, 0) =\frac{1}{2} \exp\left(-\frac{x^2}{2\sigma^2}\right) + \frac{2}{5}, 
\eq
where $\sigma = 0.1$.  This corresponds to the situation where there is
a Gaussian-shaped congestion in the middle of free-flow traffic on a
ring road.

\begin{figure}[!tbp]
  \begin{center}
    \subfigure[]{\includegraphics[width=0.30\textwidth]{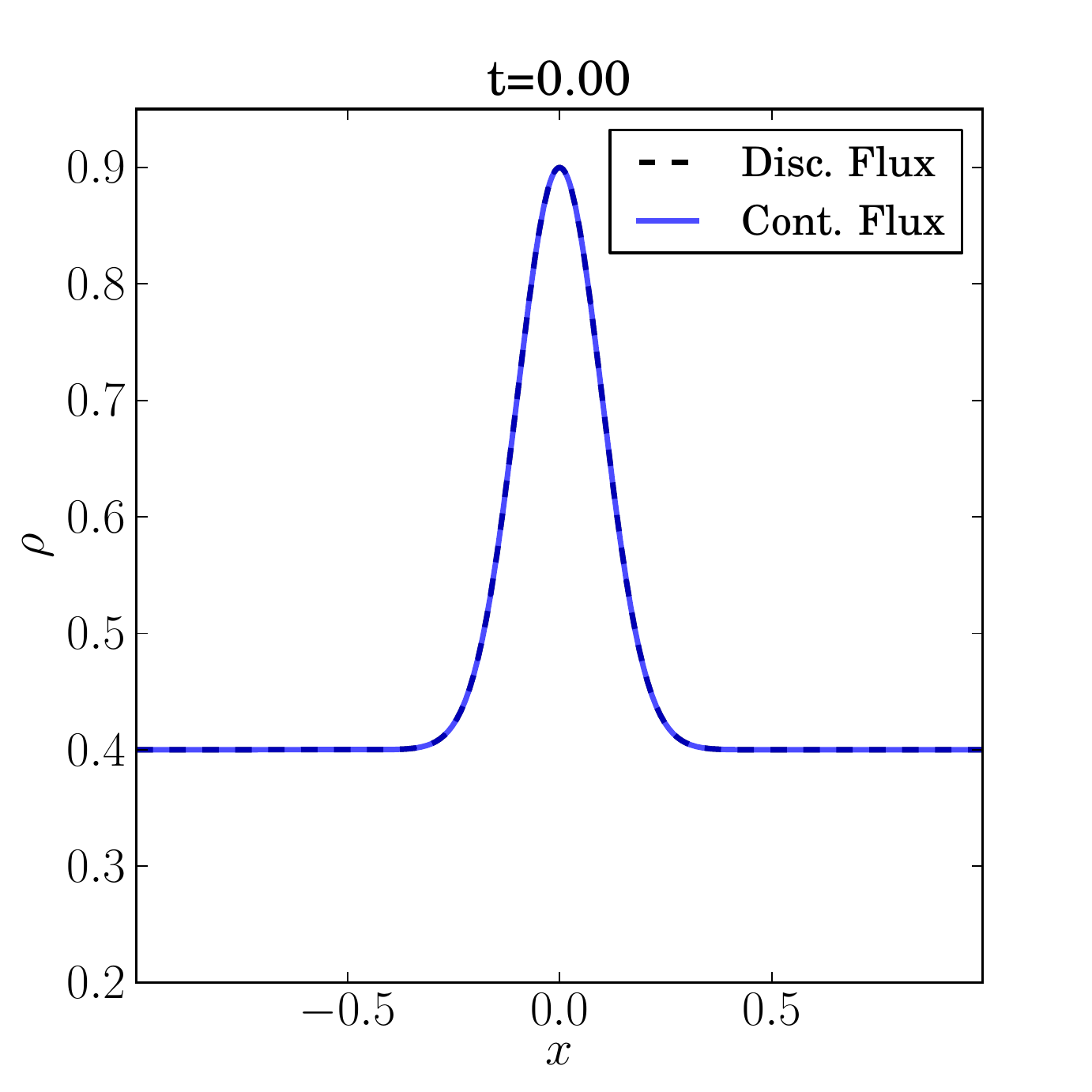}}
    \quad
    \subfigure[]{\includegraphics[width=0.30\textwidth]{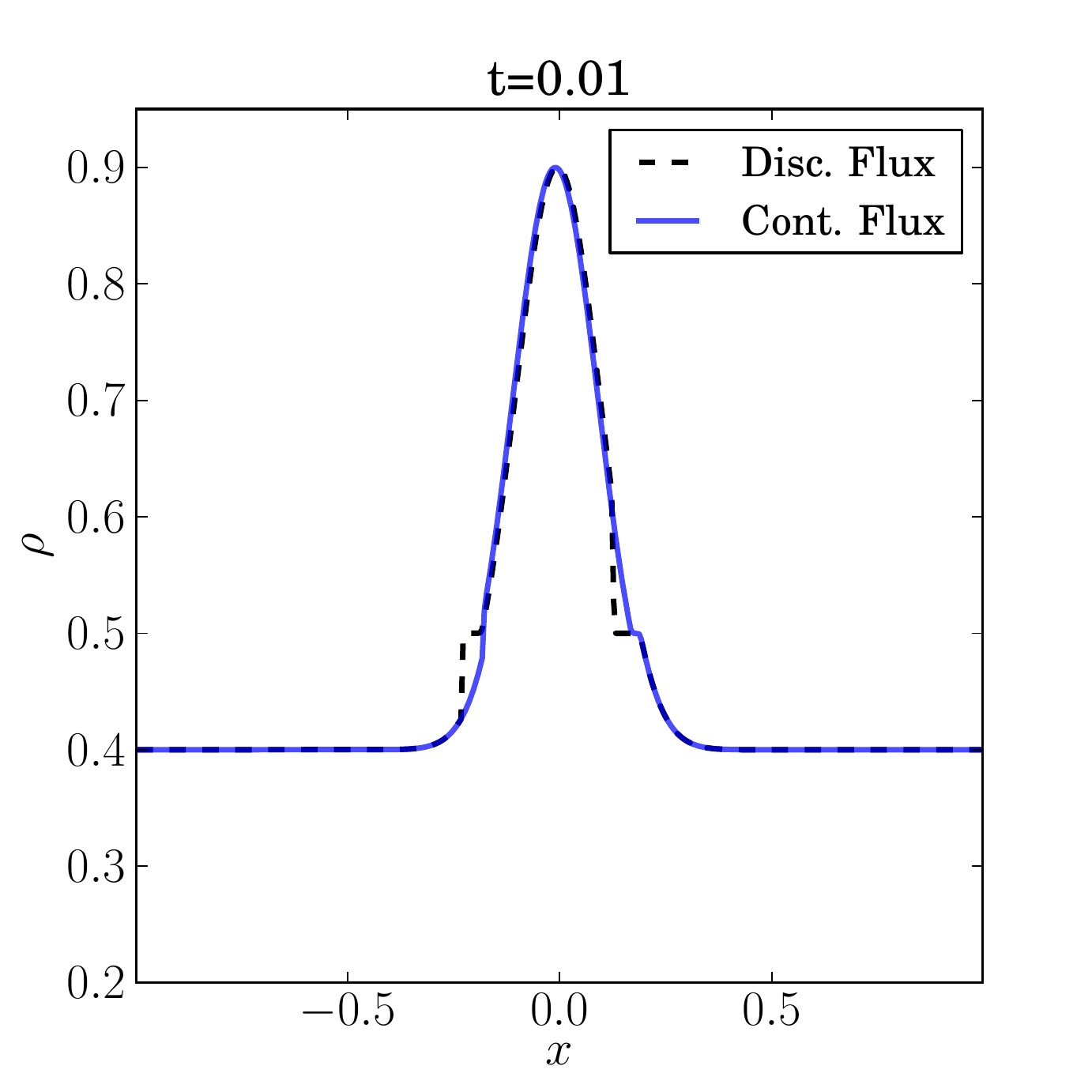}}
    \quad
    \subfigure[]{\includegraphics[width=0.30\textwidth]{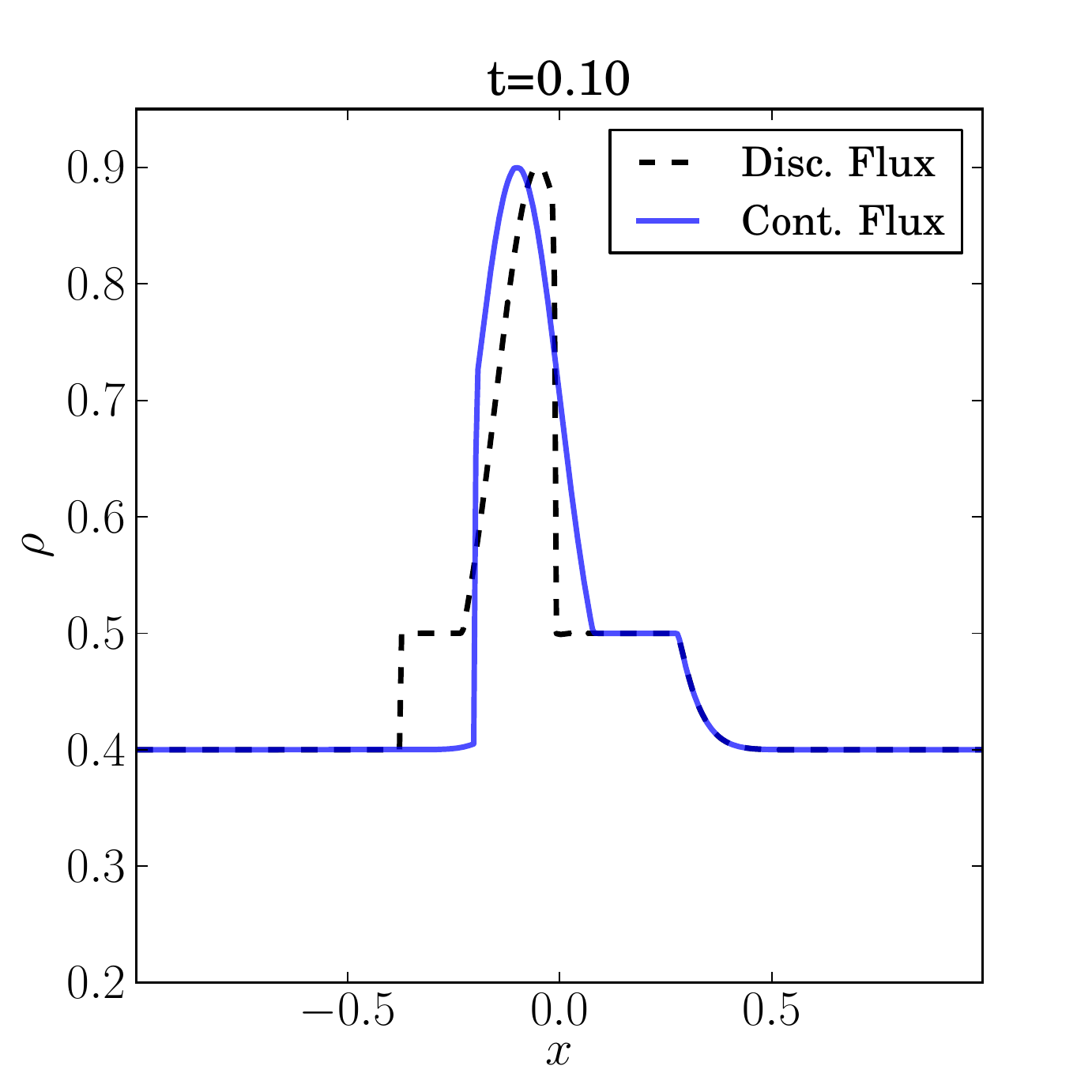}}
    \\
    \subfigure[]{\includegraphics[width=0.30\textwidth]{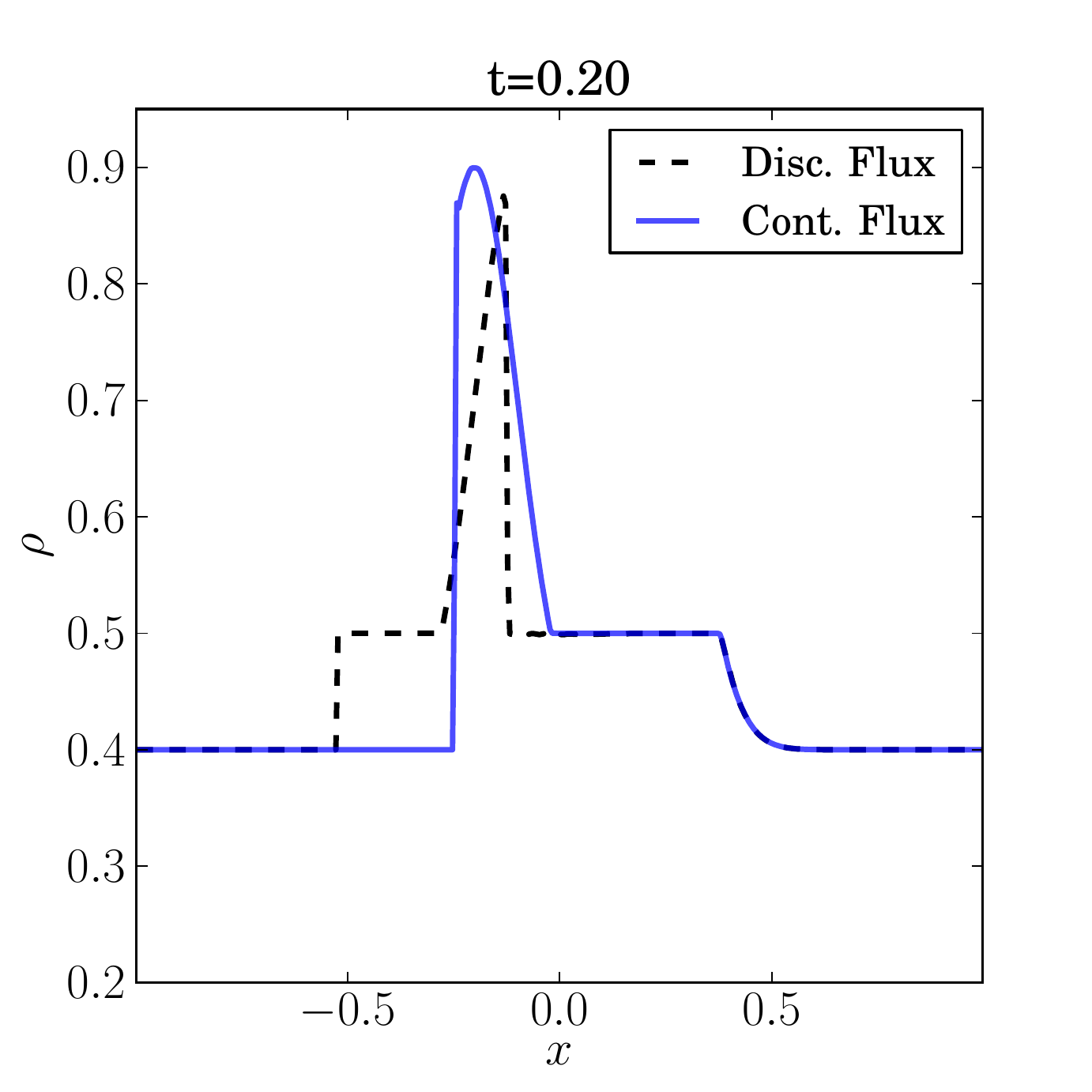}}
    \quad
    \subfigure[]{\includegraphics[width=0.30\textwidth]{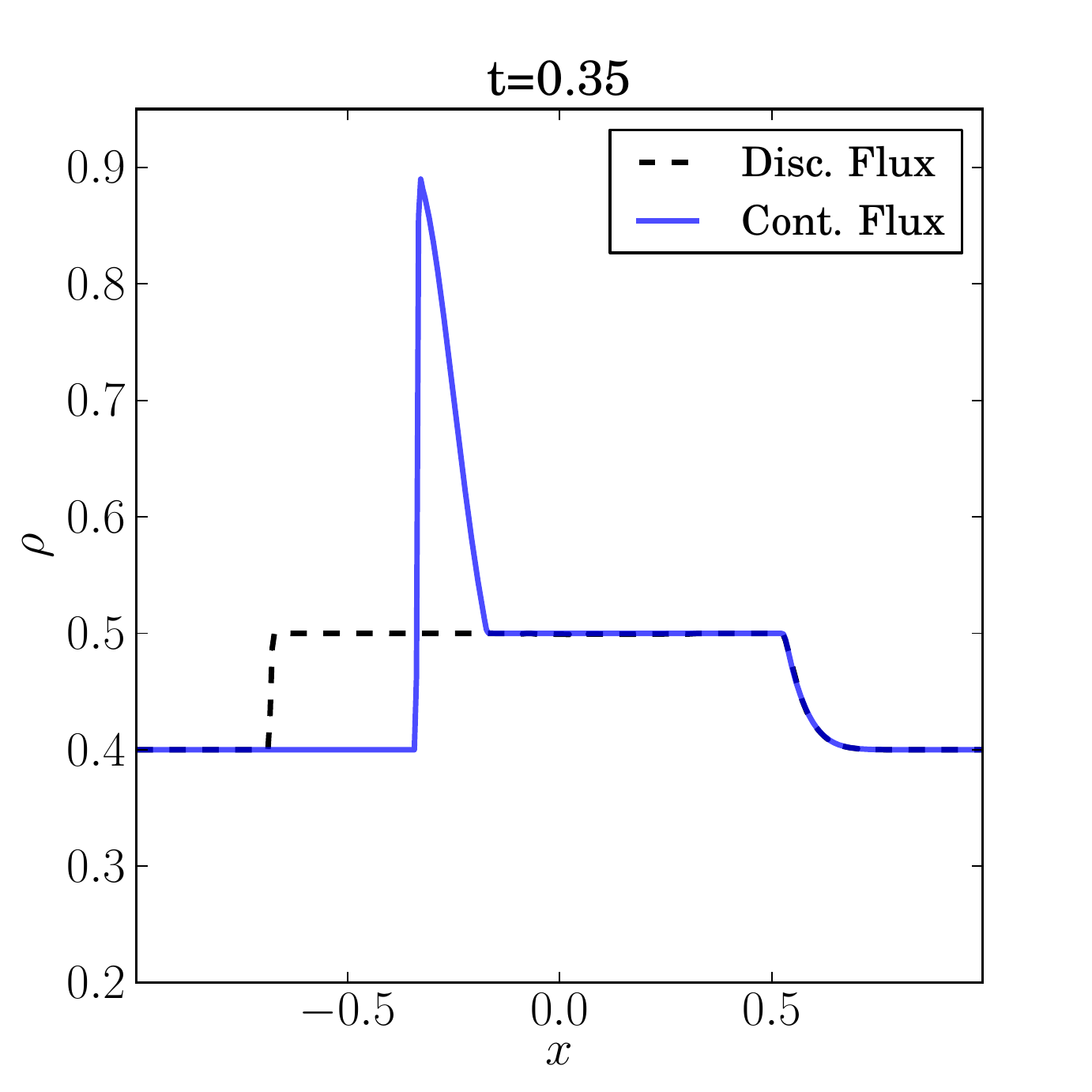}}
    \quad
    \subfigure[]{\includegraphics[width=0.30\textwidth]{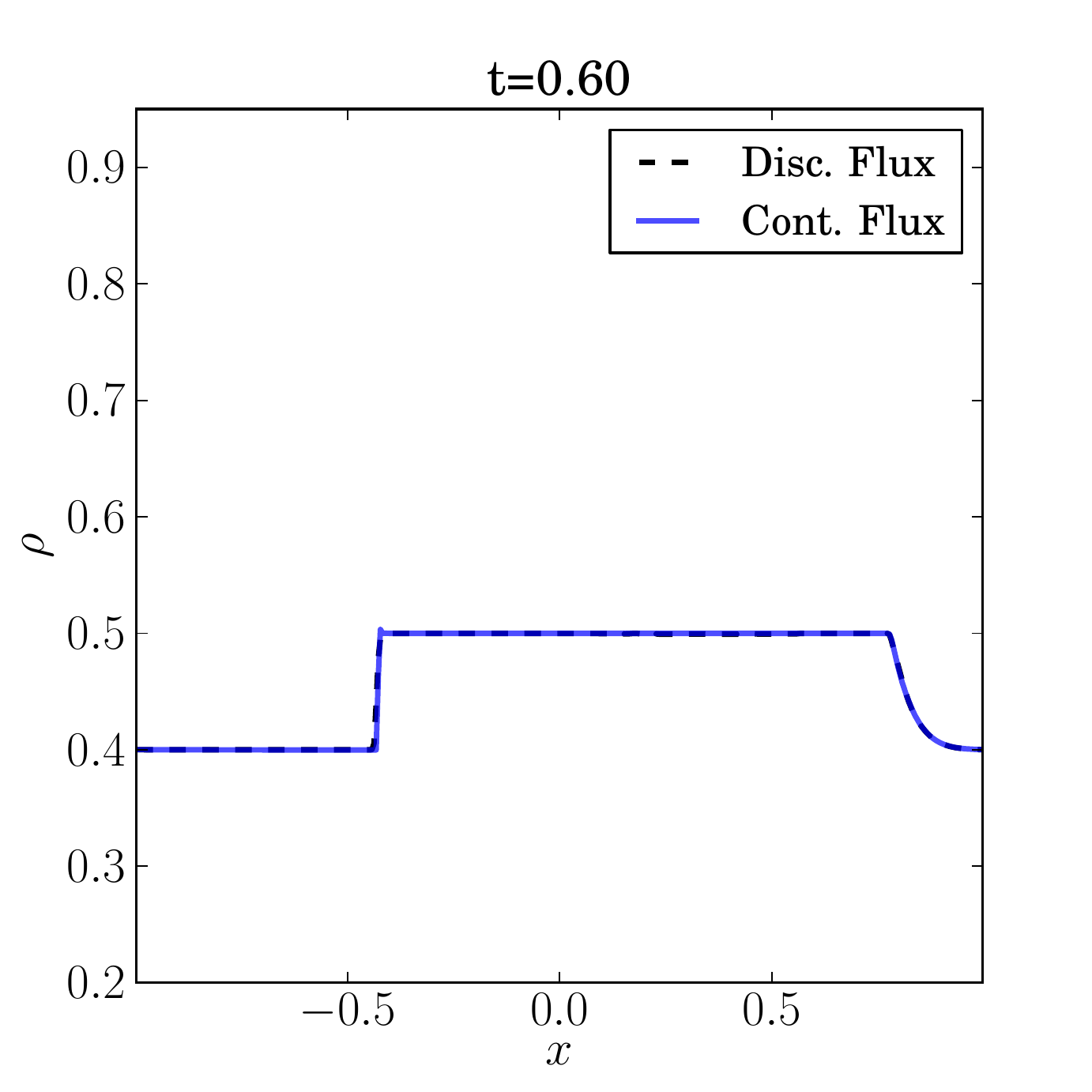}}
    \caption{Time evolution of the solution with smooth initial
      conditions using a high resolution scheme ($\delta = 10^{-5}$,
      $\Delta x = 0.005$, $\mbox{CFL}=0.9$) for both the smooth and discontinuous flux.}
    \label{fig:Gaussian2_Sim}
  \end{center}
\end{figure}

Fig.~\ref{fig:Gaussian2_Sim} depicts the time evolution of the solution
computed using a high resolution Godunov scheme for the continuous and
discontinuous fluxes.  Note that the numerical solution for the
discontinuous flux requires application of the look-ahead procedure
discussed in Section~\ref{sec:NumericalScheme} while the continuous flux
requires only the calculation of interactions between adjacent cells.
Initially, the congested region begins to spread out and a horizontal
plateau appears on the right side at a value of $\rho=\rho_m$.  On the
left edge of this plateau, a shock appears for the discontinuous flux in
comparison with a gradual continuous variation in density for the
continuous flux.    

The solutions are more drastically different when comparing the dynamics
on the left edge of the congested region.  For the discontinuous flux, a
compound wave forms to the left of the peak, resulting in the formation
of a second horizontal plateau at the maximum congested flow speed
$\rho=\rho_m$. Therefore, drivers on the far left will enter into a
optimal congested flow before being hit by the left-moving congestion
wave.  The continuous flux, on the other hand, produces only a single
left-moving congestion wave.

We conclude therefore that if the fundamental diagram is discontinuous,
we expect to observe compound congestion waves. The leading congestion
wave acts to push vehicles into some form of optimal congested flow, possibly a
synchronized traffic state. This wave is followed up by a slower
congestion wave. Note that we only expect these compound congestion waves
when the upstream free flow traffic has sufficiently high density. For
example, in Fig.~\ref{fig:Gaussian_Sim}, we observed only a single
congestion shock because the traffic density to the left is too small to
sustain this compound congestion wave.

%% file: Sec6.tex
\section{Conclusions}

In this paper, we derived analytical solutions to the Riemann problem
for a hyperbolic conservation law with a piecewise linear flux function
having a discontinuity at the point $\rho=\rho_m$.  In the case when the
left and right states in the Riemann initial data lie on either side of
the discontinuity, the solution consists of a compound wave made up of a
shock and contact line connected by a constant intermediate state at
$\rho_m$.  In the special case when either the left or right state
equals $\rho_m$, the Riemann problem gives rise to a {zero rarefaction
  wave} that propagates with infinite speed.  Even though the strength
of this wave is zero, it nonetheless has a significant impact on the
solution structure as it interacts with other elementary waves.  

Our analytical results were validated using a high resolution
Godunov-type scheme, based on our exact Riemann solver and implemented
using the wave propagation formalism of LeVeque~\cite{LeVequeRedBook}.
This approach builds the effect of {zero waves} directly into the
algorithm in a way that avoids the overly stringent CFL time step
constraint that might otherwise derive from the infinite speed of
propagation of zero waves.  We demonstrate the accuracy and efficiency
of our method using several test problems, and include comparisons with
higher-order WENO simulations.

We conclude with a brief discussion of several possible avenues for
future research.  First, a detailed convergence analysis of the
algorithm would help to identify non-smooth error components that arise
from the discontinuity in the flux, as well as the look-ahead procedure
required to determine shock speeds in the case when $\rho_m=\rho_l$ or
$\rho_r$.  Second, we would like to consider other nonlinear forms of
the piecewise discontinuous flux function since the linearity assumption
in this paper simplifies our analytical solutions considerably in that
elementary waves arising from local Riemann problems consist of shocks
and contact lines only.  For example, it would be interesting to perform
a detailed comparison with the results of Lu et al.~\cite{Lu2009} who
considered a discontinuous, piecewise quadratic flux.  Thirdly, we
mention some preliminary computations of traffic flow using a cellular
automaton model \cite{Wiens2011} that give rise to an apparent
discontinuity in the fundamental diagram.  This connection between
cellular automaton models (that only specify rules governing individual
driver behaviour) and kinematic wave models (in which the two-capacity
effect is incorporated explicitly via the flux function) merits further
study.  Finally, the situation where the flux is a discontinuous
function of the spatial variable has been analyzed much more
extensively (see the journal issue introduced by the article
\cite{BurgerKarlsen2008}, and references therein).  We would like to
draw deeper connections between this work and the problem where the
discontinuity appears in the density.